\documentclass[12pt]{amsart}
\usepackage{a4wide}
\usepackage{amssymb}
\usepackage{amsthm}
\usepackage{amsmath}
\usepackage{amscd}
\usepackage{verbatim}
\usepackage{bm}
\usepackage{enumerate}
\usepackage[all]{xy}

\numberwithin{equation}{section}

\theoremstyle{plain}
\newtheorem{thm}{Theorem}[section]
\newtheorem{prop}[thm]{Proposition}
\newtheorem{cor}[thm]{Corollary}
\newtheorem{lem}[thm]{Lemma}

\newtheorem{hypo}[thm]{Hypothesis}

\theoremstyle{remark}
\newtheorem{rmk}[thm]{Remark}

\theoremstyle{definition}
\newtheorem{defn}[thm]{Definition}
\newtheorem{ex}[thm]{Example}

\newcommand{\kzxz}[4]{\left(\begin{smallmatrix} #1 & #2 \\ #3 & #4\end{
smallmatrix}\right) }

\begin{document}

\title[Special cycles on toroidal compactifications]{Special cycles on toroidal compactifications of orthogonal Shimura varieties}

\author[Jan H.~Bruinier]{Jan Hendrik Bruinier}
\author[Shaul Zemel]{Shaul Zemel}

\address{Fachbereich Mathematik,
Technische Universit\"at Darmstadt, Schlossgartenstrasse 7, D--64289
Darmstadt, Germany}
\email{bruinier@mathematik.tu-darmstadt.de}

\address{Einstein Institute of Mathematics, the Hebrew University of Jerusalem, Edmund Safra Campus, Jerusalem 91904, Israel}
\email{ zemels@math.huji.ac.il}

\thanks{The first author is partially supported by DFG grant BR-2163/4-2 and the LOEWE research unit USAG}


\date{\today}

\begin{abstract}
We determine the behavior of automorphic Green functions along the boundary components of toroidal compactifications of orthogonal Shimura varieties. We use this analysis to define boundary components of special divisors and prove that the generating series of the resulting special divisors on a toroidal compactification is modular.
\end{abstract}

\maketitle

\section{Introduction}

Orthogonal Shimura varieties that arise from even lattices of signature $(n,2)$ carry a natural family of divisors that are called Heegner divisors, rational quadratic divisors, or special divisors in the literature. In his seminal paper \cite{[Bo1]}, Borcherds constructed meromorphic modular forms on these varieties, and then used them in the sequel \cite{[Bo2]} to establish the modularity of the generating series of these divisors in the first Chow group of such a Shimura variety. This provided a new proof of the Gross--Kohnen--Zagier theorem from \cite{[GKZ]} about the modularity of Heegner divisors on modular curves and a generalization to orthogonal Shimura varieties of arbitrary dimension and level. Taking Chern classes one obtains as a corollary also the modularity of the generating series of the cohomology classes of special divisors, such as the celebrated Hirzebruch-Zagier theorem for Hilbert modular surfaces \cite{[HZ]} and particular cases of results of Kudla-Millson on cohomology classes of special cycles on locally symmetric spaces associated with orthogonal groups \cite{KM}.

The modularity of generating series of special cycles plays a central role in the Kudla program, see, e.g.,~\cite{[K1]}, \cite{[K2]}. Once one knows modularity, the generating series can be used to construct arithmetic analogues of the theta correspondence, relating Siegel modular forms to classes of special cycles in (arithmetic) Chow groups. For applications to intersection and height pairings one is interested in special cycles on smooth \emph{toroidal compactifications} and corresponding modularity results, cf.~Problem~3 of \cite{[K2]}. The present paper is a contribution in this direction.

In order to present the results in more detail, we need some notation. Let $L$ be an even lattice of signature $(n,2)$, and let $\Gamma_{L}$ be the stable special orthogonal group of $L$, see \eqref{GammaLdef}. We write $X_{L}=\Gamma_{L}\backslash\mathcal{D}$ for the associated connected complex Shimura variety. For every element $\mu$ in the discriminant group $\Delta_{L}=L^{*}/L$ of $L$ and every positive number $m\in\frac{\mu^{2}}{2}+\mathbb{Z}$ there is a corresponding special divisor $Z(m,\mu)$ on $X_{L}$. If $\rho_{L}$ denotes the Weil representation associated with $L$, then Theorem~4.5 of \cite{[Bo2]} states the following.
\begin{thm} \label{modopen}
Let $Z(0,0)$ be any divisor on $X_{L}$ representing the line bundle of modular forms of weight $-1/2$, and denote the image of $Z(m,\mu)$ in $\operatorname{CH}^{1}(X_{L})\otimes\mathbb{Q}$ by $[Z(m,\mu)]$ for any $m$ and $\mu$. Then the formal power series \[[Z(0,0)]\mathfrak{e}_{0}+\sum_{\mu\in\Delta_{L}}\sum_{0<m\in\frac{\mu^{2}}{2}+\mathbb{Z}}[Z(m,\mu)] \cdot q^{m}\mathfrak{e}_{\mu}\] is a modular form of weight $1+n/2$ and representation $\rho_{L}$ with coefficients in $\operatorname{CH}^{1}(X_{L})\otimes\mathbb{Q}$.
\end{thm}
Note that Theorem \ref{modopen} involves the representation $\rho_{L}$, and not the complex conjugate representation $\overline{\rho}_{L}$ as in Theorem~4.5 of \cite{[Bo2]}, since we work here with the opposite signature convention. The idea of the proof is that given a weakly holomorphic modular form $f$ of weight $1-n/2$ and representation $\overline{\rho}_{L}$ with integral principal part, there exists, according to Theorem~13.3 of \cite{[Bo1]}, a Borcherds product $\Psi_{f}$ on $X_{L}$ that gives rise to a particular relation among the special divisors in $\operatorname{CH}^{1}(X_{L})$. The result can then be deduced by means of an incarnation of Serre duality, which appears as Theorem~3.1 of \cite{[Bo2]}, and which is stated here as Proposition~\ref{Serredual}.

\smallskip

The variety $X_{L}$ is typically non-compact. Following work of Satake \cite{[Sa]}, Baily and Borel constructed in \cite{[BB]} a canonical minimal compactification $X_{L}^{\mathrm{BB}}$ of $X_{L}$. When $n$ is larger than the Witt rank of $L$, an appropriate equivalent of Theorem \ref{modopen} holds also for the closures of the divisors $Z(m,\mu)$ on $X_{L}^{\mathrm{BB}}$. Indeed, meromorphic modular forms with respect to $\Gamma_{L}$, and in particular Borcherds products, are defined on $X_{L}^{\mathrm{BB}}$, and the complement of $X_{L}$ in $X_{L}^{\mathrm{BB}}$ has codimension larger than 1 by the assumption on the Witt rank. Hence the divisor of any Borcherds product yields a relation among the special divisors, and the result follows in the same manner.

The Baily-Borel compactification $X_{L}^{\mathrm{BB}}$ is typically highly singular at the boundary. In order to work on a smooth compact variety, we shall consider \emph{toroidal compactifications} as constructed in \cite{[AMRT]}. These are described explicitly in the orthogonal setting for instance in \cite{[Fi]}, and we present the results here as well in a manner that is more suitable for our purposes. However, the complement of $X_{L}$ inside a toroidal compactification $X_{L}^{\mathrm{tor}}$ is a divisor, and it turns out that Theorem \ref{modopen} does not extend to $X_{L}^{\mathrm{tor}}$ if one works with the closures of the special divisors on $X_{L}$. Instead, one has to add certain boundary divisors with appropriate multiplicities, a phenomenon which was already observed for Hilbert modular surfaces in \cite{[HZ]} and for unitary Shimura varieties in \cite{[BHY]} and \cite{BHKRY}. Here we define special divisors $Z^{\mathrm{tor}}(m,\mu)$ on $X_{L}^{\mathrm{tor}}$ which include suitable linear combinations of the boundary divisors, and we prove an analogue of Theorem \ref{modopen}.

The idea of our approach is the following. For every $\mu\in\Delta_{L}$ and every positive number $m\in\frac{\mu^{2}}{2}+\mathbb{Z}$ there is a harmonic Mass form $F_{m,\mu}$ of weight $1-\frac{n}{2}$ with representation $\overline{\rho}_{L}$ whose principal part is given by $q^{-m}(\mathfrak{e}_{\mu}+\mathfrak{e}_{-\mu})$. Its regularized theta lift against the Siegel theta function $\Theta_{L}$ of $L$ determines an automorphic Green function $\Phi_{m,\mu}^{L}$ on $X_{L}$ with a logarithmic singularity along $Z(m,\mu)$, see \cite{[Br1]} and \cite{[BFu]} for more details. We analyze the behavior of this Green function towards each boundary component of $X_{L}^{\mathrm{tor}}$. It turns out that besides a term with $\log$-$\log$-growth, the Green function has logarithmic growth along certain boundary divisors. Here the corresponding multiplicities are determined by regularized theta lifts associated to lattices of signature $(n-2,0)$ and $(n-1,1)$ which are obtained as quotients of $L$.

Every primitive rank 2 isotropic sublattice $J \subseteq L$ determines a 1-dimensional Baily-Borel boundary component and also a canonical boundary divisor $B_{J}$ of $X_{L}^{\mathrm{tor}}$. On the other hand, the lattice $D=(J^{\perp} \cap L)/J$ is positive definite of signature $(n-2,0)$. There is a linear map taking $F_{m,\mu}$ to a harmonic Maass form $F_{m,\mu}^{D}$ with representation $\overline{\rho}_{D}$. We show in Proposition~\ref{PhinearJ} that $\Phi_{m,\mu}^{L}$ has logarithmic growth along the divisor $\operatorname{mult}_{J}(m,\mu) \cdot B_{J}$, where the multiplicity $\operatorname{mult}_{J}(m,\mu)$ is determined by the regularized theta lift of $F_{m,\mu}^{D}$ against the theta function $\Theta_{D}$ of the lattice $D$, which in turn is given by $\frac{2m}{n-2}$ times the representation number of $m$ by the coset of $D$ arising from $\mu$.

Every primitive rank 1 isotropic sublattice $I \subseteq L$ determines a 0-dimensional Baily-Borel boundary component, as well as a Lorentzian lattice $K=(I^{\perp} \cap L)/I$ of signature $(n-1,1)$. In contrast to the canonical boundary divisor $B_{J}$ from the previous paragraph, the inverse image in $X_{L}^{\mathrm{tor}}$ depends on the choice of an admissible cone decomposition of the rational closure of a cone $C$ of negative norm vectors in $K_{\mathbb{R}}$. The lattice $K$ and an oriented primitive negative norm vector $\omega \in K \cap C$ that spans an inner ray in the associated cone decomposition determine a boundary divisor $B_{I,\omega}$ of $X_{L}^{\mathrm{tor}}$. There is a linear map taking $F_{m,\mu}$ to a harmonic Maass form $F_{m,\mu}^{K}$ with representation $\overline{\rho}_{K}$. We show in Proposition~\ref{Phinearomega} (or, more explicitly, in the special case presented in \eqref{expPhiomega}) that $\Phi_{m,\mu}^{L}$ has logarithmic growth along the divisor $\operatorname{mult}_{I,\omega}(m,\mu) \cdot B_{I,\omega}$, where the multiplicity $\operatorname{mult}_{I,\omega}(m,\mu)$ is given by the regularized theta lift of $F_{m,\mu}^{K}$ against the Siegel theta function $\Theta_{K}$ of $K$, evaluated at the special point associated with $\omega$ in the Grassmannian of $K$.

We define the special divisor $Z^{\mathrm{tor}}(m,\mu)\in\operatorname{Div}(X_{L}^{\mathrm{tor}})\otimes\mathbb{R}$ by the formula
\begin{equation} \label{ZtorInt}
Z^{\mathrm{tor}}(m,\mu)=Z(m,\mu)+\sum_{J \subseteq L}\operatorname{mult}_{J}(m,\mu) \cdot B_{J}+\sum_{I \subseteq L}\sum_{\mathbb{R}_{+}\omega}\operatorname{mult}_{I,\omega}(m,\mu) \cdot B_{I,\omega}.
\end{equation}
Here $J$ (resp.~$I$) runs over a set of representatives of rank 2 (resp.~rank 1) primitive isotropic sublattices of $L$ modulo $\Gamma_{L}$, and given $I$, the index $\mathbb{R}_{+}\omega$ runs over representatives for the inner rays of the chosen admissible cone decomposition. Note that while the multiplicities $\operatorname{mult}_{J}(m,\mu)$ are easily seen to be rational because of their description as representation numbers, the multiplicities $\operatorname{mult}_{I,\omega}(m,\mu)$ seem to be real in general. Only certain linear combinations of the latter quantities can be shown to be rational, cf.~Theorem~\ref{divbctor}.

Our first main result, given in Theorem~\ref{GreenPhi} below, states that $\frac{1}{2}\Phi_{m,\mu}^{L}$ is a logarithmic Green function for the divisor $Z^{\mathrm{tor}}(m,\mu)$ on $X_{L}^{\mathrm{tor}}$ with a term that is pre-$\log$-$\log$ along the boundary in the sense of \cite{[BBK]}. In particular this means that if $d^{c}$ is the complementary operator $\frac{\partial-\overline{\partial}}{4\pi i}$ then we have the current equation
\begin{equation} \label{ddceq}
dd^{c}[\tfrac{1}{2}\Phi_{m,\mu}^{L}]+\delta_{Z^{\mathrm{tor}}(m,\mu)}=[\eta^{\mathrm{tor}}(m,\mu)],
\end{equation}
where $\eta^{\mathrm{tor}}(m,\mu)$ is the sum of a smooth 2-form on $X_{L}^{\mathrm{tor}}$ and the $dd^{c}$-image of the pre-$\log$-$\log$ term.

Our second main result, Theorem~\ref{Ztormod}, describes the generating series of these special divisors.
\begin{thm} \label{main}
Assume that $n$ is larger than the Witt rank of $L$, set $Z^{\mathrm{tor}}(0,0)$ to be any divisor on $X_{L}^{\mathrm{tor}}$ that represents the line bundle of modular forms of weight $-1/2$, and write $[Z^{\mathrm{tor}}(m,\mu)]$ for the image of $Z^{\mathrm{tor}}(m,\mu)$ in $\operatorname{CH}^{1}(X_{L}^{\mathrm{tor}})\otimes\mathbb{R}$ for every pair $(m,\mu)$. Then the formal power series \[[Z^{\mathrm{tor}}(0,0)]\mathfrak{e}_{0}+\sum_{\mu\in\Delta_{L}}\sum_{0<m\in\frac{\mu^{2}}{2}+\mathbb{Z}}[Z^{\mathrm{tor}}(m,\mu)] \cdot q^{m}\mathfrak{e}_{\mu}\] is a modular form of weight $1+n/2$ and representation $\rho_{L}$ with coefficients in $\operatorname{CH}^{1}(X_{L}^{\mathrm{tor}})\otimes\mathbb{R}$.
\end{thm}

Taking advantage of our analysis of the Green functions $\Phi_{m,\mu}^{L}$, this result can now be proved in a similar way as Theorem \ref{modopen}. For every weakly holomorphic modular form $f$ of weight $1-n/2$ and representation $\overline{\rho}_{L}$ with integral principal part, we consider the associated Borcherds product $\Psi_{f}$, but now as a section of a line bundle over $X_{L}^{\mathrm{tor}}$. Since $f$ can be (essentially) uniquely written as a linear combination of the harmonic Maass forms $F_{m,\mu}$, the logarithm of the Petersson metric of $\Psi_{f}$ decomposes as a corresponding linear combination of the Green functions $\Phi_{m,\mu}^{L}$. This implies that the divisor of $\Psi_{f}$ on $X_{L}^{\mathrm{tor}}$ includes the boundary components of $X_{L}^{\mathrm{tor}}$ with multiplicities compatible with \eqref{ZtorInt}. Hence $f$ gives rise to a relation among the divisors $Z^{\mathrm{tor}}(m,\mu)$, and the theorem follows again by Serre duality. In fact, as we construct the divisors $Z^{\mathrm{tor}}(m,\mu)$ via their Green functions, our argument can also prove modularity in Chow groups arising in arithmetic settings---see Remark \ref{modCH1} below.

We remark that \cite{[P]} proves a similar result on the moduli space of polarized K3 surfaces in dimension 19, though this result does not involve the multiplicities $\operatorname{mult}_{I,\omega}(m,\mu)$ because the toroidal compactification considered there is based on fans that have no internal rays.

Our argument also serves for giving more direct proofs for some results from the literature about the expansion of the Borcherds product $\Psi_{f}$. First, given a 0-dimensional cusp and a Weyl chamber $W$, this Borcherds product involves the exponent of a linear term, and the value of this term at negative norm elements of $K$ is considered in Theorem~10.6 of \cite{[Bo1]}. In addition, Theorem~2.1 and Corollary~2.3 of \cite{[K3]} consider an expansion of $\Psi_{f}$ near a 1-dimensional cusp, involving a coordinate $q_{2}$ raised to the exponent given in Equation (4.27) there. The boundary part of the divisor of $\Psi_{f}$, using the multiplicities from \eqref{ZtorInt}, produces both of these results directly---see Remark \ref{pfBo1K3} below.

\smallskip

This paper is divided into five sections, the first being the introduction. Section \ref{TorComp} reviews quickly the general theory of toroidal compactifications, while Section \ref{CompGrass} presents the explicit description of what happens in the orthogonal case. Section \ref{Green} investigates the Green functions associated with special divisors, and explains the extensions of the Green functions to toroidal compactifications. Finally, Section \ref{Modularity} recalls some results on Borcherds products and shows how to modify the modularity theorem for special cycles on the open variety to become suitable for the toroidal case.

\smallskip

We thank Max R\"{o}ssler and Salim Tayou for useful comments.

\section{General Toroidal Compactifications \label{TorComp}}

We skim here briefly, without any proofs, the general theory of toric varieties and toroidal compactifications. More details on toric varieties can be found in \cite{[Fu]} or in \cite{[CLS]}, and the general construction of toroidal compactifications is available in \cite{[AMRT]}, \cite{[Nam]}, \cite{[Fi]}, and others. The particular parts that will be of use to us are Proposition \ref{ordTdiv}, Lemma \ref{genparstruc}, Propositions \ref{DoverXi} and \ref{subgpact}, Definition \ref{admissible}, and Theorem \ref{torcomp}.

\subsection{Toric Varieties}

There is a correspondence between algebraic tori of finite dimension and \emph{lattices}, namely free Abelian groups of finite rank. Given a lattice $\Lambda$, the corresponding torus $T_{\Lambda}$ is the one for which $T_{\Lambda}(\mathbb{C})=\Lambda_{\mathbb{C}}/\Lambda$. The dual lattice $\Lambda^{*}=\operatorname{Hom}(\Lambda,\mathbb{Z})$ is identified with $\operatorname{Hom}(T_{\Lambda},\mathbb{G}_{m})$ via
\begin{equation} \label{znudef}
\nu \mapsto z_{\nu}, \qquad z_{\nu}:\big(\zeta \in T_{\Lambda}(\mathbb{C})=\Lambda_{\mathbb{C}}/\Lambda\big)\mapsto\mathbf{e}\big((\nu,\zeta)\big),
\end{equation}
where we adopt the usual notation $\mathbf{e}(w)=e^{2\pi iw}$ for $w\in\mathbb{C}$ and the pairing $(\nu,\zeta)$ in \eqref{znudef} is the extension of the pairing between $\Lambda$ and $\Lambda^{*}$ to elements from $\Lambda_{\mathbb{C}}$. The functions $z_{\nu}$ from \eqref{znudef} satisfy $z_{\mu+\nu}=z_{\mu}z_{\nu}$ for $\mu$ and $\nu$ in $\Lambda^{*}$, a relation that we henceforth implicitly assume wherever we use this notation. It follows that as an algebraic variety we have \[T_{\Lambda}=\operatorname{Spec}\mathbb{C}[z_{\nu}|\;\nu\in\Lambda^{*}].\]

A \emph{toric variety} is an irreducible algebraic variety $X$ containing a torus $T$ as a dense Zariski open subset, such that the action of $T$ on itself extends to an action of $T$ on all of $X$. We shall only consider \emph{normal separated} toric varieties in what follows, where a variety is always understood to be locally of finite type. Recall that a \emph{rational polyhedral cone} $\sigma\subseteq\Lambda_{\mathbb{R}}$ is a set generated by finitely many elements of $\Lambda$, with coefficients in $\mathbb{R}_{+}=\{a\in\mathbb{R}|\;a\geq0\}$. The \emph{dual cone} of $\sigma$ is \[\sigma^{*}=\big\{\nu\in\Lambda^{*}_{\mathbb{R}}\big|\;(\nu,\lambda)\geq0\mathrm{\ for\ all\ }\lambda\in\sigma\big\},\] which is a rational polyhedral cone in $\Lambda^{*}_{\mathbb{R}}$. The cone $\sigma$ is called \emph{strongly convex} if it contains no full lines. This is equivalent to $\sigma^{*}$ being of full dimension $\operatorname{rk}\Lambda$ in $\Lambda^{*}_{\mathbb{R}}$. For any such cone $\sigma$, the semigroup $\sigma^{*}\cap\Lambda^{*}$ is finitely generated, and \[X_{\Lambda,\sigma}=\operatorname{Spec}\mathbb{C}[z_{\nu}|\;\nu\in\sigma^{*}\cap\Lambda^{*}]\] is a normal irreducible affine algebraic variety, containing $T_{\Lambda}(\mathbb{C})$ as the open acting torus.

We shall use the notation $\tau\leq\sigma$ for stating that $\tau$ is a \emph{face} of $\sigma$ (in the usual polyhedral sense). In this case $X_{\Lambda,\tau}$ is an open $T_{\Lambda}(\mathbb{C})$-invariant subvariety of $X_{\Lambda,\sigma}$. General (normal separated) toric varieties are then defined using the following notion.
\begin{defn} \label{fandef}
A \emph{rational fan} is a collection $\Sigma$ of strongly convex rational polyhedral cones in $\Lambda_{\mathbb{R}}$, satisfying the following conditions:
\begin{enumerate}[$(i)$]
\item If $\sigma\in\Sigma$ and $\tau\leq\sigma$ then $\tau\in\Sigma$.
\item Given $\sigma$ and $\tau$ in $\Sigma$, we have $\sigma\cap\tau\leq\sigma$ and $\sigma\cap\tau\leq\tau$.
\end{enumerate}
\end{defn}
The first main result in the theory of toric varieties is the following one.
\begin{thm} \label{torvar}
Every rational fan $\Sigma$ as in Definition \ref{fandef} gives rise to a separated normal toric variety $X_{\Lambda,\Sigma}$, on which the acting torus is $T_{\Lambda}$. Conversely, every separated normal toric variety arises, up to isomorphism, as $X_{\Lambda,\Sigma}$ from a lattice $\Lambda$ and a fan $\Sigma$ satisfying the conditions of Definition \ref{fandef}.
\end{thm}

\smallskip

Theorem \ref{torvar} indicates that properties of the toric varieties $X_{\Lambda,\sigma}$ should be expressed as properties of the defining fan $\Sigma$. We recall that the \emph{dimension} of a cone $\sigma$ is the dimension of the real vector space $\mathbb{R}\sigma=\sigma+(-\sigma)$.
\begin{defn} \label{smoothdef}
A cone $\sigma$, of some dimension $0 \leq l\leq\operatorname{rk}\Lambda$, is called \emph{simplicial} if it can be generated using precisely $l$ elements of $\Lambda$. It is called \emph{smooth} if these elements form a basis over $\mathbb{Z}$ for the group $\Lambda\cap\mathbb{R}\sigma$. The fan $\Sigma$ from Definition \ref{fandef} will be called \emph{simplicial} or \emph{smooth} if every $\sigma\in\Sigma$ has the respective property.
\end{defn}
We state the following relations between properties of $\Sigma$ and properties of $X_{\Lambda,\Sigma}$, whose proofs are rather simple and straightforward.
\begin{prop} \label{smoothtor}
The following assertions hold:
\begin{enumerate}[$(i)$]
\item The affine toric variety $X_{\Lambda,\sigma}$ is smooth if and only if the cone $\sigma$ is smooth in the sense of Definition \ref{smoothdef}. The general toric variety $X_{\Lambda,\Sigma}$ is smooth if and only if the fan $\Sigma$ is smooth.
\item If the cone $\sigma$ is simplicial then $X_{\Lambda,\sigma}$ has at most finite quotient singularities. If the fan $\Sigma$ is simplicial then the toric variety $X_{\Lambda,\Sigma}$ has at most finite quotient singularities.
\item The toric variety $X_{\Lambda,\Sigma}$ is of finite type if and only if the fan $\Sigma$ is finite.
\item The toric variety $X_{\Lambda,\Sigma}$ is proper if and only if $\Sigma$ is finite and $\bigcup_{\sigma\in\Sigma}\sigma$ is all of $\Lambda_{\mathbb{R}}$.
\end{enumerate}
\end{prop}
We remark that while our main object of interest below will be \emph{smooth} toroidal compactifications, which are obtained when the toric varieties appearing in the construction below are smooth (as characterized in Proposition \ref{smoothtor}), this fact will play no role in any of the following proofs. Our results are therefore equally valid for non-smooth toroidal compactifications.

\smallskip

The orbits in the action of $T_{\Lambda}$ on $X_{\Lambda,\Sigma}$ are in one-to-one correspondence with the cones $\sigma\in\Sigma$. More explicitly, define the annihilator of $\sigma$ to be \[\sigma^{0}=\big\{\nu\in\Lambda_{\mathbb{R}}\big|\;(\nu,\lambda)=0\mathrm{\ for\ all\ }\lambda\in\sigma\big\}=\sigma^{*}\cap(-\sigma^{*}),\] and then the \emph{orbit} $O(\sigma)$ associated with $\sigma$ and the corresponding closure $O_{\Sigma}(\sigma)$ in $X_{\Lambda,\Sigma}$ are defined to be
\begin{equation} \label{Osigmadef}
O(\sigma)=\mathcal{Z}_{X_{\Lambda,\sigma}}\big(\big\{z_{\nu}\big|\;\nu\in(\sigma^{*}\cap\Lambda^{*})\setminus\sigma^{0}\big\}\big)\quad\mathrm{and}\quad O_{\Sigma}(\sigma)=\bigcup_{\sigma\leq\eta}O(\eta).
\end{equation}
Here we have used, for an affine variety $X$, the notation $\mathcal{Z}_{X}(S)$ for the set of common zeros of a set $S\subseteq\mathbb{C}[X]$. The dimension of both $O(\sigma)$ and $O_{\Sigma}(\sigma)$ is $\operatorname{rk}\Lambda-\dim\sigma$. The special case of \eqref{Osigmadef} concerning divisors is given by \[D_{\rho}=O_{\Sigma}(\rho)\in\operatorname{Div}(X_{\Lambda,\Sigma})\] for a \emph{ray}, i.e., a 1-dimensional cone $\rho\in\Sigma$, and the generic part $O(\rho)=D_{\rho} \cap X_{\Lambda,\rho}$ of $D_{\rho}$ is smooth. Any such ray is $\mathbb{R}_{+}\omega$ for a unique, well-defined primitive element $\omega\in\Lambda$, and if we take some element $\kappa\in\Lambda^{*}$ with $(\kappa,\omega)=1$, which is possible by the primitivity of $\omega$, then $z_{\kappa}$ is a defining function for $O(\rho)$ in $X_{\Lambda,\rho}$, hence a local defining function for $D_{\rho}$ in $X_{\Lambda,\Sigma}$. We deduce the following beautiful formula from \cite{[Fu]} and others, which will be useful for us in what follows.
\begin{prop} \label{ordTdiv}
For $\nu\in\Lambda^{*}$ and a ray $\rho=\mathbb{R}_{+}\omega\in\Sigma$ as above, the order $\operatorname{ord}_{D_{\rho}}(z_{\nu})$ of the function $z_{\nu}$ for $\nu\in\Lambda^{*}$ along the the divisor $D_{\rho}$ is $(\nu,\omega)$.
\end{prop}

\smallskip

\subsection{Toroidal Compactifications}

The theory of toric varieties can be used to construct interesting (e.g., smooth) compactifications of algebraic varieties arising as quotients of Hermitian symmetric spaces by arithmetic groups. Indeed, let $\mathcal{D}$ be a simple Hermitian symmetric space of non-compact type, and let $G$ be the connected component of $\operatorname{Aut}(\mathcal{D})$. The group $G$ can be viewed as the connected component of the group of real points of a simple algebraic group $\mathcal{G}$. In case $\mathcal{G}$ is defined over $\mathbb{Q}$, the group $G$ contains arithmetic subgroups, which are discrete and of finite co-volume in $G$. By fixing such a subgroup $\Gamma$, the space $X=\Gamma\backslash\mathcal{D}$ carries the structure of a complex algebraic variety, which is usually defined over $\mathbb{Q}$ or over some number field (see, e.g., Theorem \ref{BBS} below), but is not complete in general.

Any compactification of $X$ arises by adding appropriate boundary components to $X$. The first, and most natural, way to do this is based on the boundary components of $\mathcal{D}$ as a symmetric space. Roughly, $\mathcal{D}$ admits a \emph{Harish-Chandra embedding} as a bounded open subset of a complex vector space (as in, e.g., Theorem 1 in Section III.2 of \cite{[AMRT]}), and its complement inside its closure inside this vector space consists of many isomorphic copies of Hermitian symmetric spaces of smaller dimensions, which are the boundary components in question (see the beginning of Section III.3 of \cite{[AMRT]} for more details). We denote the assertion that $\Xi$ is a boundary component of $\mathcal{D}$ by $\Xi\leq\mathcal{D}$, and the rational structure on $\mathcal{G}$ determines which boundary components of $\mathcal{D}$ are rational (this is always the case for $\mathcal{D}$ itself). One therefore defines \[\mathcal{D}^{\mathrm{BB}}=\bigcup_{\Xi\leq\mathcal{D},\ \text{$\Xi$ rational}}\Xi.\] The group $\Gamma$ acts on $\mathcal{D}^{\mathrm{BB}}$, and we set
\begin{equation} \label{XBBdef}
X^{\mathrm{BB}}=\Gamma\backslash\mathcal{D}^{\mathrm{BB}}\stackrel{open}{\supseteq}\Gamma\backslash\mathcal{D}=X.
\end{equation}
The quotient $X^{\mathrm{BB}}$ from \eqref{XBBdef}, with the Baily--Borel topology and complex structure, is called the \emph{Baily--Borel compactification} of $X$. It is a normal projective variety, which contains $X$ as a Zariski dense open subset. The boundary components of $X^{\mathrm{BB}}$, which are called \emph{cusps}, are quotients of rational boundary components of $\mathcal{D}$ by arithmetic groups. They are typically of high codimension and very singular. This compactification can also be constructed and characterized by the following result of Baily, Borel, and Satake (see, e.g., \cite{[Sa]} and \cite{[BB]}).
\begin{thm} \label{BBS}
The algebraic variety $X^{\mathrm{BB}}$ from \eqref{XBBdef} can be obtained as the $\operatorname{proj}$ construction of the ring of holomorphic modular forms with respect to $\Gamma$ (hence it is projective and defined over a number field). It is the minimal algebraic compactification of $X$.
\end{thm}
The minimality in Theorem \ref{BBS} means that any algebraic compactification of $X$ carries a unique map to $X^{\mathrm{BB}}$ that restricts to the identity map on the open subset $X$ of both varieties.

\smallskip

Smooth compactifications of $X$ can be constructed using the theory of \emph{toroidal} compactifications, obtained by gluing to $X$ open subsets of toric varieties. This theory is described in detail in \cite{[AMRT]} in general, and more explicitly for the symplectic case in \cite{[Nam]}. We will give the precise details for our case of interest, the orthogonal case, in Section \ref{CompGrass} below, slightly in parallel to \cite{[Fi]}, with a special case appearing also in Section 4.2 of \cite{[P]}. For describing the general idea, recall that the minimality in Theorem \ref{BBS} allows one to restrict attention to the pre-image of a given cusp, and any cusp is the $\Gamma$-orbit of some rational boundary component $\Xi\leq\mathcal{D}$. Let $\mathcal{P}(\Xi)$ denote the connected component of the stabilizer $\operatorname{St}_{G}(\Xi)$, which is a parabolic subgroup of $G$, for which \cite{[AMRT]} proves the following result (see Theorem 3 on Section III.3 and the beginning of Section III.4 of that reference).
\begin{lem} \label{genparstruc}
The parabolic subgroup $\mathcal{P}(\Xi)$ lies in the split short exact sequence \[1\to\mathcal{W}(\Xi)\to\mathcal{P}(\Xi) \to G_{h}(\Xi) \times G_{\ell}(\Xi) \times M(\Xi)\to1.\] The unipotent radical $\mathcal{W}(\Xi)$ of $\mathcal{P}(\Xi)$ lies in the short exact sequence \[0\to\mathcal{U}(\Xi)\to\mathcal{W}(\Xi)\to\mathcal{V}(\Xi)\to0,\] which does not split in general, where the center $\mathcal{U}(\Xi)$ of $\mathcal{W}(\Xi)$ and the quotient $\mathcal{V}(\Xi)$ are isomorphic to real vector spaces, with $\dim_{\mathbb{R}}\mathcal{V}(\Xi)$ even. The group $G_{h}(\Xi)$ is the connected component of $\operatorname{Aut}(\Xi)$. The group $G_{\ell}(\Xi)$ is a simple subgroup of $\operatorname{GL}\big(\mathcal{U}(\Xi)\big)$, whose symmetric space is a self-adjoint cone $\Omega(\Xi)$ inside $\mathcal{U}(\Xi)$ which is usually non-Hermitian, and $M(\Xi)$ is a compact group.
\end{lem}
We remark that some of the five components $G_{h}(\Xi)$, $G_{\ell}(\Xi)$, $M(\Xi)$, $\mathcal{U}(\Xi)$, or $\mathcal{V}(\Xi)$ from Lemma \ref{genparstruc} may be trivial for some $\mathcal{D}$ and $\Xi$ (see, e.g., Proposition \ref{paradim0} below). Note that the center $\mathcal{U}(\Xi)$ of $\mathcal{W}(\Xi)$ is normal also in $\mathcal{P}(\Xi)$, and we denote the quotient group $\mathcal{P}(\Xi)/\mathcal{U}(\Xi)$ by $\overline{\mathcal{P}}(\Xi)$.

The choice of $\Xi$ also gives a presentation of $\mathcal{D}$ in a specific way, in correspondence with the action of $\mathcal{P}(\Xi)$ as decomposed in Lemma \ref{genparstruc}. We refer the reader to Subsection III.4.3 of \cite{[AMRT]} for the details.
\begin{prop} \label{DoverXi}
The choice of $\Xi$ describes the Hermitian symmetric domain $\mathcal{D}$ as a topologically trivial holomorphic fiber bundle
\begin{equation} \label{holfib}
\mathcal{U}(\Xi)+i\Omega(\Xi)\to\mathcal{D}\stackrel{\pi}{\to}\widetilde{\mathcal{D}}(\Xi)
\end{equation}
over a complex manifold $\widetilde{\mathcal{D}}(\Xi)$, in which the fiber of $\mathcal{D}$ over a point in $\widetilde{\mathcal{D}}(\Xi)$ is the image of $\mathcal{U}(\Xi)+i\Omega(\Xi)$ inside an affine space underlying the complexification $\mathcal{U}(\Xi)_{\mathbb{C}}$ of the real vector space $\mathcal{U}(\Xi)$ from Lemma \ref{genparstruc}. Appropriately modifying the imaginary part of the fiber of the first fibration yields a map $\Phi:\mathcal{D}\to\Omega(\Xi)$, and with it a topologically trivial affine vector bundle structure
\begin{equation} \label{Phifib}
\mathcal{U}(\Xi)\to\mathcal{D}\stackrel{\pi\times\Phi}{\longrightarrow}\widetilde{\mathcal{D}}(\Xi)\times\Omega(\Xi).
\end{equation}
The manifold $\widetilde{\mathcal{D}}(\Xi)$ is an affine holomorphic vector bundle, whose underlying real structure is a topologically trivial affine vector bundle
\begin{equation} \label{vecfib}
\mathcal{V}(\Xi)\to\widetilde{\mathcal{D}}(\Xi)\to\Xi,
\end{equation}
in which the fiber $\mathcal{V}(\Xi)_{s}$ of $\widetilde{\mathcal{D}}(\Xi)$ over $s\in\Xi$ is the real affine space underlying $\mathcal{V}(\Xi)$, but carrying a complex structure that depends holomorphically on $s$.
\end{prop}
The action of the groups from Lemma \ref{genparstruc} on the fibrations from Proposition \ref{DoverXi} is also described in Subsection III.4.3 of \cite{[AMRT]}, and is as follows.
\begin{prop} \label{subgpact}
The center $\mathcal{U}(\Xi)$ of the unipotent radical $\mathcal{W}(\Xi)$ from Lemma \ref{genparstruc} acts additively on the affine fibers of the bundle from \eqref{Phifib}, or equivalently on the real part of those from \eqref{holfib}. The quotient $\mathcal{V}(\Xi)$ acts similarly on the fibers from \eqref{vecfib}, and the map $\Phi$ from Proposition \ref{DoverXi} is invariant under this action. In the resulting action of the Levi quotient $G_{h}(\Xi) \times G_{\ell}(\Xi) \times M(\Xi)$ from Lemma \ref{genparstruc} on the quotient space $\Xi\times\Omega(\Xi)$, the groups $G_{h}(\Xi)$ and $G_{\ell}(\Xi)$ act on their symmetric spaces $\Xi$ and $\Omega(\Xi)$ respectively, the latter one via $\Phi$, and $M(\Xi)$ acts trivially.
\end{prop}

\begin{rmk} \label{VHS}
Complex analytic varieties like $\widetilde{\mathcal{D}}(\Xi)$ show up from another point of view when one studies \emph{variation of Hodge structures}. Indeed, Proposition \ref{DoverXi} describes $\widetilde{\mathcal{D}}(\Xi)$ as a variation of Hodge structure over $\Xi$. Proposition \ref{subgpact} shows that the fibers in \eqref{vecfib} are \emph{principal homogeneous spaces} of $\mathcal{V}(\Xi)$, and those from \eqref{Phifib} are principal homogeneous spaces of $\mathcal{U}(\Xi)$. The boundary component $\Xi$ is understood to lie ``at $\infty$'', in the direction of the cone $\Omega(\Xi)$.
\end{rmk}

\smallskip

Recalling that $X=\Gamma\backslash\mathcal{D}$ for an arithmetic group $\Gamma$, we set
\begin{equation} \label{grpsindZ}
\mathcal{P}_{\mathbb{Z}}(\Xi)=\mathcal{P}(\Xi)\cap\Gamma,\qquad\mathcal{W}_{\mathbb{Z}}(\Xi)=\mathcal{W}(\Xi)\cap\Gamma,\qquad\mathcal{U}_{\mathbb{Z}}(\Xi)=\mathcal{U}(\Xi)\cap\Gamma
\end{equation}
as well as \[\mathcal{V}_{\mathbb{Z}}(\Xi)=\mathcal{W}_{\mathbb{Z}}(\Xi)/\mathcal{U}_{\mathbb{Z}}(\Xi)\subseteq\mathcal{V}(\Xi)\] and
\begin{equation} \label{quotsZ}
\overline{\mathcal{P}}_{\mathbb{Z}}(\Xi)=\mathcal{P}_{\mathbb{Z}}(\Xi)/\mathcal{U}_{\mathbb{Z}}(\Xi)\subseteq\overline{\mathcal{P}}(\Xi)=\mathcal{P}(\Xi)/\mathcal{U}(\Xi).
\end{equation}
Now, the spaces $\mathcal{U}(\Xi)$ and $\mathcal{U}(\Xi)_{\mathbb{C}}$ are obtained from $\mathcal{U}_{\mathbb{Z}}(\Xi)$ by extending the scalars to $\mathbb{R}$ and $\mathbb{C}$ respectively, and by Proposition \ref{subgpact}, the group $\mathcal{U}_{\mathbb{Z}}(\Xi)\subseteq\mathcal{U}(\Xi)$ acts on the real part of the fibers of the bundle from \eqref{holfib}. Hence $\mathcal{U}_{\mathbb{Z}}(\Xi)\backslash\mathcal{D}$ is also a fiber bundle over $\widetilde{\mathcal{D}}(\Xi)$, with the fibers being affine models of
\begin{equation} \label{affTUZ}
[\mathcal{U}(\Xi)/\mathcal{U}_{\mathbb{Z}}(\Xi)]+i\Omega(\Xi)\subseteq[\mathcal{U}(\Xi)/\mathcal{U}_{\mathbb{Z}}(\Xi)]+i\mathcal{U}(\Xi)=
\mathcal{U}(\Xi)_{\mathbb{C}}/\mathcal{U}_{\mathbb{Z}}(\Xi)=T_{\mathcal{U}_{\mathbb{Z}}(\Xi)}(\mathbb{C}).
\end{equation}
In these coordinates, the component $\Xi$ of the set $\mathcal{D}^{\mathrm{BB}}$ from \eqref{XBBdef} lies at $\infty$ in the direction of $\Omega(\Xi)$. The idea is to extend the fibers $T_{\mathcal{U}_{\mathbb{Z}}(\Xi)}$ to appropriate toric varieties, such that the point $s\in\Xi\subseteq\mathcal{D}^{\mathrm{BB}}$ will be replaced by toroidal boundary components. This will make the boundary not only a bundle over the boundary component $\Xi$ of $X^{\mathrm{BB}}$ itself, but also over $\widetilde{\mathcal{D}}(\Xi)$.

In order to do so while respecting the action of $\Gamma$, one observes that there is an open subset $\Omega_{N}(\Xi)\subseteq\Omega(\Xi)$, which is invariant under the image of $\mathcal{P}_{\mathbb{Z}}(\Xi)$ in $G_{\ell}(\Xi)$, which satisfies
\begin{equation} \label{OmegaNXi}
\Omega_{N}(\Xi)+\Omega(\Xi)\subseteq\Omega_{N}(\Xi),
\end{equation}
and such that if $\xi\in\mathcal{D}$ with $\Phi(\xi)\in\Omega_{N}(\Xi)$ and $A\in\Gamma$ does not lie in $\mathcal{P}_{\mathbb{Z}}(\Xi)$ then $\Phi(A\xi)\not\in\Omega_{N}(\Xi)$. This yields the following neighborhood of $\Xi$ in $\mathcal{D}$.
\begin{lem} \label{NXi}
The subset $N_{\Xi}=\Phi^{-1}\big(\Omega_{N}(\Xi)\big)$ of $\mathcal{D}$ is $\mathcal{P}_{\mathbb{Z}}(\Xi)$-invariant, and any element $A\in\Gamma$ which satisfies $AN_{\Xi} \cap N_{\Xi}\neq\emptyset$ must be in $\mathcal{P}_{\mathbb{Z}}(\Xi)$.
\end{lem}
Lemma \ref{NXi} and \eqref{OmegaNXi} imply that that $\mathcal{P}_{\mathbb{Z}}(\Xi) \backslash N_{\Xi}$ becomes a neighborhood of the corresponding cusp of $X$. It is clear from the definition of $N_{\Xi}$ in Lemma \ref{NXi} that the restriction of the map $\pi$ from \eqref{holfib} to $N_{\Xi}$ still surjects onto the manifold $\widetilde{\mathcal{D}}(\Xi)$ from Proposition \ref{DoverXi}. Dividing by the action of $\mathcal{U}_{\mathbb{Z}}(\Xi)$, or equivalently restricting \eqref{affTUZ}, gives the fiber bundle
\begin{equation} \label{Nfib}
[\mathcal{U}(\Xi)/\mathcal{U}_{\mathbb{Z}}(\Xi)]+i\Omega_{N}(\Xi)\to\mathcal{U}_{\mathbb{Z}}(\Xi) \backslash N_{\Xi}\to\widetilde{\mathcal{D}}(\Xi),
\end{equation}
with the fibers being smaller open subsets of $T_{\mathcal{U}_{\mathbb{Z}}(\Xi)}(\mathbb{C})$.

The toric variety extending $T_{\mathcal{U}_{\mathbb{Z}}(\Xi)}$ will therefore have to be based on a fan $\Sigma(\Xi)$ of cones that are roughly contained in $\Omega(\Xi)$ (this will be made more precise soon). But it will also have to respect the action of the image of $\mathcal{P}_{\mathbb{Z}}(\Xi)$ in $G_{\ell}(\Xi)$. In addition, consider a cusp $\Upsilon$ whose closure contains $\Xi$, a situation that we denote by $\Xi\leq\Upsilon$. After embedding the Levi quotients inside $\mathcal{P}(\Xi)$ and $\mathcal{P}(\Upsilon)$ one obtains the following result, appearing in Theorem 3 of Chapter 3 of \cite{[AMRT]}.
\begin{lem} \label{relsubgp}
Let $\Xi$ and $\Upsilon$ be cusps of $\mathcal{D}$, and assume that $\Xi\leq\Upsilon$. We then have the inclusion $G_{h}(\Xi) \subseteq G_{h}(\Upsilon)$. One the other hand, we get $\mathcal{U}(\Upsilon)\subseteq\mathcal{U}(\Xi)$ and $G_{\ell}(\Upsilon) \subseteq G_{\ell}(\Xi)$ in the other direction, with the latter embedding being not canonical, and the cone $\Omega(\Upsilon)$ is a rational boundary component of $\Omega(\Xi)$. All the rational boundary components of $\Omega(\Xi)$ are $\Omega(\Upsilon)$ for $\Xi\leq\Upsilon$.
\end{lem}
The exact properties of each fan $\Sigma(\Xi)$, as well as of the collection of fans $\{\Sigma(\Xi)\}_{\Xi}$, are expressed in the following definition, taken from \cite{[AMRT]} and \cite{[Nam]} (see, e.g., Subsection II.4.3 and III.5 of the former reference).
\begin{defn} \label{admissible}
For $\mathcal{D}$, $\Gamma$, and $\Xi$ as above we say that a fan $\Sigma(\Xi)$ is an \emph{admissible cone decomposition} of (the closure of) $\Omega(\Xi)$ if the following three conditions are satisfied:
\begin{enumerate}[$(i)$]
\item $\bigcup_{\sigma\in\Sigma(\Xi)}\sigma=\bigcup_{\Xi\leq\Upsilon}\Omega(\Upsilon)$. The union on the right hand side is the \emph{closure} of $\Omega(\Xi)$.
\item If $\sigma\in\Sigma(\Xi)$ and $B \in G_{\ell}(\Xi)$ is in the image of the group $\overline{\mathcal{P}}_{\mathbb{Z}}(\Xi)$ from \eqref{quotsZ} then $B\sigma\in\Sigma(\Xi)$, i.e., $\Sigma(\Xi)$ is $\overline{\mathcal{P}}_{\mathbb{Z}}(\Xi)$-invariant.
\item The action of $\overline{\mathcal{P}}_{\mathbb{Z}}(\Xi)$ on $\Sigma(\Xi)$ has finitely many orbits.
\end{enumerate}
For a collection $\{\Sigma(\Xi)\}_{\Xi\leq\mathcal{D}}$ of fans to be called \emph{admissible} the following three following must hold:
\begin{enumerate}[$(i)$]
\item $\Sigma(\Xi)$ is an admissible cone decomposition of the closure of $\Omega(\Xi)$ for every $\Xi\leq\mathcal{D}$.
\item If $\Xi\leq\Upsilon$ then $\Sigma(\Upsilon)=\big\{\sigma\in\Sigma(\Xi)\big|\;\sigma\subseteq\bigcup_{\Upsilon\leq\Psi}\Omega(\Psi)\big\}$.
\item For $\Xi\leq\mathcal{D}$ and $A\in\Gamma$ we have $\Sigma(A\Xi)=\big\{A\sigma\big|\;\sigma\in\Sigma(\Xi)\big\}$.
\end{enumerate}
\end{defn}
For two cusps $\Xi$ and $\Upsilon$ with $\Xi\leq\Upsilon$, we consider $\Omega(\Upsilon)$, as well as $\Omega(\Psi)$ for every $\Upsilon\leq\Psi$, as identified with its image as a boundary component of $\Omega(\Xi)$ as in Lemma \ref{relsubgp}. Note that in the second condition on $\Sigma(\Xi)$ in Definition \ref{admissible} we consider $\overline{\mathcal{P}}_{\mathbb{Z}}(\Xi)$-invariance, and not the equivalent $\mathcal{P}_{\mathbb{Z}}(\Xi)$-invariance, since we already divided by $\mathcal{U}_{\mathbb{Z}}(\Xi)$ in order to be able to work with $T_{\mathcal{U}_{\mathbb{Z}}(\Xi)}$. This condition indicates that $\Sigma(\Xi)$ will usually contain infinitely many cones, since the image of $\overline{\mathcal{P}}_{\mathbb{Z}}(\Xi)$ in $G_{\ell}(\Xi)$ is typically infinite.

The existence of admissible cone decompositions for $\Omega(\Xi)$, as well as admissible collections of fans $\{\Sigma(\Xi)\}_{\Xi\leq\mathcal{D}}$, is not hard to prove: One way is by slightly modifying the argument leading to the Main Theorem II in Section III.7 of \cite{[AMRT]} (the assumption that $\Gamma$ is neat there poses no restriction, since any arithmetic group contains a neat subgroup of finite index).

\smallskip

Given an admissible cone decomposition $\Sigma(\Xi)$ for some $\Xi$, recall the fiber bundles from \eqref{affTUZ} and \eqref{Nfib}, and that the fibers are open subsets of affine spaces over the torus $T_{\mathcal{U}_{\mathbb{Z}}(\Xi)}(\mathbb{C})$. Embedding this torus inside $X_{\mathcal{U}_{\mathbb{Z}}(\Xi),\Sigma(\Xi)}$, one gets a decomposition
\begin{equation} \label{intclos}
\overline{[\mathcal{U}(\Xi)/\mathcal{U}_{\mathbb{Z}}(\Xi)]+i\Omega(\Xi)}^{o}=
\bigcup_{\sigma\in\Sigma(\Xi)}\Big[\overline{[\mathcal{U}(\Xi)/\mathcal{U}_{\mathbb{Z}}(\Xi)]+i\Omega(\Xi)}^{o} \cap O(\sigma)\Big],
\end{equation}
where the set on the left hand side is the relative interior of the closure of the fiber $[\mathcal{U}(\Xi)/\mathcal{U}_{\mathbb{Z}}(\Xi)]+i\Omega(\Xi)$ from \eqref{affTUZ} in $X_{\mathcal{U}_{\mathbb{Z}}(\Xi),\Sigma(\Xi)}$. It is clear that none of the sets in the union from \eqref{intclos} is empty, that the intersection associated with any $\sigma\in\Sigma(\Xi)$ that does not lie in any $\Sigma(\Upsilon)$ for a boundary component $\Upsilon\neq\Xi$ with $\Xi\leq\Upsilon$ is the full orbit $O(\sigma)$, and that similar assertions holds for the fiber $[\mathcal{U}(\Xi)/\mathcal{U}_{\mathbb{Z}}(\Xi)]+i\Omega_{N}(\Xi)$ from \eqref{Nfib}, with all the terms in which $\sigma\neq\{0\}$ being the same. Let now $\mathcal{D}_{\Xi}$ be the manifold having the fiber bundle structure
\begin{equation} \label{fibDXi}
\bigcup_{\sigma\in\Sigma(\Xi)}\Big[\overline{[\mathcal{U}(\Xi)/\mathcal{U}_{\mathbb{Z}}(\Xi)]+i\Omega(\Xi)}^{o} \cap O(\sigma)\Big]\to\mathcal{D}_{\Xi}\to\widetilde{\mathcal{D}}(\Xi).
\end{equation}
In particular $\mathcal{D}_{\mathcal{D}}$ in \eqref{fibDXi} is just $\mathcal{D}$ again. Since $\Sigma(\Xi)$ is an admissible cone decomposition in the sense of Definition \ref{admissible}, the action of $\overline{\mathcal{P}}_{\mathbb{Z}}(\Xi)$ extends to $\mathcal{D}_{\Xi}$, and the relative interior $\overline{\mathcal{P}_{\mathbb{Z}}(\Xi) \backslash N_{\Xi}}^{o}$ of the closure of $\mathcal{P}_{\mathbb{Z}}(\Xi) \backslash N_{\Xi}$ in $\overline{\mathcal{P}}_{\mathbb{Z}}(\Xi)\big\backslash\mathcal{D}_{\Xi}$ is
\begin{equation} \label{NmodPZ}
\big(\mathcal{P}_{\mathbb{Z}}(\Xi) \backslash N_{\Xi}\big)\cup\bigg\{\overline{\mathcal{P}}_{\mathbb{Z}}(\Xi)\bigg\backslash\bigg[\bigcup_{\{0\}\neq\sigma\in\Sigma(\Xi)}\Big[\overline{[\mathcal{U}(\Xi)/\mathcal{U}_{\mathbb{Z}}(\Xi)]+i\Omega(\Xi)}^{o} \cap O(\sigma)\Big]\bigg]\bigg\}.
\end{equation}
The union in \eqref{NmodPZ} is essentially a finite union of $T_{\mathcal{U}_{\mathbb{Z}}(\Xi)}$-orbits of positive codimension (each one perhaps divided by the action of some finite group). A technical important result, in which we denote the image of $\overline{\mathcal{P}}_{\mathbb{Z}}(\Xi)$ in $G_{\ell}(\Xi)$ by $\Gamma_{\ell}(\Xi)$, is the following one (see Proposition 2 in Subsection III.6.3 of \cite{[AMRT]} and the surrounding discussion).
\begin{thm} \label{propdisc}
Cones in $\Sigma(\Xi)$ that are not contained in $\Sigma(\Upsilon)$ for any other cusp $\Upsilon$ with $\Xi\leq\Upsilon$ have finite stabilizers in $\Gamma_{\ell}(\Xi)$. This group therefore operates properly discontinuously on the torus $T_{\mathcal{U}_{\mathbb{Z}}(\Xi)}$, hence also on the fibers from \eqref{Nfib}. Moreover, when $\Sigma(\Xi)$ is an admissible cone decomposition, this action extends to a properly discontinuous action on $X_{\mathcal{U}_{\mathbb{Z}}(\Xi),\Sigma(\Xi)}$ and on the fibers of $\mathcal{D}_{\Xi}$ from \eqref{fibDXi}, taking boundary components to boundary components.
\end{thm}
Theorem \ref{propdisc} is complemented by an equivalent to part $(iii)$ in Proposition \ref{smoothtor}, which states that $\Gamma_{\ell}(\Xi) \backslash X_{\mathcal{U}_{\mathbb{Z}}(\Xi),\Sigma(\Xi)}$ is an open subset of a variety that is of finite type, because of the finiteness of the number of orbits in Definition \ref{admissible}. We also deduce from that theorem that if the cones in $\Sigma(\Xi)$ are smooth as in Definition \ref{smoothdef}, then this quotient is smooth up to finite quotient singularities.

Assume now that the collection $\{\Sigma(\Xi)\}_{\Xi}$ is admissible as in Definition \ref{admissible}, and take two cusps $\Xi$ and $\Upsilon$ with $\Xi\leq\Upsilon$. The first condition in that definition and Lemma \ref{relsubgp} allow us to identify $\big(\mathcal{U}_{\mathbb{Z}}(\Xi)/\mathcal{U}_{\mathbb{Z}}(\Upsilon)\big)\big\backslash\mathcal{D}_{\Upsilon}$ as an open subset of $\mathcal{D}_{\Xi}$, and write
\begin{equation} \label{DUpsDXi}
\mathcal{D}_{\Xi}=\big[\big(\mathcal{U}_{\mathbb{Z}}(\Xi)/\mathcal{U}_{\mathbb{Z}}(\Upsilon)\big)\big\backslash\mathcal{D}_{\Upsilon}\big]\cup
\Big\{\bigcup_{\sigma\in\Sigma(\Xi)\setminus\Sigma(\Upsilon)}\Big[\overline{[\mathcal{U}(\Xi)/\mathcal{U}_{\mathbb{Z}}(\Xi)]+i\Omega(\Xi)}^{o} \cap O(\sigma)\Big]\Big\}.
\end{equation}
Recall that $\mathcal{W}(\Upsilon)$ acts trivially on $\Upsilon$, so that it leaves all the boundary components of $\Upsilon$ invariant, and in particular fixes $\Xi$. On the other hand, $\mathcal{W}(\Xi)$ acts trivially on $\Omega(\Xi)$, so that it fixes the boundary components of that cone that are described in Lemma \ref{relsubgp}, including, in particular, $\Omega(\Upsilon)$. We thus obtain the inclusions
\begin{equation} \label{intpar}
\mathcal{U}_{\mathbb{Z}}(\Upsilon)\subseteq\mathcal{U}_{\mathbb{Z}}(\Xi)\subseteq\mathcal{W}_{\mathbb{Z}}(\Xi)\subseteq\mathcal{P}_{\mathbb{Z}}(\Upsilon)\cap\mathcal{P}_{\mathbb{Z}}(\Xi)
\supseteq\mathcal{W}_{\mathbb{Z}}(\Upsilon),
\end{equation}
all of which holding also without the index $\mathbb{Z}$. The inclusion from \eqref{DUpsDXi} thus yields
\[\big[\big(\mathcal{P}_{\mathbb{Z}}(\Upsilon)\cap\mathcal{P}_{\mathbb{Z}}(\Xi)\big)/\mathcal{U}_{\mathbb{Z}}(\Upsilon)\big]\big\backslash\mathcal{D}_{\Upsilon}\stackrel{open}{\subseteq}
\big[\big(\mathcal{P}_{\mathbb{Z}}(\Upsilon)\cap\mathcal{P}_{\mathbb{Z}}(\Xi)\big)/\mathcal{U}_{\mathbb{Z}}(\Xi)\big]\big\backslash\mathcal{D}_{\Xi},\] using which we can glue the appropriate neighborhoods of the cusps in the sets $\overline{\mathcal{P}}_{\mathbb{Z}}(\Xi)\backslash\mathcal{D}_{\Xi}$ and $\overline{\mathcal{P}}_{\mathbb{Z}}(\Upsilon)\backslash\mathcal{D}_{\Upsilon}$ properly when $\Xi\leq\Upsilon$. Definition \ref{admissible} also implies that if two rational boundary components lie over the same cusp of $X^{\mathrm{BB}}$, which means that they are identified via $\Gamma$, then their toroidal pre-images coincide. Then \eqref{NmodPZ} and the defining properties of the set $N_{\Xi}$ from Lemma \ref{NXi} form the basis for the main result in the theory of toroidal compactifications, given as Main Theorem I in Section III.5 of \cite{[AMRT]}, which we now state.
\begin{thm} \label{torcomp}
Let $\mathcal{D}$ be a simple Hermitian symmetric space, let $\Gamma \leq G=\operatorname{Aut}(\mathcal{D})$ be an arithmetic subgroup, and let $\{\Sigma(\Xi)\}_{\Xi\leq\mathcal{D}}$ be admissible as in Definition \ref{admissible}. Set $X^{\mathrm{tor}}=X^{\mathrm{tor}}_{\{\Sigma(\Xi)\}_{\Xi\leq\mathcal{D}}}$ to be the variety obtained by gluing, for each $\Xi\leq\mathcal{D}$, the open subset $\mathcal{P}_{\mathbb{Z}}(\Xi) \backslash N_{\Xi}$ of $X=\Gamma\backslash\mathcal{D}$ with its natural image in $\overline{\mathcal{P}}_{\mathbb{Z}}(\Xi)\backslash\mathcal{D}_{\Xi}$ as in \eqref{NmodPZ}. Then $X^{\mathrm{tor}}$ is a compactification of $X$, which is a complete algebraic variety defined over the same field as $X$. Assume, in addition, that every cone in every $\Sigma(\Xi)$ is smooth as in Definition \ref{smoothdef}, and that if $\xi \in X^{\mathrm{tor}}$ is in the image of $\mathcal{D}_{\Xi}$ but not in $\mathcal{D}_{\Upsilon}$ for any other $\Xi\leq\Upsilon$ then the stabilizer of $\xi$ in $\overline{\mathcal{P}}_{\mathbb{Z}}(\Xi)$ is trivial. In this case $X^{\mathrm{tor}}$ is a smooth variety, hence a compact complex manifold.
\end{thm}
The main drawback of the construction from Theorem \ref{torcomp} is that it is not canonical---it depends on the choice of the admissible cone decompositions $\Sigma(\Xi)$ for $\Xi\leq\mathcal{D}$. However, we shall assume that the choice of admissible cone decompositions is fixed from the beginning, and write just $X^{\mathrm{tor}}$ in what follows, for avoiding the heavy notation $X^{\mathrm{tor}}_{\{\Sigma(\Xi)\}_{\Xi}}$.

\begin{rmk}
We remark that a standard argument involving barycentric subdivisions implies that every fan has a smooth refinement, and since for a $\Gamma_{\ell}(\Xi)$-invariant fan this can be done in a $\Gamma_{\ell}(\Xi)$-equivariant manner, which also respects the relations between the fans associated with different cusps as in Definition \ref{admissible}, admissible collections of smooth admissible cone decompositions always exist. On the other hand, the condition about triviality of stabilizers depend on the arithmetic group $\Gamma$, and indeed for some groups $\Gamma$ the open variety $X$ from \eqref{XBBdef} may itself be not smooth. The latter issue can be resolved, however, by taking a finite index subgroup of $\Gamma$ and the corresponding cover of $X$. \label{smexist}
\end{rmk}

\section{Compactifications of Grassmannian Quotients \label{CompGrass}}

This section describes the explicit form of the toroidal compactifications in the orthogonal case. The construction parallels \cite{[Fi]} a bit, but with emphasis on the parts that are required for our goals. The most important results here are the definition of the toroidal divisors in \eqref{BJBIomega}, the explicit admissibility conditions in Corollary \ref{admcdecom} and Lemma \ref{GammaL0dim}, and the construction in Theorem \ref{formtc}.

\subsection{Grassmannians with a Complex Structure}

Let $L$ be an even integral lattice of signature $(n,2)$. We shall denote the image of two vectors $\lambda$ and $\mu$ in $L$ under the associated bilinear form by $(\lambda,\mu)$, and we shorthand $(\lambda,\lambda)$ to $\lambda^{2}$ for every $\lambda \in L$. The quadratic form on $L$ therefore sends $\lambda \in L$ to $\frac{\lambda^{2}}{2}$. The group of automorphisms of the real quadratic space $L_{\mathbb{R}}$ is the Lie group $\operatorname{O}(L_{\mathbb{R}})\cong\operatorname{O}(n,2)$, and we denote the connected component of its identity by $\operatorname{SO}^{+}(L_{\mathbb{R}})\cong\operatorname{SO}^{+}(n,2)$.

The symmetric space of $\operatorname{O}(L_{\mathbb{R}})$ and of $\operatorname{SO}^{+}(L_{\mathbb{R}})$ is isomorphic to the \emph{Grassmannian}
\begin{equation} \label{Grassdef}
\operatorname{Gr}(L)=\big\{v \subseteq L_{\mathbb{R}}\big|\;v<0,\ \dim v=2\}=\big\{v \subseteq L_{\mathbb{R}}\big|\;v<0,\ v^{\perp}>0\}.
\end{equation}
It is a real manifold of dimension $2n$, but it carries a natural complex structure, derived from the following description (see Section 13 of \cite{[Bo1]}, Section 3.2 of \cite{[Br1]}, Section 1.2 of \cite{[Ze1]}, Section 3.1 of \cite{[Fi]}, and others, but note that here we work with the opposite signature). We set \[(L_{\mathbb{C}})_{0}=\mathcal{Z}_{L_{\mathbb{C}}}(Z_{L}^{2})=\big\{Z_{L} \in L_{\mathbb{C}}\big|\;Z_{L}^{2}=0\big\}\] to be the zero set of the quadratic form on $L_{\mathbb{C}}$, as well as
\begin{equation} \label{bunGV}
P=\big\{Z_{L}\in(L_{\mathbb{C}})_{0}\big|\;(Z_{L},\overline{Z}_{L})<0\big\},
\end{equation}
which is an open subset of $(L_{\mathbb{C}})_{0}$. We have a natural map
\begin{equation} \label{fiboverGVR}
P\to\operatorname{Gr}(L),\quad Z_{L}=X_{L}+iY_{L} \mapsto v=\mathbb{R}X_{L}\oplus\mathbb{R}Y_{L},
\end{equation}
which is clearly surjective. Moreover, the set $P$ from \eqref{bunGV} is the disjoint union $P^{+} \cup P^{-}$ of two connected components, distinguished by the orientation on the basis $(X_{L},Y_{L})$ for the space $v$ from \eqref{fiboverGVR}, and the restriction of the map from that equation to one component, say $P^{+}$, is still surjective. The natural action of $\operatorname{O}(L_{\mathbb{R}})$ on $L_{\mathbb{C}}$ preserves $P$, the subgroup $\operatorname{SO}^{+}(L_{\mathbb{R}})$ preserves the two components $P^{+}$ and $P^{-}$, and the map from \eqref{fiboverGVR} is $\operatorname{O}(L_{\mathbb{R}})$-equivariant. In total, we get the following result.
\begin{prop} \label{GVcomp}
The Grassmannian $\operatorname{Gr}(L)$ from \eqref{Grassdef} is diffeomorphic to the image of the connected component $P^{+}$ from \eqref{bunGV} in the projective space $\mathbb{P}(L_{\mathbb{C}})$. This image is an analytically open subset of the projective quadric $\mathbb{P}(L_{\mathbb{C}})_{0}=\mathcal{Z}_{\mathbb{P}(L_{\mathbb{C}})}(Z_{L}^{2})$. Therefore $\operatorname{Gr}(L)$ carries the structure of an $n$-dimensional complex manifold, on which $\operatorname{SO}^{+}(L_{\mathbb{R}})$ operates holomorphically and transitively with compact stabilizers.
\end{prop}

\smallskip

Following, e.g., the description in Section 2 of \cite{[BFr]}, one finds that the closure $\overline{\operatorname{Gr}(L)}$ of $\operatorname{Gr}(L)$ in $\mathbb{P}(L_{\mathbb{C}})_{0}$ is obtained by adding appropriate boundary components to $\operatorname{Gr}(L)$, which are associated with isotropic subspaces of $L_{\mathbb{R}}$. To understand them, we first observe that the closure $\overline{P^{+}}$ of $P^{+}$ in $(L_{\mathbb{C}})_{0}$ is \[\big\{Z_{L}=X_{L}+iY_{L} \in L_{\mathbb{C}}\big|\;Z_{L}^{2}=0,\ (Z_{L},\overline{Z}_{L})\leq0,\ (X_{L},Y_{L})\mathrm{\ linearly\ dependent\ or\ oriented}\big\}.\] In the projective image we omit the point $Z_{L}=0$, and the boundary arises from the complement of $P^{+}$ itself, so that the boundary of $\operatorname{Gr}(L)$ in $\mathcal{Z}_{\mathbb{P}(L_{\mathbb{C}})}(Z_{L}^{2})$ is the projective image of
\begin{equation} \label{bdP+}
\big\{0 \neq Z_{L} \in L_{\mathbb{C}}\big|\;Z_{L}^{2}=0,\ (Z_{L},\overline{Z}_{L})=0,\ (X_{L},Y_{L})\mathrm{\ linearly\ dependent\ or\ oriented}\big\}.
\end{equation}
Recalling the map from \eqref{fiboverGVR}, taking a vector to the space spanned by its real and imaginary part, we consider its behavior on the set from \eqref{bdP+}. Since now $X_{L}$ and $Y_{L}$ are perpendicular and isotropic, the space that they span, whose complexification is the span of $Z_{L}$ and $\overline{Z}_{L}$ over $\mathbb{C}$, is isotropic. Its dimension cannot be 0 (since $Z_{L}=0$ is excluded in \eqref{bdP+}), hence this space is either a line or a plane. It is a line if and only if $X_{L}$ and $Y_{L}$ are linearly dependent over $\mathbb{R}$, or equivalently $Z_{L}$ and $\overline{Z}_{L}$ are linearly dependent over $\mathbb{C}$. Then all the possible generators $Z_{L}$ of the complexification of this line have the same projective image. The resulting point in $(L_{\mathbb{C}})_{0}$ is a 0-dimensional boundary component of $\operatorname{Gr}(L)$.

Otherwise the dimension is 2, and we choose a fixed oriented basis $(z,w)$ for the resulting isotropic plane. Since $X_{L}$ and $Y_{L}$ are linearly independent, we find that $Z_{L}$ cannot be spanned over $\mathbb{C}$ by $w$ alone, and can therefore we be written as $\beta(\tau w+z)$ for some complex $\beta\neq0$ and $\tau=x+iy$. Moreover, the linear independence shows that $\tau\not\in\mathbb{R}$, and the orientation implies that $y>0$, i.e., $\tau$ lies in the \emph{upper half-plane} $\mathcal{H}=\{\tau\in\mathbb{C}|\;\Im\tau>0\}$. In the projective image we may ignore only the factor $\beta$, and hence the boundary component of $\overline{\operatorname{Gr}(L)}$ that is associated with our isotropic plane is the projective image of the set
\begin{equation} \label{coor1dcusp}
\{\tau w+z|\;\tau\in\mathcal{H}\},
\end{equation}
which is 1-dimensional over $\mathbb{C}$. Note that taking $\tau$ to the boundary $\mathbb{P}^{1}(\mathbb{R})$ yields the 0-dimensional boundary component of $\operatorname{Gr}(L)$ that is associated with a line that is contained in our plane, namely the one spanned by $\tau w+z$ for finite $\tau\in\mathbb{R}$ and by $w$ for $\tau=\infty\in\mathbb{P}^{1}(\mathbb{R})$. We can summarize this in the following result.
\begin{prop} \label{boundR}
The boundary of $\operatorname{Gr}(L)$ inside the projective quadric from Proposition \ref{GVcomp} consists of a 1-dimensional copy of $\mathcal{H}$ for every 2-dimensional isotropic subspace $J_{\mathbb{R}}$ of $L_{\mathbb{R}}$ and a point for every isotropic line $I_{\mathbb{R}}$ in $L_{\mathbb{R}}$. The latter boundary components exists in any non-trivial case, i.e., wherever $n\geq1$, while the former ones are present only if $n\geq2$. All these boundary components are mutually disjoint, and the point associated with $I_{\mathbb{R}}$ lies in the closure of the 1-dimensional component associated with $J_{\mathbb{R}}$ if and only if $I_{\mathbb{R}} \subseteq J_{\mathbb{R}}$.
\end{prop}
The last assertion of Proposition \ref{boundR} is the incarnation of the statement, from \cite{[AMRT]} and others, that ``a boundary component of a boundary component is a boundary component'', in the orthogonal case. The reason for the notation $I_{\mathbb{R}}$ for an isotropic line and $J_{\mathbb{R}}$ for an isotropic plane will become clear in the next subsection.

\smallskip

Recalling the rational structure arising from the fact that $L_{\mathbb{R}}$ is the extension of scalars of the lattice $L$ to $\mathbb{R}$, the find that the Baily--Borel space $\mathcal{D}^{\mathrm{BB}}$ is obtained by adding to $\mathcal{D}=\operatorname{Gr}(L)$ those boundary components from Proposition \ref{boundR} that are associated with isotropic subspaces of $L_{\mathbb{R}}$ that arise as extensions of scalars from isotropic sublattices of $L$. In addition, the lattice $L$ has a \emph{dual lattice} \[L^{*}=\operatorname{Hom}(L,\mathbb{Z})\cong\big\{\nu \in L_{\mathbb{R}}\big|\;(\nu,L)\subseteq\mathbb{Z}\big\} \subseteq L_{\mathbb{R}},\] and we thus denote this subgroup of $L_{\mathbb{R}}$ also by $L^{*}$. It contains $L$ with finite index, and the quotient $\Delta_{L}=L^{*}/L$ is called the \emph{discriminant group} of $L$. It carries the $\mathbb{Q}/\mathbb{Z}$-valued quadratic form $\mu\mapsto\frac{\mu^{2}}{2}\in\mathbb{Q}/\mathbb{Z}$ and the non-degenerate $\mathbb{Q}/\mathbb{Z}$-valued bilinear form sending $\mu$ and $\nu$ to $(\mu,\nu)\in\mathbb{Q}/\mathbb{Z}$. The group $\operatorname{O}(L)$ of elements of $\operatorname{O}(L_{\mathbb{R}})$ that preserve $L$ also preserves $L^{*}$ and thus acts on $\Delta_{L}$. We put $\operatorname{SO}^{+}(L)=\operatorname{O}(L)\cap\operatorname{SO}^{+}(L_{\mathbb{R}})$, and define its  \emph{discriminant kernel} of $L$, sometimes called the \emph{stable special orthogonal group} of $L$, by
\begin{equation} \label{GammaLdef}
\Gamma_{L}=\ker\big(\operatorname{SO}^{+}(L)\to\operatorname{Aut}(\Delta_{L})\big).
\end{equation}
This is the arithmetic subgroup that we will consider throughout the rest of this paper.

\subsection{Parabolic Subgroups of Orthogonal Groups}

Given a boundary component of $\operatorname{Gr}(L)$, the parabolic subgroup of $\operatorname{SO}^{+}(L_{\mathbb{R}})$ stabilizing this component is the stabilizer of the corresponding isotropic subspace. There are (at most) two types of parabolic subgroups: Ones that are associated with isotropic lines in $L_{\mathbb{R}}$, and ones that arise from isotropic planes in $L_{\mathbb{R}}$. Having the rational structure, we work only with extensions of scalars of primitive isotropic sublattices of $L$. Hence every rational isotropic line in $L_{\mathbb{R}}$ is $I_{\mathbb{R}}$ for a primitive isotropic sublattice $I$ of rank 1 in $L$, for which $I=I_{\mathbb{R}} \cap L$, and we set
\begin{equation} \label{IKdef}
I_{L^{*}}=I_{\mathbb{R}} \cap L^{*},\quad I^{\perp}_{L}=I^{\perp} \cap L,\quad I^{\perp}_{L^{*}}=I^{\perp} \cap L^{*},\quad\mathrm{and}\quad K=I^{\perp}_{L}/I_{L}.
\end{equation}
Then $K_{\mathbb{R}}=I_{\mathbb{R}}^{\perp}/I_{\mathbb{R}}$ is Lorentzian (of signature $(n-1,1)$ in our convention), in which $K$ is an even lattice, and one verifies that \[\big\{\eta \in K_{\mathbb{R}}\big|\;(\eta,K)\subseteq\mathbb{Z}\big\} \cong K^{*}=\operatorname{Hom}(K,\mathbb{Z}) \cong I^{\perp}_{L^{*}}/I_{L^{*}}\] and the discriminant group of $K$ is given by \[\Delta_{K}=K^{*}/K=I^{\perp}_{L^{*}}/(I^{\perp}_{L}+I_{L^{*}})=L^{*}_{I}/(L+I_{L^{*}}),\] where $L^{*}_{I}$ stands for the subgroup
\begin{equation} \label{L*I}
L+I^{\perp}_{L^{*}}=\big\{\lambda \in L^{*}\big|\;\text{$\exists\nu \in L$ such that $(\lambda,\gamma)=(\nu,\gamma)$ for all $\gamma \in I$} \big\}
\end{equation}
of $L^{*}$. The lattice $L^{*}_{I}$ from \eqref{L*I} is the one denoted by $L_{0}'(N)$ in \cite{[Br1]}, and we define the projection
\begin{equation} \label{projDelta}
p_{K}:I^{\perp}_{L^{*}}/I^{\perp}_{L}=L^{*}_{I}/L\to\Delta_{K},\qquad p_{K}\mu=\mu+I_{L^{*}}/I.
\end{equation}
We also recall from \cite{[Br2]} the operators, that we denote by $\uparrow^{L}_{K}:\mathbb{C}[\Delta_{K}]\to\mathbb{C}[\Delta_{L}]$ and $\downarrow^{L}_{K}:\mathbb{C}[\Delta_{L}]\to\mathbb{C}[\Delta_{K}]$, that commute with the Weil representations, and are defined, for $\delta\in\Delta_{K}$ and $\mu\in\Delta_{L}$, by
\begin{equation} \label{arrowops}
\uparrow^{L}_{K}(\mathfrak{e}_{\delta})=\sum_{\mu \in p_{K}^{-1}(\delta)}\mathfrak{e}_{\mu}\qquad\mathrm{and}\qquad\downarrow^{L}_{K}(\mathfrak{e}_{\mu})=\begin{cases} \mathfrak{e}_{p_{K}(\mu)} & \mu \in L^{*}_{I}/L \\ 0 & \mathrm{otherwise}.\end{cases}
\end{equation}

We denote by $\mathcal{P}_{I,\mathbb{R}}$ the identity component of $\operatorname{St}_{\operatorname{SO}^{+}(L_{\mathbb{R}})}(I_{\mathbb{R}})$. A simple analysis of its structure yields (see, e.g., \cite{[Nak]}, \cite{[Bo2]}, or \cite{[Ze3]}, among others) the following result.
\begin{lem} \label{stdim1}
Consider the natural action of $\mathbb{R}_{+}^{\times}\times\operatorname{SO}^{+}(K_{\mathbb{R}})$ on $K_{\mathbb{R}}$, in which the former multiplier acts by ordinary scalar multiplication. Then there is a split exact sequence \[0\to(K_{\mathbb{R}},+)\to\mathcal{P}_{I,\mathbb{R}}\to\mathbb{R}_{+}^{\times}\times\operatorname{SO}^{+}(K_{\mathbb{R}})\to1.\]
\end{lem}

\smallskip

The rational isotropic planes in $L_{\mathbb{R}}$, which may exist only when $n\geq2$, are all of the form $J_{\mathbb{R}}$, where $J$ is a primitive isotropic sublattice of rank 2 in $L$. Then $J=J_{\mathbb{R}} \cap L$, and as in \eqref{IKdef} we define
\begin{equation} \label{JDdef}
J_{L^{*}}=J_{\mathbb{Q}} \cap L^{*},\quad J^{\perp}_{L}=J^{\perp} \cap L,\quad J^{\perp}_{L^{*}}=J^{\perp} \cap L^{*},\quad\mathrm{and}\quad D=J^{\perp}_{L}/J.
\end{equation}
Here $D$ is an even lattice inside the positive definite $n-2$-dimensional space $D_{\mathbb{R}}=J_{\mathbb{R}}^{\perp}/J_{\mathbb{R}}$, and we have descriptions of the dual lattice $D^{*} \subseteq D_{\mathbb{R}}$, the discriminant group $\Delta_{D}$, the lattice $L^{*}_{J}$, the projection $p_{D}$, and the operators $\uparrow^{L}_{D}$ and $\downarrow^{L}_{D}$ that are similar to those from $I$ and $K$ above, as in, e.g., \eqref{L*I}, \eqref{projDelta}, and \eqref{arrowops}. Since $D_{\mathbb{R}}$ is definite, the stabilizer $\operatorname{St}_{\operatorname{SO}^{+}(L_{\mathbb{R}})}(J_{\mathbb{R}})$ is connected, and we denote it by $\mathcal{P}_{J,\mathbb{R}}$. For presenting its structure we recall the following definition.
\begin{defn} \label{Heisdef}
The \emph{Heisenberg group} $\operatorname{H}(D_{\mathbb{R}})$ associated with $D_{\mathbb{R}}$ is the group $D_{\mathbb{R}} \times D_{\mathbb{R}}\times\mathbb{R}$ with the product rule \[(\alpha,\beta,t)\cdot\big(\gamma,\delta,s)=(\alpha+\gamma,\beta+\delta,t+s+(\alpha,\delta)-(\beta,\gamma)\big).\] It is a non-trivial central extension of $D_{\mathbb{R}} \times D_{\mathbb{R}}$ by $\operatorname{Z}\big(\operatorname{H}(D_{\mathbb{R}})\big)=\{(0,0,t)|\;t\in\mathbb{R}\}\cong(\mathbb{R},+)$, as expressed by the non-split short exact sequence \[0\to\operatorname{Z}\big(\operatorname{H}(D_{\mathbb{R}})\big)\to\operatorname{H}(D_{\mathbb{R}}) \to D_{\mathbb{R}} \oplus D_{\mathbb{R}}\to0.\]
\end{defn}
The group $\operatorname{GL}_{2}^{+}(\mathbb{R})\times\operatorname{SO}(D_{\mathbb{R}})$ operates on $\operatorname{H}(D_{\mathbb{R}})$ via the formula
\begin{equation} \label{actonHeis}
[M,\gamma]:(\alpha,\beta,t)\mapsto(a\gamma\alpha+b\gamma\beta,c\gamma\alpha+d\gamma\beta,\det M \cdot t),
\end{equation}
where $M=\big(\begin{smallmatrix} a & b \\ c & d\end{smallmatrix}\big)$. The action is compatible with the group structure on $\operatorname{H}(D_{\mathbb{R}})$. The structure of $\mathcal{P}_{J,\mathbb{R}}$ is now given (see, e.g., \cite{[K3]} or \cite{[Ze3]}) in the following lemma.

\begin{lem} \label{stdim2}
After choosing an oriented basis $(z,w)$ for $J_{\mathbb{R}}$, the action of $\operatorname{GL}_{2}^{+}(\mathbb{R})\times\operatorname{SO}(D_{\mathbb{R}})$ on $\operatorname{H}(D_{\mathbb{R}})$ in \eqref{actonHeis} defines $\mathcal{P}_{J,\mathbb{R}}$ via the split short exact sequence \[1\to\operatorname{H}(D_{\mathbb{R}})\to\mathcal{P}_{J,\mathbb{R}}\to\operatorname{GL}_{2}^{+}(\mathbb{R})\times\operatorname{SO}(D_{\mathbb{R}})\to1.\]
\end{lem}
In correspondence with the general structure theorem from, e.g., Lemma 1.1 and Propositions 1.4 and 1.6 of \cite{[Ze3]}, the choice of the basis $(z,w)$ for $J_{\mathbb{R}}$ yields the isomorphisms \[\operatorname{GL}(J_{\mathbb{R}})\cong\operatorname{GL}_{2}(\mathbb{R}),\quad\operatorname{Hom}_{\mathbb{R}}(D_{\mathbb{R}},J_{\mathbb{R}}) \cong D_{\mathbb{R}} \oplus D_{\mathbb{R}},\quad\text{and}\quad\operatorname{Hom}_{\mathbb{R}}^{as}(J_{\mathbb{R}}^{*},J_{\mathbb{R}})\cong\mathbb{R},\] where for the middle one we also used the self-duality of $D_{\mathbb{R}}$.

\smallskip

In order to put these results in the framework of the theory of toroidal compactifications, we shall henceforth write $\Xi$ for a 0-dimensional boundary component and $\Upsilon$ for a 1-dimensional boundary component. The symmetric space $\mathcal{D}$ is $\operatorname{Gr}(L)$, and its rational completion is \[\mathcal{D}^{\mathrm{BB}}=\operatorname{Gr}(L)\cup\Big(\bigcup_{I \subseteq L}\Xi_{I}\Big)\cup\Big(\bigcup_{J \subseteq L}\Upsilon_{J}\Big),\] where $I$ (resp. $J$) will always mean a primitive isotropic sublattice of $L$ of rank 1 (resp. 2). Note that it is possible that the union over $J$ here, and even the union over $I$, may be empty even in cases where Proposition \ref{boundR} does not imply its emptiness. On the other hand, Meyer's Theorem (see, e.g., Corollary 2 on page 43 of \cite{[Se]}) assures us that any indefinite lattice of dimension at least 5 contains a non-zero isotropic vector. It follows that if $n\geq3$ then rational 0-dimensional boundary components always exist, and when $n\geq5$ we also have the existence of 1-dimensional such components.

Now, if $\Xi=\Xi_{I}$ is the 0-dimensional boundary component of $\mathcal{D}^{\mathrm{BB}}\subseteq\overline{\operatorname{Gr}(L)}$ that is associated with the rank 1 lattice $I \subseteq L$, then $\mathcal{P}(\Xi)$ is the group $\mathcal{P}_{I,\mathbb{R}}$. We relate its structure in Lemma \ref{stdim1} with the general structure from Proposition \ref{genparstruc} as follows.
\begin{prop} \label{paradim0}
The unipotent radical $\mathcal{W}(\Xi)$ of $\mathcal{P}(\Xi)=\mathcal{P}_{I,\mathbb{R}}$, as well as its center $\mathcal{U}(\Xi)$, is the subgroup $(K_{\mathbb{R}},+)$ from Lemma \ref{stdim1}. The quotient $\mathcal{V}(\Xi)$ is trivial, and so are the group $G_{h}(\Xi)=\operatorname{Aut}(\Xi)$ and the compact group $M(\Xi)$. The group $G_{\ell}(\Xi)$ is the full quotient $\mathbb{R}_{+}^{\times}\times\operatorname{SO}^{+}(K_{\mathbb{R}})$, namely the full Levi subgroup.
\end{prop}
Indeed, $\Xi$ is a point, and the Levi subgroup from Lemma \ref{stdim1} has no compact factors, yielding the triviality of $G_{h}(\Xi)$ and $M(\Xi)$.

For the explicit description of the fibrations from Proposition \ref{DoverXi}, choose an orientation on $I$, and let $z$ be the oriented generator of $I$. As presented in Section 13 of \cite{[Bo1]}, Section 3.2 of \cite{[Br1]}, Section 2.2 of \cite{[Ze1]}, or Section 1 of \cite{[Ze2]}, the condition of pairing to 1 with $z$ produces a section of the restriction of the projection from \eqref{fiboverGVR} to the component $P^{+}$. The orientation on elements $v\in\operatorname{Gr}(L)$ and on $I$ given by the choices of $P^{+}$ and $z$ determine a cone $C$ of negative norm vectors in the Lorentzian space $K_{\mathbb{R}}$, such that if \[K_{\mathbb{R}}^{1}=\big\{\xi \in L_{\mathbb{R}}\big|\;(\xi,z)=1\big\}\big/I_{\mathbb{R}}\] is the affine model of $K_{\mathbb{R}}$ arising from $z$, then our section defines a biholomorphic map
\begin{equation} \label{GVRKRiC}
\operatorname{Gr}(L) \to K_{\mathbb{R}}^{1}+iC\stackrel{open}{\subseteq}K_{\mathbb{R}}^{1}+iK_{\mathbb{R}}=K_{\mathbb{C}}^{1}=\big\{\xi \in L_{\mathbb{C}}\big|\;(\xi,z)=1\big\}\big/I_{\mathbb{C}}.
\end{equation}
The inverse of the map from \eqref{GVRKRiC}, which implies that it is a biholomorphism, is easily determined by the fact that different lifts to $L_{\mathbb{C}}$ of an element of $K_{\mathbb{C}}^{1}$, which has a non-zero pairing with the isotropic vector $z$, have different norms, and hence only one such lift is also isotropic. The cone $C$ and the complementary cone $-C$ are self-adjoint in the terminology of \cite{[AMRT]} or \cite{[Nam]}, and the model $K_{\mathbb{R}}^{1}+iC$ for $\operatorname{Gr}(L)$ from \eqref{GVRKRiC} is a \emph{tube domain}. The explicit description of the objects from Proposition \ref{DoverXi} is now as follows.
\begin{prop} \label{struc0dim}
The basis $\Xi$ and the fiber bundle $\widetilde{D}(\Xi)$ are trivial, and $\mathcal{D}=\operatorname{Gr}(L)$ is isomorphic (as the fiber of a map to a point) to the affine tube domain $K_{\mathbb{R}}^{1}+iC$ from \eqref{GVRKRiC}. The map $\Phi$ takes the element $v\in\operatorname{Gr}(L)$, whose image in $P^{+}$ under the section is $Z_{L}$, to the projection $Y \in C$ of $Y_{L}=\Im Z_{L} \in I_{\mathbb{R}}^{\perp}$ in $K_{\mathbb{R}}$, and the cone $\Omega(\Xi)$ is $C$.
\end{prop}
In addition, $\mathcal{U}(\Xi)=K_{\mathbb{R}}$ acts additively on the real part of $\operatorname{Gr}(L) \cong K_{\mathbb{R}}^{1}+iC$ from \eqref{GVRKRiC}, and the quotient $\mathbb{R}_{+}^{\times}\times\operatorname{SO}^{+}(K_{\mathbb{R}})$ acts on its symmetric space, the cone $C$, via the map $\Phi$ described in Proposition \ref{struc0dim}, as Proposition \ref{subgpact} predicts.

We also need the intersections from \eqref{grpsindZ}, with $\Gamma$ being $\Gamma_{L}$ from \eqref{GammaLdef}. The following result, in which $\Gamma_{K}$ is defined analogously to $\Gamma_{L}$, is proved in Corollary~3.6 of \cite{[Ze3]}, though it was known earlier.
\begin{prop} \label{dim1}
The intersection $\mathcal{P}_{\mathbb{Z}}(\Xi)=\mathcal{P}_{I,\mathbb{R}}\cap\Gamma_{L}$ lies in a split short exact sequence \[0\to(K,+)\to\mathcal{P}_{I,\mathbb{R}}\cap\Gamma_{L}\to\Gamma_{K}\to1.\] Hence $\mathcal{W}_{\mathbb{Z}}(\Xi)=\mathcal{U}_{\mathbb{Z}}(\Xi)=K$, the group $\mathcal{P}_{\mathbb{Z}}(\Xi)$ is the semi-direct product, and its image $\Gamma_{\ell}(\Xi)$ in $G_{\ell}(\Xi)$ is $\Gamma_{K}$.
\end{prop}

As the action of $\mathcal{U}_{\mathbb{Z}}(\Xi)=K$ on  $\operatorname{Gr}(L) \cong K_{\mathbb{R}}^{1}+iC$ from \eqref{GVRKRiC} is additive on the real part, we get that as the fiber over a point from Proposition \ref{struc0dim} we have
\begin{equation} \label{GVRmodK}
K \backslash \operatorname{Gr}(L)\cong(K_{\mathbb{R}}^{1}/K)+iC\stackrel{open}{\subseteq}(K_{\mathbb{R}}^{1}/K)+iK_{\mathbb{R}}=K_{\mathbb{C}}^{1}/K=T_{K}^{1}(\mathbb{C}),
\end{equation}
which is indeed a principal homogenous space over $T_{K}(\mathbb{C})=T_{\mathcal{U}_{\mathbb{Z}}(\Xi)}(\mathbb{C})$, as expected from \eqref{affTUZ}. In \cite{[Bo1]}, \cite{[Br1]}, \cite{[Ze1]}, \cite{[Ze2]}, and others, the affine tube domain from \eqref{GVRKRiC} is replaced, via the choice of $\zeta \in L^{*}$ with $(\zeta,z)=1$, by a subset of $K_{\mathbb{C}}$ itself, and modulo $K$ one obtains a subset of $T_{K}(\mathbb{C})$ in \eqref{GVRmodK}. However, this is not necessary for most applications, and we stay with our affine models, since they are more canonical. We summarize:
\begin{cor} \label{D/UZXi}
The set $\mathcal{U}_{\mathbb{Z}}(\Xi)\backslash\mathcal{D}=K \backslash \operatorname{Gr}(L)$ is isomorphic, as the fiber over the trivial space $\widetilde{\mathcal{D}}(\Xi)$, to an open subset of a principal homogenous space over $T_{K}(\mathbb{C})$, with imaginary part in the direction of $C$.
\end{cor}

\smallskip

We now compare the structure of the groups and fibrations that are associated with a 1-dimensional boundary component to the components from Section \ref{TorComp}. For a rank 2 primitive isotropic sublattice $J \subseteq L$, we denote the 1-dimensional boundary component $\Upsilon_{J}$ of $\mathcal{D}^{\mathrm{BB}}\subseteq\overline{\operatorname{Gr}(L)}$ by $\Upsilon$, and then $\mathcal{P}(\Upsilon)$ is the group $\mathcal{P}_{J,\mathbb{R}}$ considered in Lemma \ref{stdim2}. For putting this lemma in the framework of Proposition \ref{genparstruc}, recall that the action of $\operatorname{GL}_{2}^{+}(\mathbb{R})$ on $\mathcal{H}$ and its factor of automorphy are defined as
\begin{equation} \label{GL2actH}
M=\begin{pmatrix} a & b \\ c & d \end{pmatrix}\in\operatorname{GL}_{2}^{+}(\mathbb{R}):\tau\in\mathcal{H} \mapsto M\tau=\frac{a\tau+b}{c\tau+d},\qquad\mathrm{and}\qquad j(M,\tau)=c\tau+d
\end{equation}
respectively, and that $\operatorname{GL}_{2}^{+}(\mathbb{R})$ is the direct product of the connected component $\mathbb{R}_{+}^{\times}$ of $\operatorname{Z}\big(\operatorname{GL}_{2}^{+}(\mathbb{R})\big)$ with $\operatorname{SL}_{2}(\mathbb{R})$. The decomposition of $\mathcal{P}(\Upsilon)$ is now as follows.
\begin{prop} \label{paradim1}
The unipotent radical $\mathcal{W}(\Upsilon)$ of the group $\mathcal{P}(\Upsilon)=\mathcal{P}_{J,\mathbb{R}}$ is the Heisenberg group $\operatorname{H}(D_{\mathbb{R}})$, its center $\mathcal{U}(\Upsilon)$ is $\operatorname{Z}\big(\operatorname{H}(D_{\mathbb{R}})\big)=\mathbb{R}$, and the quotient $\mathcal{V}(\Upsilon)$ is $D_{\mathbb{R}} \oplus D_{\mathbb{R}}$. Inside the Levi subgroup, $G_{h}(\Upsilon)$ is $\operatorname{SL}_{2}(\mathbb{R})$, $G_{\ell}(\Upsilon)$ is the group of $\mathbb{R}_{+}^{\times}$ of scalars, and $M(\Upsilon)$ is the compact group $\operatorname{SO}(D_{\mathbb{R}})$.
\end{prop}

Take now, following Section 1 of \cite{[K3]}, an oriented $\mathbb{Z}$-basis $(z,w)$ for $J$ from \eqref{JDdef}, and lift $\operatorname{Gr}(L)$ into $P^{+}$ via the section of the map from \eqref{fiboverGVR} that is associated with $z$. Assume that $v\in\operatorname{Gr}(L)$ is taken to $Z_{L} \in P^{+}$, and then the fact that $Y_{L}^{2}<0$ for $Y_{L}=\Im Z_{L}$ implies that $Y_{L} \not\in J_{\mathbb{R}}^{\perp}$ and hence $(Y_{L},w)\neq0$. The orientation then implies that
\begin{equation} \label{secwithz}
x+iy=\tau=(Z_{L},w)\in\mathcal{H}
\end{equation}
(note that in spite of the similar notation, this is \emph{not} the coordinate for the cusp $\Upsilon$ itself, appearing in \eqref{coor1dcusp}). We now define, similarly to $K_{\mathbb{R}}^{1}$ and $K_{\mathbb{C}}^{1}$ from \eqref{GVRKRiC}, the affine spaces \[D_{\mathbb{R}}^{1,x}=\big\{\xi \in L_{\mathbb{R}}\big|\;(\xi,z)=1,\ (\xi,w)=x\big\}\big/J_{\mathbb{R}},\quad D_{\mathbb{R}}^{0,y}=\big\{\xi \in L_{\mathbb{R}}\big|\;(\xi,z)=0,\ (\xi,w)=y\big\}\big/J_{\mathbb{R}},\] and the complexification \[D_{\mathbb{C}}^{1,\tau}=D_{\mathbb{R}}^{1,x}+iD_{\mathbb{R}}^{0,y}=\big\{\xi \in L_{\mathbb{C}}\big|\;(\xi,z)=1,\ (\xi,w)=\tau\big\}\big/J_{\mathbb{C}}.\] These are principal homogenous spaces over $D_{\mathbb{R}}$ and an affine model of $D_{\mathbb{C}}$ respectively. The subgroups from Proposition \ref{paradim1} are given in \cite{[Fi]} in terms of explicit matrices, under the choice of a representing matrix for the quadratic form on $L$, but the more canonical description of the fibrations from Proposition \ref{DoverXi} uses the following presentation.
\begin{lem} \label{actU}
With the choice of the basis $(z,w)$ for $J \subseteq J_{\mathbb{R}}$, the action of $\mathcal{W}(\Upsilon)$ on $L_{\mathbb{R}}$ is as follows. First, the action of the element $\alpha\in\mathbb{R}\cong\mathcal{U}(\Upsilon)$ on $L_{\mathbb{R}}$ is defined by \[\alpha:\lambda\mapsto\lambda-\alpha(\lambda,z)w+\alpha(\lambda,w)z.\] Consider now a pair $(\mu,\nu) \in D_{\mathbb{R}} \oplus D_{\mathbb{R}}=\mathcal{V}(\Upsilon)$, and let $\eta$ denote an element of $J_{\mathbb{R}}^{\perp}$ and $\lambda+J_{\mathbb{R}}$ denote the image of $\lambda \in L_{\mathbb{R}}$ in $L_{\mathbb{R}}/J_{\mathbb{R}}$. Then the action of $(\mu,\nu)$ is defined by the formulas \[\eta\mapsto\eta-(\eta,\mu)w-(\eta,\nu)z\qquad\mathrm{and}\qquad\lambda+J_{\mathbb{R}}\mapsto\lambda+(\lambda,w)\mu+(\lambda,z)\nu+J_{\mathbb{R}}.\] Complexifying, and taking $Z_{L}$ and $\tau$ are as in \eqref{secwithz}, we deduce that \[\alpha:Z_{L}+\mathbb{C}z \mapsto Z_{L}-\alpha w+\mathbb{C}z,\qquad\mathrm{and}\qquad(\mu,\nu):Z_{L}+J_{\mathbb{R}} \mapsto Z_{L}+\tau\mu+\nu+J_{\mathbb{R}}.\]
\end{lem}

\begin{proof}
In general, the unipotent radical $\mathcal{W}(\Upsilon)$ is defined as the kernel of the action of $\mathcal{P}(\Upsilon)=\mathcal{P}_{J,\mathbb{R}}$ on $J_{\mathbb{R}}$ via $\operatorname{GL}^{+}(J_{\mathbb{R}})\cong\operatorname{GL}_{2}^{+}(\mathbb{R})$ and on $D_{\mathbb{R}}$ via $\operatorname{SO}(D_{\mathbb{R}})$ (see, e.g., \cite{[Ze3]}). Moreover, as $\dim J_{\mathbb{R}}>1$, the center $\mathcal{U}(\Upsilon)$ is the subgroup acting trivially on $J_{\mathbb{R}}^{\perp}$, and in dimension 2 the asserted action is the only orthogonal one fixing $J_{\mathbb{R}}^{\perp}$ pointwise (this is in correspondence with the space $\operatorname{Hom}_{\mathbb{R}}^{as}(J_{\mathbb{R}}^{*},J_{\mathbb{R}})$ of \cite{[Ze3]}, or equivalently $\bigwedge^{2}J_{\mathbb{R}}$ in the notation of \cite{[L]}, being 1-dimensional). It follows that the quotient $\mathcal{V}(\Upsilon)$ has a well-defined action on elements of $J_{\mathbb{R}}^{\perp}$, which is given by subtracting elements of $J_{\mathbb{R}}$, but on general elements of $L_{\mathbb{R}}$ the action is only defined modulo $J_{\mathbb{R}}$ (by subtracting elements of the quotient $D_{\mathbb{R}}$). With the appropriate normalization, the orthogonality condition yields again the required formulas. For the last assertion we just consider the image of $\alpha$ modulo $\mathbb{C}z$, and substitute the values $(Z_{L},z)=1$ and $(Z_{L},w)=\tau$ from \eqref{secwithz}. This proves the lemma.
\end{proof}

The structure of $\mathcal{D}=\operatorname{Gr}(L)$ that is associated with $\Upsilon$ in Proposition \ref{DoverXi} can now be described.
\begin{prop} \label{struc1dim}
The affine vector bundle $\widetilde{\mathcal{D}}(\Upsilon)$ from \eqref{vecfib} is $\bigcup_{\tau\in\mathcal{H}}D_{\mathbb{C}}^{1,\tau}$, with $D_{\mathbb{C}}^{1,\tau}$ being the fiber over $\tau\in\Upsilon\cong\mathcal{H}$. The fiber of $\mathcal{D}$ over any point in $\widetilde{\mathcal{D}}(\Upsilon)$, as in \eqref{holfib}, is a translated upper half-plane. The cone $\Omega(\Upsilon)$ is the ray $\mathbb{R}_{+}$ inside $\mathcal{U}(\Upsilon)=\mathbb{R}$, and if $v\in\operatorname{Gr}(L)$ is represented by $Z_{L} \in P^{+}$ and $\tau$ is as in \eqref{secwithz} then $\Phi(v)=-\frac{(\Im Z_{L})^{2}}{2y}$.
\end{prop}

\begin{proof}
Lemma \ref{actU} shows that the fiber $D_{\mathbb{C}}^{1,\tau}$ over $\tau\in\Upsilon$ is a principal homogenous space over $\mathcal{V}(\Upsilon)$, with the complex structure varying with $\tau$, as in \eqref{vecfib}. In addition, we have seen that if $v\in\operatorname{Gr}(L)$ is sent to $Z_{L} \in P^{+}$ then its image on $\widetilde{\mathcal{D}}(\Upsilon)$ is just $Z_{L}+J_{\mathbb{C}} \in D_{\mathbb{C}}^{1,\tau}$, for $\tau$ from \eqref{secwithz}. Hence two images of our section that map to the same point in $D_{\mathbb{C}}^{1,\tau}$ can only differ by complex multiples of $w$ and $z$. But inverting the map from \eqref{GVRKRiC} shows that the condition $Z_{L}^{2}=0$ determines the required multiple of $z$ once that of $w$ is chosen, so that it remains to see which multiples of $w$ satisfy the condition from \eqref{bunGV}. Write $Z_{L}$ as $X_{L}+iY_{L}$ as above, and as $Y_{L} \perp z$ in the image of the section, the value of $Y_{L}^{2}$ is evaluated using only the multiple of $w$, and subtracting $tw$ from $Y_{L}$ subtracts $2yt$ from $Y_{L}^{2}$. Therefore $t$ must be bounded from below for obtaining section images from $\operatorname{Gr}(L)$, and there is no restriction on the multiple of $w$ subtracted from $X_{L}$. The associated fiber is therefore an affine translated copy of $\mathcal{H}$, and in particular it is non-empty, so that the map from $\operatorname{Gr}(L)$ to $\bigcup_{\tau\in\mathcal{H}}D_{\mathbb{C}}^{1,\tau}$ is indeed surjective. As $\mathcal{U}(\Upsilon)$ is 1-dimensional, and the oriented cone $\Omega(\Upsilon)$ can only be the directed ray in $\mathcal{U}(\Upsilon)$, this is indeed as described in \eqref{holfib}.

Now, Lemma \ref{actU} also shows that $\mathcal{U}(\Upsilon)$ acts transitively on the real parts of the fibers of from \eqref{holfib}, and the asserted $\Phi$-image $-\frac{Y_{L}^{2}}{2y}$ was seen to be a translate of the imaginary part of these fibers. Moreover, as this $\Phi$-image depends only on the imaginary part of $Z_{L}$, the presentation from \eqref{Phifib} also follows. In addition, our $\Phi$-image lies, by definition, in the directed cone $\Omega(\Upsilon)=\mathbb{R}_{+}^{\times}$ inside $\mathcal{U}(\Upsilon)=\mathbb{R}$, and since Lemma \ref{actU} shows that both $\tau$ and $z$ are invariant under $\mathcal{W}(\Upsilon)$, so is the value of $\Phi$. Now, the action of the Levi quotient on $\widetilde{\mathcal{D}}(\Upsilon)$ over $\Upsilon$ and on $\Omega(\Upsilon)$ via $\Phi$ is given in \eqref{actonHeis} via the interpretations from Proposition \ref{DoverXi}, and if $M\in\operatorname{GL}_{2}^{+}(\mathbb{R})=G_{h}(\Upsilon) \times G_{\ell}(\Upsilon)$ and $Z_{L}$ is in the image of the section and has projection in $D_{\mathbb{C}}^{1,\tau}$ for $\tau$ from \eqref{secwithz}, then we have to divide the result by the factor of automorpy $j(M,\tau)$ from \eqref{GL2actH} in order to remain in the image of the section. It follows that the action of $G_{h}(\Upsilon)=\operatorname{SL}_{2}(\mathbb{R})$ and $M(\Upsilon)=\operatorname{SO}(D_{\mathbb{R}})$ are trivial on $\Omega(\Upsilon)$, and that the group $G_{\ell}(\Upsilon)=\mathbb{R}_{+}^{\times}$ from Proposition \ref{paradim1} is indeed the automorphism group of $\Omega(\Upsilon)$, with trivial stabilizers. On the other hand, on $\Upsilon\cong\mathcal{H}$ the groups $G_{\ell}(\Upsilon)$ and $M(\Upsilon)$ operate trivially, and the action of $G_{h}(\Upsilon)$ is the natural one from \eqref{GL2actH}. This proves the proposition.
\end{proof}
Note that \cite{[K3]} and \cite{[Fi]} essentially use a complement for $J_{\mathbb{R}}^{\perp}$ in $L_{\mathbb{R}}$, describing $Z_{L}$ as in \eqref{coorJ} below. However, we prefer the affine presentation from Proposition \ref{struc1dim}, since it is more canonical and does not depend on the complement.

As a more geometric interpretation of the space $\widetilde{\mathcal{D}}(\Upsilon)$, relating to Remark \ref{VHS}, we recall the following definition.
\begin{defn} \label{ShimVHS}
Over $\mathcal{H}$ lies \emph{Shimura's variation of Hodge structure} $\mathcal{V}_{1}$, which has dimension 2. Its fiber is $\mathbb{C}^{2}$, and over $\tau\in\mathcal{H}$ it is defined by the vector $\binom{\tau}{1}$ having Hodge weight $(1,0)$, while $\binom{\overline{\tau}}{1}$ is of weight $(0,1)$. It is $\operatorname{SL}_{2}(\mathbb{R})$-equivariant with respect to the standard representation of that group on $\mathbb{C}^{2}$, and has pure weight 1. If $\pi:\mathcal{E}\to\mathcal{H}$ is the universal elliptic curve over $\mathcal{H}$ then we have $\mathcal{V}_{1} \cong R^{1}\pi_{*}\mathbb{C}$, and the integral structure on the local system is based on the subgroup $\mathbb{Z}^{2}\subseteq\mathbb{C}^{2}$.
\end{defn}
In terms of Definition \ref{ShimVHS}, the short exact sequence form \eqref{vecfib} describes the space $\widetilde{\mathcal{D}}(\Upsilon)$, via Proposition \ref{struc1dim}, as an affine model of the vector bundle $\mathcal{V}_{1} \otimes D$ over $\Upsilon\cong\mathcal{H}$.

\smallskip

As \cite{[Ze3]} shows, the structure of the intersection $\mathcal{P}_{\mathbb{Z}}(\Upsilon)=\mathcal{P}_{J,\mathbb{R}}\cap\Gamma_{L}$ can be rather delicate in general. Indeed, \cite{[Fi]} makes the rather restrictive assumption that $L$ splits two hyperbolic planes over $\mathbb{Z}$, a case in which most of the difficulties disappear. We shall describe this group under some assumptions (less restrictive than those of \cite{[Fi]}), though the main results are valid for arbitrary lattices.

Recall that for any primitive sublattice $M \subseteq L$ there is a \emph{unimodular complement}, namely a subgroup of $L^{*}$ whose direct sum with $M^{\perp} \cap L^{*}$ yields the full lattice $L^{*}$ (the name arises from the fact that such a complement becomes naturally isomorphic to $\operatorname{Hom}(M,\mathbb{Z})$, since $M$ is primitive). From now until Remark \ref{genlat} below we shall make the following hypothesis on the lattice $J$ from \eqref{JDdef}.
\begin{hypo} \label{assonJ}
The lattice $J$ has a unimodular complement $\tilde{J}$ that is isotropic and whose intersection with the lattice $L^{*}_{J}$ defined as in \eqref{L*I} is contained in $L$, and the pairing $(\lambda,\mu)$ is even for every $\lambda$ and $\mu$ in the definite lattice $D=J^{\perp}_{L}/J$ from \eqref{JDdef}.
\end{hypo}
Note that our assumption on the pairing in $D$ in Hypothesis \ref{assonJ} is stronger than $D$ being an even lattice (which is always the case when $L$ is even), and is required for the assertion about $\operatorname{H}(D)$ in Theorem \ref{PZUpsilon} below (see Equation (14) of \cite{[Ze3]}). It is always satisfied when $L$ is the rescaling by 2 of an integral lattice. The properties of the complement $\tilde{J}$ from Hypothesis \ref{assonJ} are equivalent to the vanishing of the maps $\alpha$ and $\iota$ from \cite{[Ze3]}. Equation (16) of \cite{[Ze3]} then shows, after substituting $\alpha=0$ there, that \[L=J\oplus(\tilde{J} \cap L^{*}_{J})\oplus\big[(J\oplus\tilde{J})^{\perp} \cap L\big]\quad\mathrm{and}\quad L^{*}=J_{L^{*}}\oplus\tilde{J}\oplus\big[(J\oplus\tilde{J})^{\perp} \cap L^{*}\big],\] with $(J\oplus\tilde{J})^{\perp} \cap L$ (resp. $(J\oplus\tilde{J})^{\perp} \cap L^{*}$) mapping isomorphically onto $D$ (resp. $D^{*}$). It is easy to see that a complement as in Hypothesis \ref{assonJ} always exists when $L$ is the rescaling of a maximal lattice, since a maximal lattice of Witt rank 2 always splits two hyperbolic planes, but there are other examples of lattices where such complements exist. Note that a similar condition for the rank 1 lattices $I$ (which does \emph{not} follow from that on the lattices $J$, because a complement $\tilde{I}$ and the lattice $L^{*}_{I}$ from \eqref{L*I} are larger than those associated with $J$) is not required, since Proposition \ref{dim1} holds without any assumption on $I$.

As the basis $(z,w)$ that we took for $J_{\mathbb{R}}$ spans $J$ over $\mathbb{Z}$, the resulting isomorphism between $\operatorname{GL}_{2}^{+}(\mathbb{R})$ and $\operatorname{GL}^+(J_{\mathbb{R}})$, appearing implicitly in Lemma \ref{stdim2}, identifies $\operatorname{SL}_{2}(\mathbb{Z})$ with $\operatorname{SL}(J)$. We define $\Gamma_{J}$ to be the congruence subgroup of $\operatorname{SL}_{2}(\mathbb{Z})$ that corresponds under this isomorphism to \[\operatorname{SL}(J_{L^{*}},J)=\big\{M\in\operatorname{SL}(J)\big|\;M\lambda-\lambda \in J\mathrm{\ for\ every\ }\lambda \in J_{L^{*}}\big\}.\] We also define the \emph{integral Heisenberg group} $\operatorname{H}(D)$ to be $D \times D\times\mathbb{Z}$ with the product rule from Definition \ref{Heisdef}, and we get from \cite{[Ze3]} the following result.
\begin{thm} \label{PZUpsilon}
The group $\operatorname{H}(D)$ is embedded in $\operatorname{H}(D_{\mathbb{R}})$ as the coordinates $D \times D\times\mathbb{Z}$, and one has the short exact sequence \[1\to\operatorname{H}(D)\to\mathcal{P}_{J,\mathbb{R}}\cap\Gamma_{L}\to\Gamma_{J}\times\Gamma_{D}\to1,\] which is split over $\mathbb{Z}$, with the action from \eqref{actonHeis}.
\end{thm}
The groups with index $\mathbb{Z}$ here take the following form, in which we define $u_{\Upsilon}$ to be the unique oriented primitive element of the cyclic group $\mathcal{U}_{\mathbb{Z}}(\Upsilon)$ that lies in the cone $\Omega(\Upsilon)$.
\begin{prop} \label{dim2}
The integral unipotent radical $\mathcal{W}_{\mathbb{Z}}(\Upsilon)$ of $\mathcal{P}_{\mathbb{Z}}(\Upsilon)=\mathcal{P}_{J,\mathbb{R}}\cap\Gamma_{L}$ is $\operatorname{H}(D)$, its center $\mathcal{U}_{\mathbb{Z}}(\Upsilon)$ is cyclic and generated by $u_{\Upsilon}$, and the quotient $\mathcal{U}_{\mathbb{Z}}(\Upsilon)$ is $D \oplus D$. In addition, the image $\Gamma_{\ell}(\Upsilon)$ of $\mathcal{P}_{\mathbb{Z}}(\Upsilon)$ in $G_{\ell}(\Upsilon)$ is trivial, in $G_{h}(\Upsilon)$ is it $\Gamma_{J}$, and in $M(\Upsilon)$ it is the finite group $\Gamma_{D}$.
\end{prop}
The cyclic group $\mathcal{U}_{\mathbb{Z}}(\Upsilon)$ from Proposition \ref{dim2}, oriented via $\Omega(\Upsilon)$, is $\operatorname{Hom}^{as}(J^{*},J)$ in the notation of \cite{[Ze3]} and $\bigwedge^{2}J$ in the notation of \cite{[L]}. The triviality of $\Gamma_{\ell}(\Upsilon)$ implies that there is no question of stabilizers in the fan $\Sigma(\Upsilon)$. Recall that if $\Gamma_{L}$ is neat then neither $\Gamma_{L}$ nor the group $\Gamma_{K}$ from Proposition \ref{dim1} nor $\Gamma_{J}$ have non-trivial elements of finite order, and the finite group $\Gamma_{D}$ is trivial. Dividing by the groups from Proposition \ref{dim2} thus yields, by Lemma \ref{actU} and Proposition \ref{struc1dim}, the following result.
\begin{cor} \label{D/UZUpsilon}
The fibers of $\mathcal{U}_{\mathbb{Z}}(\Upsilon)\backslash\mathcal{D}=\mathbb{Z}u_{\Upsilon} \backslash \operatorname{Gr}(L)$ over the space $\widetilde{\mathcal{D}}(\Upsilon)$ from Proposition \ref{struc1dim} are open subsets of affine models of the 1-dimensional torus $T_{\mathcal{U}_{\mathbb{Z}}(\Upsilon)}(\mathbb{C})$. The space $\mathcal{W}_{\mathbb{Z}}(\Upsilon) \backslash \operatorname{Gr}(L)$ is a similar bundle over the universal family $\mathcal{E} \otimes D$ over $\mathcal{H}$ from Definition \ref{ShimVHS}. If $\Gamma_{L}$ is neat then the quotient $\mathcal{P}_{\mathbb{Z}}(\Upsilon) \backslash \operatorname{Gr}(L)$ is also a bundle with the same fibers, but now over the open Kuga--Sato variety $W_{J}$ defined by $\mathcal{E} \otimes D$ over $\Gamma_{J}\backslash\mathcal{H}$. The same assertions hold if we replace $\operatorname{Gr}(L)$ by an open subset $N_{\Upsilon}$ as defined in Lemma \ref{NXi}.
\end{cor}
As the exponential map of the corresponding normalized coordinate from $\mathcal{H}$ takes the open subset of $T_{\mathcal{U}_{\mathbb{Z}}(\Upsilon)}^{1}(\mathbb{C})$ arising either from quotients of $\operatorname{Gr}(L)$ or of $N_{\Upsilon}$ onto a punctured open disc in $\mathbb{C}$, all the bundles from Corollary \ref{D/UZUpsilon} are punctured disc bundles.

\subsection{The Toroidal Boundary Components}

We begin by considering the toroidal boundary components over 1-dimensional cusps. These boundary components are canonical, and their description does not require the gluing with the boundary components lying over the 0-dimensional cusps. Indeed, we have seen in Propositions \ref{dim2} and \ref{struc1dim} that $\mathcal{U}_{\mathbb{Z}}(\Upsilon)$ has rank 1 and hence $\dim T_{\mathcal{U}_{\mathbb{Z}}(\Upsilon)}=1$, and the cone $\Omega(\Upsilon)$ is a ray in $\mathcal{U}(\Upsilon)$. There is therefore only one possible cone decomposition: We must take $\Sigma(\Upsilon)$ to be $\big\{\Omega(\Upsilon)\cup\{0\},\{0\}\big\}$. It follows that if $\nu_{\Upsilon}$ is the element of $\mathcal{U}_{\mathbb{Z}}(\Upsilon)^{*}$ taking the generator $u_{\Upsilon}$ from Proposition \ref{dim2} to 1, then the toric variety $X_{\mathcal{U}_{\mathbb{Z}}(\Upsilon),\Sigma(\Upsilon)}$ is $\operatorname{Spec}\mathbb{C}[z_{\nu_{\Upsilon}}]\cong\mathbb{A}^{1}$, and the set $\mathcal{D}_{\Upsilon}$ from \eqref{fibDXi} is the union of $\mathcal{U}_{\mathbb{Z}}(\Upsilon) \backslash \operatorname{Gr}(L)$ and $\widetilde{\mathcal{D}}(\Upsilon)$. Note that the $\mathbb{A}^{1}$-coordinate $z_{\nu_{\Upsilon}}$ is the coordinate denoted by $q_{2}=\mathbf{e}(\tau_{2}')$ in \cite{[K3]}, so that the equation $z_{\nu_{\Upsilon}}=0$ for that boundary component, given in the proof of Proposition \ref{ordTdiv}, is just $q_{2}=0$ from that reference. As Proposition \ref{dim2} shows that $\mathcal{P}_{\mathbb{Z}}(\Upsilon)$ maps trivially to $G_{\ell}(\Upsilon)$, no quotient of the fiber $X_{\mathcal{U}_{\mathbb{Z}}(\Upsilon),\Sigma(\Upsilon)}$ has to be taken in the construction of $X^{\mathrm{tor}}$. In particular, Theorem \ref{propdisc} becomes trivial for our $\Upsilon$.

Combining this information with Corollary \ref{D/UZUpsilon} yields the form of $\mathcal{D}_{\Upsilon}$ and its quotients. The set $N_{\Upsilon}$ from Lemma \ref{NXi} is defined by the $\Phi$-image being large enough in the 1-dimensional ray $\Omega(\Upsilon)$. We recall that \eqref{quotsZ} denotes the quotient $\mathcal{P}_{\mathbb{Z}}(\Upsilon)/\mathcal{U}_{\mathbb{Z}}(\Upsilon)$ by $\overline{\mathcal{P}}_{\mathbb{Z}}(\Upsilon)$, that Lemma \ref{NXi} defines the open subset $N_{\Upsilon}$ of $\mathcal{D}$, and the notation for the relative interior of the closure appearing in \eqref{intclos} and after it, and obtain the following description of the sets $\mathcal{D}_{\Upsilon}$ and $\overline{\mathcal{U}_{\mathbb{Z}}(\Upsilon) \backslash N_{\Upsilon}}^{o}$ and their quotients in this case.
\begin{lem} \label{DXidim1}
The set $\mathcal{D}_{\Upsilon}$ and the open subset $\overline{\mathcal{U}_{\mathbb{Z}}(\Upsilon) \backslash N_{\Upsilon}}^{o}$ are full, non-punctured open disc bundles over $\widetilde{\mathcal{D}}(\Upsilon)$ from Proposition \ref{struc1dim}, where the zero section of the disc defines an analytic divisor $O(\Upsilon)$. The quotient $\mathcal{V}_{\mathbb{Z}}(\Upsilon)\backslash\mathcal{D}_{\Upsilon}$ and its open subset $\overline{\mathcal{W}_{\mathbb{Z}}(\Upsilon) \backslash N_{\Upsilon}}^{o}$ are similar disc bundles over the family $\mathcal{E} \otimes D$ over $\mathcal{H}$ from Corollary \ref{D/UZUpsilon}. If $\Gamma_{L}$ is neat then $\overline{\mathcal{P}}_{\mathbb{Z}}(\Upsilon)\backslash\mathcal{D}_{\Upsilon}$, as well as the open subset $\overline{\mathcal{P}_{\mathbb{Z}}(\Upsilon) \backslash N_{\Upsilon}}^{o}$, are open disc bundles over the open Kuga--Sato variety $W_{J}$.
\end{lem}
We deduce the following consequence (still under Hypothesis \ref{assonJ} and for neat $\Gamma_{L}$).
\begin{prop} \label{torbcdim1}
The inverse image in $\overline{\mathcal{P}}_{\mathbb{Z}}(\Upsilon)\big\backslash\mathcal{D}_{\Upsilon}$ of the image of $\Upsilon$ in $\mathcal{P}_{\mathbb{Z}}(\Upsilon)\backslash\mathcal{D}^{\mathrm{BB}}$ has a smooth neighborhood that is biholomorphic to the product of $W_{J}$ with a disc, in which the boundary component itself is defined by the vanishing of the disc coordinate. Therefore to $\Upsilon$ one associates a canonical divisor in $\overline{\mathcal{P}}_{\mathbb{Z}}(\Upsilon)\big\backslash\mathcal{D}_{\Upsilon}$, which is isomorphic to $W_{J}$.
\end{prop}
Proposition \ref{torbcdim1} follows immediately from the fact that the bundles in Lemma \ref{DXidim1} are trivial. We remark again that Lemma \ref{DXidim1} and Proposition \ref{torbcdim1} do not depend on choices of fans, since $\Sigma(\Upsilon)$ is canonical.
\begin{rmk} \label{genlat}
Without Hypothesis \ref{assonJ}, the form of the boundary divisor of $\overline{\mathcal{P}}_{\mathbb{Z}}(\Upsilon)\big\backslash\mathcal{D}_{\Upsilon}$ can be significantly more complicated---the interested reader may consult Section 4 of \cite{[Ze3]} for the details. As mentioned above, Hypothesis \ref{assonJ} is only used for the simpler description of the parts of Corollary \ref{D/UZUpsilon}, Lemma \ref{DXidim1}, and Proposition \ref{torbcdim1} involving the Kuga--Sato variety, and all the results that we present from now on are equally valid without it.
\end{rmk}

\medskip

Consider now a 0-dimensional cusp $\Xi$, associated with the rank 1 isotropic sublattice $I \subseteq L$, with $K$ as in \eqref{IKdef}. The toroidal boundary components lying over $\Xi$ depend on a choice of a fan $\Sigma(\Xi)$, which has to be an admissible cone decomposition in the sense of Definition \ref{admissible}. These components must lie, in correspondence with Corollary \ref{D/UZXi}, in the direction of $\Omega(\Xi)$, which equals $C$ by Proposition \ref{struc0dim}, but we shall also need to consider the images $\Omega(\Upsilon)$ for boundary components $\Upsilon$ of $\mathcal{D}^{\mathrm{BB}}$ with $\Xi\leq\Upsilon$. These are $\Upsilon=\mathcal{D}=\operatorname{Gr}(L)$, with $\Omega(\mathcal{D})=\{0\}$, and those components $\Upsilon$ associated with primitive isotropic rank 2 sublattices $J \subseteq L$ such that $I \subseteq J$. Given such $J$, note that $J \subseteq I^{\perp}_{L}$, and the quotient $J/I$ is a rank 1 primitive isotropic sublattice of $K$, generating the isotropic line $J_{\mathbb{R}}/I_{\mathbb{R}} \subseteq K_{\mathbb{R}}$. Assuming that the orientations and the bases are chosen appropriately, we get a generator $z$ for $I$, a generator $\omega_{J,I}$ for $J/I$ that lies on the boundary of the cone $C$, and an element $w$ of $J$ such that $(z,w)$ is an oriented basis for $J$ over $\mathbb{Z}$ and $\omega_{J,I}=-w+I$. This extra sign comes from the fact that $w$ pairs positively with elements of $C$ by \eqref{secwithz}, while the vector $\omega_{J,I}$, which lies in the closure of the cone $C$ of negative norm vectors, pairs negatively with elements of $C$.

Lemmas \ref{stdim1} and \ref{stdim2} and Propositions \ref{dim1} and \ref{dim2} allow us to neatly express the intersections $\mathcal{P}(\Xi)\cap\mathcal{P}(\Upsilon)=\mathcal{P}_{I,\mathbb{R}}\cap\mathcal{P}_{J,\mathbb{R}}$ and $\mathcal{P}_{\mathbb{Z}}(\Xi)\cap\mathcal{P}_{\mathbb{Z}}(\Upsilon)$. In particular, the image of the latter inside the Levi quotient from Lemma \ref{stdim1} is the stabilizer $\operatorname{St}_{\Gamma_{K}}(J/I)$, for which an analogue of Proposition \ref{dim1} yields the split short exact sequence
\begin{equation} \label{seqofint}
0\to(D,+)\to\operatorname{St}_{\Gamma_{K}}(J/I)\to\Gamma_{D}\to1.
\end{equation}
with the quotient $\Gamma_{D}$ being trivial in case $\Gamma_{L}$ is neat. Moreover, the inclusions from Lemma \ref{relsubgp} express themselves in the following relations, in which we note that the subgroup $\{0\} \times D_{\mathbb{R}}\times\mathbb{R}$ of $\operatorname{H}(D_{\mathbb{R}})$, and its intersection $\{0\} \times D\times\mathbb{Z}$ with $\operatorname{H}(D)$, are Abelian.
\begin{prop} \label{compint}
The image of $\mathcal{U}(\Xi)$ in the short exact sequence from Lemma \ref{stdim2} consists of the part $\{0\} \times D_{\mathbb{R}}\times\mathbb{R}$ of $\operatorname{H}(D_{\mathbb{R}})$, and of the upper triangular unipotent part of $\operatorname{SL}_{2}(\mathbb{R})$. It is a direct sum of these two groups. After writing $\mathcal{U}(\Xi)=\mathcal{W}(\Xi)\cong(K_{\mathbb{R}},+)$ as in Proposition \ref{paradim0}, we get
\[\mathcal{U}(\Xi)\cap\mathcal{W}(\Upsilon)=\{0\} \times D_{\mathbb{R}}\times\mathbb{R}\cong(J_{\mathbb{R}}/I_{\mathbb{R}})^{\perp}\] as well as \[\mathcal{U}(\Upsilon)=\mathbb{R}u_{\Upsilon} \cong J_{\mathbb{R}}/I_{\mathbb{R}}=\mathbb{R}\omega_{J,I}.\] Hence $\Omega(\Upsilon)$ is the isotropic ray $\mathbb{R}_{+}\omega_{J,I}$ on the boundary of $\Omega(\Xi)=C$. On the other hand, $G_{h}(\Xi)$ is the trivial subgroup of $G_{h}(\Upsilon)$. Finally, the description of $\mathcal{U}_{\mathbb{Z}}(\Xi)=K$ from Proposition \ref{dim1} is obtained by replacing every $\mathbb{R}$ above by $\mathbb{Z}$, and the spaces $K_{\mathbb{R}}$, $D_{\mathbb{R}}$, $J_{\mathbb{R}}/I_{\mathbb{R}}$, and $(J_{\mathbb{R}}/I_{\mathbb{R}})^{\perp}$ by the lattices $K$, $D$, $J/I$, and  $(J/I)^{\perp}_{K}$ respectively.
\end{prop}
The proof of Proposition \ref{compint} is a simple application of Lemma \ref{actU} and the action of the unipotent part of $G_{h}(\Upsilon)$ in Proposition \ref{struc1dim}. This also shows that the action of $\mathcal{U}(\Xi)$ by translations of the real part in the model of $\operatorname{Gr}(L)$ as $K_{\mathbb{R}}^{1}+iC$ from \eqref{GVRKRiC} becomes, in the bundle structure appearing in Proposition \ref{struc1dim}, an additive action on the real parts of the 1-dimensional fiber of $\operatorname{Gr}(L)$ over $\widetilde{\mathcal{D}}(\Upsilon)$, of the fibers of $\mathcal{V}_{1} \otimes D$, and of the variable $\tau\in\mathcal{H}\cong\Upsilon$.

As Proposition \ref{dim1} shows that $\overline{\mathcal{P}}_{\mathbb{Z}}(\Xi)$ is just $\Gamma_{K}$, Proposition \ref{compint} yields the following explicit description of the admissibility condition of $\Sigma(\Xi)$ in Definition \ref{admissible}, for $\Xi=\Xi_{I}$.
\begin{cor} \label{admcdecom}
The fan $\Sigma(\Xi)$ is an admissible cone decomposition as in Definition \ref{admissible} if and only if the following two conditions are satisfied:
\begin{enumerate}[$(i)$]
\item We have \[\bigcup_{\sigma\in\Sigma(\Xi)}\sigma=C\cup\{0\}\cup\Big(\bigcup_{\Xi\leq\Upsilon}\Omega(\Upsilon)\Big)=C\cup\{0\}\cup\Big(\bigcup_{I \subseteq J \subseteq I^{\perp}_{L}}\mathbb{R}_{+}\omega_{J,I}\Big),\]
where $\Upsilon$ is always a 1-dimensional cusp and $J$ is always a rank 2 primitive isotropic sublattice of $L$.
\item $\Sigma(\Xi)$ is $\Gamma_{K}$-invariant, and the number of $\Gamma_{K}$-orbits in $\Sigma(\Xi)$ is finite.
\end{enumerate}
\end{cor}
For proving Corollary \ref{admcdecom}, it only remains to observe that every primitive isotropic vector in $K$ that lies in the closure of $C$ must be the form $\omega_{J,I}$ for a $2$-dimensional primitive isotropic lattice $J\subseteq L$ with $I \subseteq J \subseteq I^{\perp}_{L}$. Note that when $V$ has Witt rank 1, there are no 1-dimensional boundary components $\Upsilon$ and $K$ is anisotropic, which means that the union $\bigcup_{\sigma\in\Sigma(\Xi)}\sigma$ from Corollary \ref{admcdecom} is just $C\cup\{0\}$. In any case, the action of $\Gamma_{K}$ sends boundary components to boundary components.

\smallskip

Let thus $\Sigma(\Xi)$ be admissible as in Corollary \ref{admcdecom}, and consider the toric variety $X_{K,\Sigma(\Xi)}$. The properly discontinuous action of $\Gamma_{K}$ on $X_{K,\Sigma(\Xi)}$, predicted by Theorem \ref{propdisc}, can be verified directly. First, since the orthogonal complement in $K$ of the generator of a ray in $C$ is positive definite, the stabilizer of every such ray in $\Gamma_{K}$ is finite, and the stabilizer of any cone in $\Sigma(\Xi)$ operates on the finitely many rays that generate it. This proves this finiteness of the stabilizer of any cone $\sigma$ containing a ray that lies inside $C$, and with it the fact that the action of $\Gamma_{K}$ on the associated open subvariety $X_{\sigma}$ of $X_{K,\Sigma(\Xi)}$ is properly discontinuous. Next, consider the stabilizer $\operatorname{St}_{\Gamma_{K}}(J/I)$ of the ray $\Omega(\Upsilon)=\mathbb{R}_{+}\omega_{J,I}$ from Proposition \ref{compint} for $\Upsilon$ with $\Xi\leq\Upsilon$ that is associated with $I \leq J \leq I^{\perp}_{L}$. The action of the infinite part $(D,+)$ of that stabilizer, given in \eqref{seqofint}, does not leave any vector on the boundary of $C$ invariant, other than the multiples $\omega_{J,I}$. Hence any cone of dimension at least 2 in $X_{K,\Sigma(\Xi)}$ still has a finite stabilizer in $\Gamma_{K}$, and the action of the larger group $\mathcal{P}_{\mathbb{Z}}(\Upsilon)$ on the orbit $O\big(\Omega(\Upsilon)\big)$ was already verified to be properly discontinuous, essentially in Lemma \ref{DXidim1}. It particular, if $\Gamma_{L}$ is neat then all of these stabilizers are trivial.

Substituting \eqref{GVRmodK} into \eqref{intclos} shows that for $(K_{\mathbb{R}}^{1}/K)+iC \subseteq X_{K,\Sigma(\Xi)}$ we have
\begin{equation} \label{XKSigma}
\overline{(K_{\mathbb{R}}^{1}/K)+iC}^{o}=\bigcup_{\sigma\in\Sigma(\Xi)}\Big(\overline{(K_{\mathbb{R}}^{1}/K)+iC}^{o} \cap O(\sigma)\Big),
\end{equation}
where for every $\sigma\in\Sigma(\Xi)$ that is neither $\{0\}$ nor $\mathbb{R}_{+}\omega_{J,I}$ for some $J$, the intersection is all of $O(\sigma)$. This is also the form of $\mathcal{D}_{\Xi}$ from \eqref{fibDXi}, by the triviality of $\widetilde{\mathcal{D}}(\Xi)$ from Proposition \ref{struc0dim}. The most natural way of defining $\Omega_{N}(\Xi)\subseteq\Omega(\Xi)=C$ for the set $N_{\Xi}$ from Lemma \ref{NXi} is simply by requiring that $|Y^{2}|$ be large enough. Using \eqref{XKSigma} and Proposition \ref{dim1}, we deduce from \eqref{NmodPZ} that the subset $\overline{\mathcal{P}_{\mathbb{Z}}(\Xi) \backslash N_{\Xi}}^{o}$ of $\Gamma_{K} \backslash X_{K,\Sigma(\Xi)}$ is
\begin{equation} \label{neighofXi}
\big(\mathcal{P}_{\mathbb{Z}}(\Xi) \backslash N_{\Xi}\big)\cup\Big\{\Gamma_{K}\Big\backslash\Big[\bigcup_{\{0\}\neq\sigma\in\Sigma(\Xi)}\Big([(K_{\mathbb{R}}^{1}/K)+iC] \cap O(\sigma)\Big)\Big]\Big\},
\end{equation}
with the union consisting of finitely many orbits.

One can also visualize \eqref{DUpsDXi} in our case: Since the comparison of any boundary component with $\mathcal{D}$ itself is immediate in general, the only case to consider is where $\dim\Xi=0$, $\dim\Upsilon=1$, and $\Xi\leq\Upsilon$. Expressing $\mathcal{D}_{\Xi}$ as in \eqref{XKSigma}, we then have to verify that $\big(\mathcal{U}_{\mathbb{Z}}(\Xi)/\mathcal{U}_{\mathbb{Z}}(\Upsilon)\big)\big\backslash\mathcal{D}_{\Upsilon}$ is the intersection of $\mathcal{D}_{\Xi}$ with the affine toric subvariety $X_{K,\Omega(\Upsilon)}$ of $X_{K,\Sigma(\Xi)}$. But indeed, Lemma \ref{DXidim1} shows that $\big(\mathcal{U}_{\mathbb{Z}}(\Xi)/\mathcal{U}_{\mathbb{Z}}(\Upsilon)\big)\big\backslash\mathcal{D}_{\Upsilon}$ is a disc bundle over $\big(\mathcal{U}_{\mathbb{Z}}(\Xi)/\mathcal{U}_{\mathbb{Z}}(\Upsilon)\big)\big\backslash\widetilde{\mathcal{D}}(\Upsilon)$, which by Proposition \ref{compint} is itself fibered over the punctured disc $\operatorname{St}_{\Gamma_{J}}(\binom{1}{0}\big)\backslash\mathcal{H}$, with the fibers being affine models of $D_{\mathbb{C}}/D=T_{D}(\mathbb{C})$. But extending the disc in the fibers from Lemma \ref{DXidim1} to $\mathbb{C}$ and the punctured disc in the base space to $\mathbb{C}^{\times}$ yields an affine structure of $T_{K}$, with the completion of the disc corresponding to the coordinates $z_{\nu}$ for those $\nu \in K^{*}$ having a positive pairing with $u_{\Upsilon}=\omega_{J,I}$ being allowed to vanish. Hence $\big(\mathcal{U}_{\mathbb{Z}}(\Xi)/\mathcal{U}_{\mathbb{Z}}(\Upsilon)\big)\big\backslash\mathcal{D}_{\Upsilon}$ is indeed the desired open subset of $X_{K,\Omega(\Upsilon)}$.

As the set from \eqref{neighofXi} maps isomorphically onto a neighborhood of $\Xi$ in $\mathcal{P}_{\mathbb{Z}}(\Upsilon)\backslash\mathcal{D}^{\mathrm{BB}}$, we apply \eqref{DUpsDXi} and deduce the following complement of Proposition \ref{torbcdim1}.
\begin{prop} \label{torbcdim0}
Denote by $X_{K,\Sigma(\Xi)}^{\mathrm{fib}}$ the complement in $X_{K,\Sigma(\Xi)}$ of the union \[T_{K}\cup\bigcup_{I \subseteq J \subseteq I^{\perp}_{L}}O(\mathbb{R}_{+}\omega_{J,I}),\] with $J$ as usual. It is a closed $\Gamma_{K}$-invariant subset of $X_{K,\Sigma(\Xi)}$, which is contained in $\mathcal{D}_{\Xi}$. Then the inverse image in $\Gamma_{K} \backslash X_{K,\Sigma(\Xi)}$ of the image of $\Xi$ in $\mathcal{P}_{\mathbb{Z}}(\Xi)\backslash\mathcal{D}^{\mathrm{BB}}$ is just $\Gamma_{K} \backslash X_{K,\Sigma(\Xi)}^{\mathrm{fib}}$. It is a closed subset of $\Gamma_{K} \backslash X_{K,\Sigma(\Xi)}$ with a finite stratification, in which each stratum is the quotient of a toric $T_{K}$-orbit modulo a finite group. If $\Gamma_{L}$ is neat and $\Sigma(\Xi)$ is smooth as in Definition \ref{smoothdef}, then this inverse image has a smooth neighborhood in $\Gamma_{K} \backslash X_{K,\Sigma(\Xi)}$.
\end{prop}

We shall need the divisors that are contained in the inverse image from Proposition \ref{torbcdim0}.
\begin{cor} \label{tordivdim1}
There are finitely many $\Gamma_{K}$-orbits of rays $\rho=\mathbb{R}_{+}\omega$ in $\Sigma(\Xi)$, with primitive $\omega \in K$ that are not of the form $\Omega(\Upsilon)$ for a 1-dimensional boundary component $\Upsilon$ with $\Xi\leq\Upsilon$. The vector $\omega \in K \cap C$ satisfies $\omega^{2}=-2N_{\omega}$ for some positive integer $N_{\omega}$, which is constant on $\Gamma_{K}$-orbits. Therefore the boundary from Proposition \ref{torbcdim0} contains finitely many irreducible divisors, each of which is the image of a divisor of the form $D_{\rho}=O_{\Sigma(\Xi)}(\mathbb{R}_{+}\omega)$, and to such a divisor is attached the associated integer $N_{\omega}$.
\end{cor}
It may happen that for some cusp $\Xi$, the set of divisors from Corollary \ref{tordivdim1} will be empty: See Example \ref{Siegel} below, as well as the case considered in \cite{[P]}. In such a case the inverse image from Proposition \ref{torbcdim0} has codimension 2 in $\Gamma_{K} \backslash X_{K,\Sigma(\Xi)}$.

\smallskip

We need to see when is our collection of fans admissible as well. The corresponding condition from Definition \ref{admissible} is simpler in our case than in general.
\begin{lem} \label{GammaL0dim}
For admissibility of the set of fans $\{\Sigma(\Xi)\}_{\Xi}$, it suffices that for two 0-dimensional cusps that are related by $\Gamma_{L}$, say $\Xi$ and $A\Xi$ for $A\in\Gamma_{L}$, we have the equality $\Sigma(A\Xi)=\big\{A\sigma\big|\;\sigma\in\Sigma(\Xi)\big\}$.
\end{lem}
Indeed, the fans $\Sigma(\Upsilon)=\big\{\Omega(\Upsilon),\{0\}\big\}$ for 1-dimensional boundary components $\Upsilon$ are canonical, which in particular implies their $\Gamma_{L}$-invariance, and they are included in $\Sigma(\Xi)$ wherever $\Xi\leq\Upsilon$ in Corollary \ref{admcdecom}.

Since our arithmetic group depends on the choice of $L$, our open Shimura variety and its Baily--Borel compactification from \eqref{XBBdef} are now
\[X_{L}=\Gamma_{L}\backslash\mathcal{D}=\Gamma_{L}\backslash\operatorname{Gr}(L)\stackrel{open}{\subseteq}X_{L}^{\mathrm{BB}}=\Gamma_{L}\backslash\mathcal{D}^{\mathrm{BB}}.\] The toroidal compactification $X_{L}^{\mathrm{tor}}$ of $X_{L}$, and the map onto $X_{L}^{\mathrm{BB}}=\Gamma_{L}\backslash\mathcal{D}^{\mathrm{BB}}$, are therefore constructed, via Theorem \ref{torcomp}, as follows.
\begin{thm} \label{formtc}
Let $\{\Sigma(\Xi)\}_{\Xi\leq \mathcal{D}}$ be an admissible collection of fans as in Definition \ref{admissible}, or equivalently in Corollary \ref{admcdecom} and Lemma \ref{GammaL0dim}. Then $X_{L}$ and the neighborhoods from Lemma \ref{DXidim1} and \eqref{neighofXi} glue together to a normal compact complex space, which is the set of complex points of a complete algebraic variety $X_{L}^{\mathrm{tor}}=X^{\mathrm{tor}}_{\{\Sigma(\Xi)\}_{\Xi}}$ over $\mathbb{C}$. It is smooth if $\Gamma_{L}$ is neat and each fan $\Sigma(\Xi)$ is smooth in the sense of Definition \ref{smoothdef}
\end{thm}
Note that if $\Gamma_{L}$ is not neat, finite quotient singularities may appear at every part of $X_{L}^{\mathrm{tor}}$.

Let $J$ and $\Upsilon$ be as in Proposition \ref{torbcdim1}, and recall that the boundary divisor $O(\Upsilon)$ of $\mathcal{D}_{\Upsilon}$ has the neighborhood $\overline{\mathcal{U}_{\mathbb{Z}}(\Upsilon) \backslash N_{\Upsilon}}^{o}$ from Lemma \ref{DXidim1}. Dividing by $\overline{\mathcal{P}}_{\mathbb{Z}}(\Upsilon)$, we get the neighborhood $\overline{\mathcal{P}_{\mathbb{Z}}(\Upsilon) \backslash N_{\Upsilon}}^{o}$ inside $\overline{\mathcal{P}}_{\mathbb{Z}}(\Upsilon)\big\backslash\mathcal{D}_{\Upsilon}$, which embeds as an open subset of $X_{L}^{\mathrm{tor}}$. We denote the maps \[\pi_{\Upsilon}:\overline{\mathcal{U}_{\mathbb{Z}}(\Upsilon) \backslash N_{\Upsilon}}^{o}\twoheadrightarrow\overline{\mathcal{P}_{\mathbb{Z}}(\Upsilon) \backslash N_{\Upsilon}}^{o}\qquad\mathrm{and}\qquad\iota_{\Upsilon}:\overline{\mathcal{P}_{\mathbb{Z}}(\Upsilon) \backslash N_{\Upsilon}}^{o}\stackrel{open}{\hookrightarrow}X_{L}^{\mathrm{tor}}.\] In addition, if $I$, $K$, $\Xi$, and $\Sigma(\Xi)$ are as above, then dividing the neighborhood $N_{\Xi}$ from Lemma \ref{NXi} by $\mathcal{U}_{\mathbb{Z}}(\Xi)=K$ from Proposition \ref{dim1} gives an open neighborhood $\overline{\mathcal{U}_{\mathbb{Z}}(\Xi) \backslash N_{\Xi}}^{o}$ of the fiber $X_{K,\Sigma(\Xi)}^{\mathrm{fib}}$ from Proposition \ref{torbcdim0} inside $X_{K,\Sigma(\Xi)}$, and dividing further by $\overline{\mathcal{P}}_{\mathbb{Z}}(\Xi)=\Gamma_{K}$ gives the open subset $\overline{\mathcal{P}_{\mathbb{Z}}(\Xi) \backslash N_{\Xi}}^{o}$ of $\Gamma_{K} \backslash X_{K,\Sigma(\Xi)}$. Also here we set the notation \[\pi_{\Xi}:\overline{\mathcal{U}_{\mathbb{Z}}(\Xi) \backslash N_{\Xi}}^{o}\twoheadrightarrow\overline{\mathcal{P}_{\mathbb{Z}}(\Xi) \backslash N_{\Xi}}^{o}\qquad\mathrm{and}\qquad\iota_{\Xi}:\overline{\mathcal{P}_{\mathbb{Z}}(\Xi) \backslash N_{\Xi}}^{o}\stackrel{open}{\to}X_{L}^{\mathrm{tor}}.\] Note that both $\iota_{\Upsilon}$ and $\iota_{\Xi}$ have open images in $X_{L}^{\mathrm{tor}}$, but while $\iota_{\Upsilon}$ is a homeomorphism onto its image (whence the notation there), $\iota_{\Xi}$ is not because the orbits $O(\mathbb{R}_{+}\omega_{J,I})$ are divided by $\mathcal{W}_{\mathbb{Z}}(\Upsilon)/\mathcal{U}_{\mathbb{Z}}(\Xi) \cong D$. Taking $\omega$ and $\rho$ to be as in Corollary \ref{tordivdim1}, we define the irreducible toroidal divisors
\begin{equation}
\begin{split} \label{BJBIomega}
B_{I,\omega} &=(\iota_{\Xi}\circ\pi_{\Xi})(D_{\rho})=(\iota_{\Xi}\circ\pi_{\Xi})\big(O_{\Sigma(\Xi)}(\mathbb{R}_{+}\omega)\big), \\ B_{J} &=\overline{(\iota_{\Upsilon}\circ\pi_{\Upsilon})\big(O(\Upsilon)\big)}.
\end{split}
\end{equation}
Then the boundary of $X_{L}^{\mathrm{tor}}$ consists of a finite number of such divisors: One divisor $B_{J}$ for every $J$ in a set of representatives for the $\Gamma_{L}$-orbits of isotropic sublattices of rank 2 in $L$, and one divisor $B_{I,\omega}$ for every $I$ in a similar set of representatives having rank 1, and for $\omega$ in a set of representatives for the $\Gamma_{K}$-orbits of primitive elements generating internal rays in the fan corresponding to $I$, where $K=I^{\perp}_{L}/I$. Clearly the image in $X_{L}^{\mathrm{tor}}$ of each divisor does not depend on the representatives $J$ or $I$ and $\omega$, and in the latter case the value $-2N_{\omega}$ of $\omega^{2}$ is also independent of the representative. We conclude by recalling that the divisor $B_{I,\omega}$ from \eqref{BJBIomega} contains also the $\iota_{\Xi}$-images of orbits $O(\sigma)$ for $\sigma\in\Sigma(\Xi)$ with $\rho\leq\sigma$, and as for $B_{J}$, wherever $\Xi\leq\Upsilon$ and $\sigma\in\Sigma(\Xi)$ contains $\Omega(\Upsilon)$ as a face, this divisor $B_{J}$ contains $\iota_{\Xi}\big(O(\sigma)\big)$ as well.

\begin{rmk} \label{congGamma}
The choice of the lattice $L$ determined the group $\Gamma_{L}$ from \eqref{GammaLdef} to be the arithmetic group by which we divide $\mathcal{D}=\operatorname{Gr}(L)$ for obtaining $X_{L}$. In some situations in may be useful to consider covers of $X_{L}$ obtained by dividing by (normal) congruence subgroups of $\Gamma_{L}$, especially when $\Gamma_{L}$ is not neat and one wishes to avoid the resulting finite quotient singularities in $X_{L}$. All of these constructions carry over to such covers, yielding compactifications of them, but as results like Proposition \ref{dim1} or Theorem \ref{PZUpsilon} no longer necessarily hold for intersections with subgroups of $\Gamma_{L}$, the precise analysis of the partial quotients with which we work in our investigation may become more difficult.
\end{rmk}

\begin{ex} \label{Siegel}
Let $L$ be the even lattice of signature $(3,2)$ that is obtained as the direct sum of two unimodular hyperbolic planes and the rank 1 lattice generated by a vector of norm 2. Then the Grassmannian from \eqref{Grassdef} can be identified with the \emph{Siegel upper half-plane} $\mathcal{H}_{2}$ of degree 2. Moreover, the group $\operatorname{SO}^{+}(L_{\mathbb{R}})$ is isomorphic to the projective symplectic group $\operatorname{PSp}_{4}(\mathbb{R})$ (though for using explicit matrices one works with the symplectic group $\operatorname{Sp}_{4}(\mathbb{R})$ itself, which is the \emph{spin group} in the orthogonal setting), and the action on $\mathcal{H}_{2}$ is the natural one by ``generalized fractional linear transformations''. Then $\Gamma_{L}$ is the projective image of $\operatorname{Sp}_{4}(\mathbb{Z})$, and it acts transitively on the rank 1 primitive isotropic lattices in $L$, as well as on the rank 2 primitive isotropic lattices in $L$. Hence there is a single 0-dimensional cusp, and a single 1-dimensional cusp.

If we represent the 1-dimensional cusp by the rank 2 lattice $J$, then the fact that $\Delta_{L}$ is anisotropic implies that $J_{L^{*}}=J$ and thus $\Gamma_{J}$ is all of $\operatorname{SL}_{2}(\mathbb{Z})$, and $\Gamma_{D}$ is trivial since $D$ has rank 1. Therefore the 1-dimensional cusp is isomorphic to the classical $j$-plane, and Proposition \ref{torbcdim1} implies that its inverse image is the Kuga--Sato variety $W_{1}$ of level 1, with 1-dimensional fibers. Note that $\operatorname{Sp}_{4}(\mathbb{Z})$ is not neat, which is visible in $\Gamma_{J}=\operatorname{SL}_{2}(\mathbb{Z})$ having torsion elements, and thus the Kuga--Sato variety $W_{1}$ is not a universal family.

Representing the 0-dimensional cusp by the rank 1 lattice $I$, the resulting lattice $K$ is isomorphic to traceless $2\times2$ matrices over $\mathbb{Z}$ with the quadratic form being $-\det$, and the cone $C$ consists of the real traceless matrices $M$ with $\det M>0$ and $\operatorname{Tr}(SM)<0$, where $S$ is the inversion matrix $\binom{0\ \ -1}{1\ \ \ 0\ }$ from \eqref{SL2Zgen} below. We have $\Gamma_{K}=\operatorname{PSL}_{2}(\mathbb{Z})$, so that the corresponding spin group is $\operatorname{SL}_{2}(\mathbb{Z})$. The first rational cone decomposition $\Sigma$ that we take is the one consisting of the $\operatorname{SL}_{2}(\mathbb{Z})$-translates of the cone \[\sigma=\mathbb{R}_{+}\begin{pmatrix} 0 & -1 \\ 0 & 0\end{pmatrix}+\mathbb{R}_{+}\begin{pmatrix} 0 & 0 \\ 1 & 0\end{pmatrix}+\mathbb{R}_{+}\begin{pmatrix} 1 & -1 \\ 1 & -1\end{pmatrix}\] and their faces, as in page 10 of \cite{[Nam]} or Definition 3.117 of \cite{[HKW]}. Note that $\Gamma_{K}$ acts transitively on the cones of fixed dimension, and that only the cones of dimensions 2 and 3 intersect $C$ (i.e., are not contained in the boundary), so that over the cusp of the basis of $W_{1}$ consists of a 1-dimensional orbit and a 0-dimensional orbit, each modulo its stabilizer. The stabilizer in $\Gamma_{K}$ of each 3-dimensional cone has order 3, and it identifies the three 2-dimensional faces. Each 2-dimensional cone, like the face $\tau$ of $\sigma$ spanned by the first two generators, has a stabilizer of order 2, generated by $S$. Since this stabilizer acts as inversion on the multiplicative coordinate of $O(\tau)\cong\mathbb{C}^{\times}$, the quotient modulo its stabilizer is isomorphic to $\mathbb{C}$. The closure $O_{\Sigma}(\tau)$ from \eqref{Osigmadef} is obtained by adding $O(\sigma)$ as well as the point $O(\hat{\sigma})$ associated with the other cone $\hat{\sigma}$ with $\tau\leq\hat{\sigma}$, with the third generator $\binom{-1\ \ -1}{\ \ 1\ \ \ \ 1}$. The stabilizer identifies these two points, and thus this construction indeed produces the compactification of $X_{L}$ by adding the natural compactification of $W_{1}$ by a projective line over the cusp.

For having a cone decomposition with internal rays, we consider the subdivision of $\sigma$ that is obtained by adding the ray $\rho=\mathbb{R}_{+}\omega$ for $\omega=\binom{1\ \ -2}{2\ \ -1}$. This element is invariant under the elements $(TS)^{\pm1}$ stabilizing $\sigma$, and we have $N_{\omega}=3$. It follows that if $\tilde{\sigma}$ is the cone generated by $\tau$ and $\rho$ then the resulting rational cone decomposition $\tilde{\Sigma}$ consists of the $\operatorname{SL}_{2}(\mathbb{Z})$-translates of $\tilde{\sigma}$ and their faces. Note that now there are 1-dimensional cones that are isotropic and anisotropic, and the 2-dimensional cones split into those that are spanned by two isotropic rays and those having a non-isotropic ray as a face. As for the orbits, $\Gamma_{K}$ again acts transitively on the 3-dimensional cones, here with no stabilizers. It acts transitively on the rays of each type, with the anisotropic ones having stabilizers of order 3. The action on each type of 2-dimensional cones is also transitive, with those like $\tau$ keeping their stabilizers and the others having trivial stabilizers. In the resulting compactification, $O_{\tilde{\Sigma}}(\tau)$ yields a similar projective line, but we also have the image of $O_{\tilde{\Sigma}}(\rho)$, which intersects it only at the point $O(\tilde{\sigma})$. As the latter closed orbit consists of a 2-dimensional torus bounded by a chain of three projective lines, the quotient is the union of a double cover of a plane from $O(\rho)$, a copy of $\mathbb{C}^{\times}$ from $O(\tilde{\tau})$, and the intersection point $O(\tilde{\sigma})$ completing the latter to a projective line with 0 and $\infty$ glued together. The explicit form of the quotient of $O(\rho)$ is isomorphic to \[\operatorname{Spec}\mathbb{C}[X,Y,Z]/\big(Z^{2}-(XY-3)Z+(X^{3}+Y^{3}-6XY+9)\big).\] Note that in both cases the compactification does not look locally like a toric variety, because of the non-trivial stabilizers.
\end{ex}

\section{Green Functions for Special Divisors \label{Green}}

In this section we introduce the special divisors on $X_{L}$ and their Green functions, and investigate their extensions to the toroidal compactification $X_{L}^{\mathrm{tor}}$.

\subsection{Special Divisors and Theta Lifts of Harmonic Maass Forms}

Given an element $\lambda \in L^{*}$ with $\lambda^{2}>0$, the set
\begin{equation} \label{nuperp}
\lambda^{\perp}=\big\{v\in\operatorname{Gr}(L)\big|\;v\perp\lambda\big\},
\end{equation}
whose pre-image in $P^{+}$ is $\{Z_{L} \in P^{+}|\;(Z_{L},\lambda)=0\big\}$, is a sub-Grassmannian of complex codimension 1 in $\operatorname{Gr}(L)$, hence an analytic divisor. We shall be interested in the following combinations of the divisors from \eqref{nuperp}, which are called \emph{special divisors}.
\begin{defn} \label{divdef}
For an element $\mu\in\Delta_{L}$ and an number $0<m\in\frac{\mu^{2}}{2}+\mathbb{Z}$ we define the \emph{special divisor of discriminant $(m,\mu)$ on $X_{L}$} to be \[Z(m,\mu)=\Gamma_{L}\bigg\backslash\sum_{\substack{\lambda \in \mu+L \\ \lambda^{2}=2m}}\lambda^{\perp}.\] We shall use $Z(m,\mu)$ also for the closure of this divisor in $X_{L}^{\mathrm{BB}}$, in $X_{L}^{\mathrm{tor}}=X^{\mathrm{tor}}_{\{\Sigma(\Xi)\}_{\Xi}}$, or in any analytic manifold containing a quotient of $\operatorname{Gr}(L)$ as a dense open subset.
\end{defn}
Since the number of $\Gamma_{L}$-orbits associated with any $m$ and $\mu$ in Definition \ref{divdef} is finite, $Z(m,\mu)$ is an algebraic divisor on $X_{L}$. The multiplicities of the irreducible components of $Z(m,\mu)$ are 2 if $2\mu=0\in\Delta_{L}$, and 1 otherwise.

\smallskip

Recall now the factor of automorphy from \eqref{GL2actH}, and let \[\operatorname{Mp}_{2}(\mathbb{Z})=\big\{(M,\varphi)\big|\;M \in \operatorname{SL}_{2}(\mathbb{Z}),\ \varphi:\mathcal{H}\to\mathbb{C},\ \varphi^{2}(\tau)=j(M,\tau)\big\}\] be the metaplectic extension of $\operatorname{SL}_{2}(\mathbb{Z})$, with the product $(M,\varphi)(N,\psi)=(MN,(\varphi \circ N)\cdot\psi)$. It is generated by the elements
\begin{equation} \label{SL2Zgen}
T=\bigg(\begin{pmatrix} 1 & 1 \\ 0 & 1\end{pmatrix},1\bigg)\qquad\mathrm{and}\qquad S=\bigg(\begin{pmatrix} 0 & -1 \\ 1 & 0\end{pmatrix},\sqrt{\tau}\bigg).
\end{equation}
Consider our even lattice $L$, of signature $(n,2)$, and let $\{\mathfrak{e}_{\mu}\}_{\mu\in\Delta_{L}}$ denote the canonical basis of the group ring $\mathbb{C}[\Delta_{L}]$. We denote by $\rho_{L}$ the \emph{Weil representation} of $\operatorname{Mp}_{2}(\mathbb{Z})$ on $\mathbb{C}[\Delta_{L}]$, normalized as in \cite{[Bo1]}, \cite{[Br1]}, or \cite{[Ze1]}. The dual representation $\rho_{L}^{*}$ of $\rho_{L}$ is canonically isomorphic to the complex conjugate representation $\overline{\rho}_{L}$.

Let $k\in\frac{1}{2}\mathbb{Z}$ be a weight, and denote the weight $k$ \emph{hyperbolic Laplacian}, in the variable $\tau=x+iy\in\mathcal{H}$, by \[\Delta_{k}=-y^{2}(\partial_{x}^{2}+\partial_{y}^{2})+ky(\partial_{y}-i\partial_{x})=-\xi_{2-k}\xi_{k},\] where $\xi_{k}$ is the operator from \cite{[BFu]} which is defined by $\xi_{k}f=2iy^{k}\overline{\partial_{\overline{\tau}}f}$. Recall that a \emph{harmonic Maass form} of weight $k$ and representation $\overline{\rho}_{L}$ is a real-analytic function $f:\mathcal{H}\to\mathbb{C}[\Delta_{L}]$ satisfying the functional equation
\begin{equation} \label{moddef}
f(A\tau)=\varphi(\tau)^{2k}\overline{\rho}_{L}(A,\varphi)f(\tau)
\end{equation}
for every $(A,\varphi)\in\operatorname{Mp}_{2}(\mathbb{Z})$ and $\tau\in\mathcal{H}$, the differential equation $\Delta_{k}f=0$, and such that the growth of $f$ towards the cusp $\infty$ is at most linear exponential. The $\xi_{k}$-image of such a modular form is a \emph{weakly holomorphic} modular form of weight $2-k$ and representation $\rho_{L}$ (i.e., it is holomorphic on $\mathcal{H}$, with at most a pole at $\infty$, and satisfies the functional equation from \eqref{moddef} with $2k$ replaced by $4-2k$ and $\overline{\rho}_{L}$ replaced by $\rho_{L}$), see \cite{[BFu]}. We denote by $M_{k}^{!}(\overline{\rho}_{L})$ the space of weakly holomorphic modular forms of weight $k$ and representation $\overline{\rho}_{L}$, by $M_{k}(\overline{\rho}_{L})$ the subspace consisting of holomorphic modular forms, by $S_{2-k}(\rho_{L})$ the space of cusp forms of weight $2-k$ and representation $\rho_{L}$, and by $H_{k}(\overline{\rho}_{L})$ the space of harmonic Maass form of weight $k$ and representation $\overline{\rho}_{L}$ whose $\xi_{k}$-image is cuspidal. Corollary 3.8 of \cite{[BFu]} provides us with the exactness of the sequence \[0 \to M_{k}^{!}(\overline{\rho}_{L}) \to H_{k}(\overline{\rho}_{L})\stackrel{\xi_{k}}{\to}S_{2-k}(\rho_{L})\to0.\]

The modularity condition with respect to $T$ implies that $F \in H_{k}(\overline{\rho}_{L})$ admits the Fourier expansion
\begin{equation} \label{Fourier}
F(\tau)=\sum_{\delta\in\Delta_{L}}\sum_{l\in-\frac{\delta^{2}}{2}+\mathbb{Z}}c^{+}(\delta,l)q^{l}\mathfrak{e}_{\delta}+
\sum_{\delta\in\Delta_{L}}\sum_{0>l\in-\frac{\delta^{2}}{2}+\mathbb{Z}}c^{-}(\delta,l)\Gamma(1-k,4\pi|n|y)q^{l}\mathfrak{e}_{\delta},
\end{equation}
with $q=\mathbf{e}(\tau)$ as usual, such that $c^{+}(\delta,l)=0$ for $l\ll0$, and where $\Gamma(s,x)=\int_{x}^{\infty}t^{s}e^{-t}\frac{dt}{t}$ is the \emph{incomplete Gamma function}. We denote by $F^{+}(\tau)$ the first sum in \eqref{Fourier}, called the \emph{holomorphic part} of $F(\tau)$, and the second sum, which is the \emph{non-holomorphic part} of $F(\tau)$, by $F^{-}(\tau)$. Moreover, restricting the sum defining the holomorphic part $F^{+}(\tau)$ to include only the indices $l<0$ yields the (finite) \emph{principal part} of $F$. Proposition~3.11 of \cite{[BFu]} then yields the following result.
\begin{prop} \label{Poincare}
For every $\mu\in\Delta_{L}$ and $0<m\in\frac{\mu^{2}}{2}+\mathbb{Z}$ there is a harmonic Maass form $F_{m,\mu} \in H_{1-\frac{n}{2}}(\overline{\rho}_{L})$, which is unique up to an element of $M_{1-\frac{n}{2}}(\overline{\rho}_{L})$, whose principal part is $q^{-m}(\mathfrak{e}_{\mu}+\mathfrak{e}_{-\mu})$.
\end{prop}
We remark that when $n\geq3$, where the harmonic Maass form $F_{m,\mu}$ from Proposition \ref{Poincare} is uniquely defined, it is alternatively constructed in \cite{[Br1]} as a Poincar\'{e} series. When $n=2$ it is unique up to a $\overline{\rho}_{L}$-invariant vector, and when $n=1$ up to a larger degree of freedom, and the Hecke trick extends Poincar\'{e} series also to these cases, with some particular choice of the element from $M_{1-\frac{n}{2}}(\overline{\rho}_{L})$. In correspondence with the notation from \cite{[Br1]} and \cite{[BFu]}, the specialization of \eqref{Fourier} when $F=F_{m,\mu}$ expands $F_{m,\mu}(\tau)$ as
\begin{equation} \label{PoinFour}
q^{-m}(\mathfrak{e}_{\mu}+\mathfrak{e}_{-\mu})+\sum_{\delta\in\Delta_{L}}\sum_{0 \leq l\in-\frac{\delta^{2}}{2}+\mathbb{Z}}b_{m,\mu}^{+}(\delta,l)q^{l}\mathfrak{e}_{\delta}+
\sum_{\delta\in\Delta_{L}}\sum_{0>l\in-\frac{\delta^{2}}{2}+\mathbb{Z}}b_{m,\mu}^{-}(\delta,l)\Gamma(1-k,4\pi|l|y)q^{l}\mathfrak{e}_{\delta}.
\end{equation}
The coefficients $b_{m,\mu}^{\pm}(\delta,l)$ are evaluated in Proposition 1.10 of \cite{[Br1]}, and they are real numbers. This is true, by the Hecke trick, also when $n=2$ and $l\neq0$, and we assume that $b_{m,\mu}^{+}(\delta,0)\in\mathbb{R}$ as well in this case. While some of the results below hold equally well for $n=1$, the fact that $F_{m,\mu}$ can be far from unique in this case can create some complications. On the other hand, the boundary components are just the classical 0-dimensional cusps in this case, and the theta lifts are investigated in detail in \cite{[BO]}, so that we do not lose much by restricting attention to $n\geq2$.

Given $\tau=x+iy\in\mathcal{H}$ and $v\in\operatorname{Gr}(L)$, we recall the \emph{Siegel theta function}
\begin{equation} \label{Thetadef}
\Theta_{L}(\tau,v)=y\sum_{\lambda \in L^{*}}\mathbf{e}\Big(\tau\tfrac{\lambda_{v^{\perp}}^{2}}{2}+\overline{\tau}\tfrac{\lambda_{v}^{2}}{2}\Big)\mathfrak{e}_{\lambda+L},
\end{equation}
where $\lambda_{v}$ and $\lambda_{v^{\perp}}$ are the projections of $\lambda$ onto the corresponding subspaces of $L_{\mathbb{R}}$. Theorem 4.1 of \cite{[Bo1]} implies that the function $\tau\mapsto\Theta_{L}(\tau,v)$ is modular of weight $\frac{n}{2}-1$ and representation $\rho_{L}$ for any fixed $v\in\operatorname{Gr}(L)$. We recall the pairing sending two vectors $a=\sum_{\delta\in\Delta_{L}}a_{\delta}\mathfrak{e}_{\delta}$ and $b=\sum_{\delta\in\Delta_{L}}b_{\delta}\mathfrak{e}_{\delta}$ to $\langle a,b \rangle_{L}=\sum_{\delta\in\Delta_{L}}a_{\delta}b_{\delta}$, and consider the regularized theta lift
\begin{equation} \label{liftdef}
\Phi_{m,\mu}^{L}(v)=\int_{\operatorname{Mp}_{2}(\mathbb{Z})\backslash\mathcal{H}}^{\mathrm{reg}}\langle\Theta_{L}(\tau,v),F_{m,\mu}(\tau)\rangle_{L}d\mu(\tau)
\end{equation}
of the function $F_{m,\mu}$ from Proposition \ref{Poincare}, where $d\mu(\tau)=\frac{dx\,dy}{y^{2}}$ is the invariant measure on $\mathcal{H}$, with the regularization defined either in \cite{[Bo1]} or in \cite{[Br1]}, and where we can replace $v$ by the complex variable from Proposition \ref{GVcomp}. Proposition 2.11 of the latter reference implies that the difference between the two regularizations is only in the value of the constant $C_{m,\mu}$ below, and is therefore irrelevant to our purposes. Recall the decomposition of the exterior derivative $d$ as $\partial+\overline{\partial}$, as well as the complementary operator $d^{c}=\frac{\partial-\overline{\partial}}{4\pi i}$, so that $dd^{c}=\frac{i}{2\pi}\partial\overline{\partial}$. In addition, if $Z$ is a divisor on $X_{L}$ then we denote by $\delta_{Z}$ the current on $X_{L}$ that takes an element $\eta$ of the space $\mathcal{A}^{2n-2}_{c}(X_{L})$ of compactly supported $(2n-2)$-forms on $X_{L}$ to $\int_{Z}\eta$. Some important properties of the theta lift $\Phi_{m,\mu}^{L}$ are given in the following theorem.
\begin{thm} \label{PhiGreen}
The invariant Laplace operator on $\operatorname{Gr}(L)$, normalized as in \cite{[Br1]}, takes the theta lift $\Phi_{m,\mu}^{L}$ from \eqref{liftdef} to $\frac{n}{4}b_{m,\mu}^{+}(0,0)$. Moreover, the function $\frac{1}{2}\Phi_{m,\mu}^{L}$ is smooth on $X_{L} \setminus Z(m,\mu)$ with a logarithmic singularity along $Z(m,\mu)$. The function $\frac{1}{2}\Phi_{m,\mu}^{L}$ also satisfies the current equation \[dd^{c}[\tfrac{1}{2}\Phi_{m,\mu}^{L}]+\delta_{Z(m,\mu)}=[\eta(m,\mu)]\] as currents on $\mathcal{A}^{2n-2}_{c}(X_{L})$, where $\eta(m,\mu)$ is a smooth harmonic 2-form representing the Chern class of the line bundle that is associated with $Z(m,\mu)$ on $X_{L}$.
\end{thm}
The first assertion in Theorem \ref{PhiGreen} is Theorem 4.7 of \cite{[Br1]} (see \cite{[Ze1]} for analogues of this statement, without contributions from the constant terms, in some different weights). For the other assertions we refer to Theorem~2.11 of {[Br1]} and Theorem~7.3 of \cite{[BFu]}. In particular, $\frac{1}{2}\Phi_{m,\mu}^{L}$ is a logarithmic \emph{Green function} for $Z(m,\mu)$ on $X_{L}$ in the sense of Arakelov geometry, see Section 1.2.4 of \cite{[GS]} and Section 5.6 of \cite{[BKK]}. We remark that \cite{[OT]} presents a different approach for the construction of such \emph{automorphic} Green functions, using spherical functions for symmetric pairs and Poincar\'e series constructions.

\subsection{Green Functions Near the Toroidal Boundary}

We will now investigate the growth behavior of $\Phi_{m,\mu}^{L}$ along the boundary of the toroidal compactification $X_{L}^{\mathrm{tor}}$.

For expressing $\Phi_{m,\mu}^{L}$ near the boundary we shall need some notation. First, for the non-triviality of the boundary of $\mathcal{D}^{\mathrm{BB}}$ there must exist an isotropic sublattice $I$ of rank 1 in $L$, with an oriented generator $z$ and an associated 0-dimensional boundary component $\Xi\in\mathcal{D}^{\mathrm{BB}}$, and we take $K$ to be as in \eqref{IKdef} and $\tilde{I} \subseteq L^{*}$ to be a unimodular complement for $I$, with an oriented generator $\zeta$ pairing to 1 with $z$. Subtracting the image of $\zeta$ transforms the affine tube domain model $K_{\mathbb{R}}^{1}+iC$ from \eqref{GVRKRiC} into the open subset $K_{\mathbb{R}}+iC$ of $K_{\mathbb{C}}$, so that the pairing of an element $Z$ of the latter space with some element $\eta \in K^{*}$ is well-defined. Given $Y \in C$, we write $|Y|$ for $\sqrt{-Y^{2}}$, so that $\frac{Y}{|Y|}$ represents a normalized oriented generator of a negative definite subspace in $K_{\mathbb{R}}$. We also observe that $\zeta \in L^{*}$ has well-defined $\mathbb{Q}/\mathbb{Z}$-pairings with elements of $\Delta_{L}$.

Now, if $\mu \in I^{\perp}_{L^{*}}/I^{\perp}_{L}$ and $\nu \in I^{\perp}_{L^{*}}\cap(\mu+L)$ with $\nu^{2}=2m$ then the closure of the analytic divisor $\nu^{\perp}$ from \eqref{nuperp} in $\mathcal{D}^{\mathrm{BB}}$ contains $\Xi$. Moreover, there are finitely many $\Gamma_{L}$-orbits of such divisors $\nu^{\perp}$, which form a hyperplane arrangement in the terminology of \cite{[L]}. Omitting the intersection of their imaginary parts with the cone $C$ produces a disjoint union of open subsets of $C$ that are called \emph{Weyl chambers}. By a slight abuse of terminology we also call a Weyl chamber the subset $K_{\mathbb{R}}^{1}+iW$ of $K_{\mathbb{R}}^{1}+iC$ for a Weyl chamber $W$ in $C$, and we recall that these Weyl chambers are associated with $m$ and $\mu$, or equivalently with the form $F_{m,\mu}$.

The Grassmannian $\operatorname{Gr}(K)$ that is associated with the Lorentzian lattice $K$, of signature $(n-1,1)$, is defined as in \eqref{Grassdef}, but with the dimension of the negative definite subspaces of $K_{\mathbb{R}}$ being 1. As each such subspace is generated by a unique element $u_{1} \in C$ with $u_{1}^{2}=-1$, we can define the Siegel theta function $\Theta_{K}$ as in \eqref{Thetadef}, but with the Grassmannian variable $v$ replaced by such a vector $u_{1}$ and $y$ replaced by $\sqrt{y}$. The function $\tau\mapsto\Theta_{K}(\tau,u_{1})$ is then modular of the same weight $\frac{n}{2}-1$ and representation $\rho_{K}$ for every such $u_{1}$. Recall the projection $p_{K}$ from \eqref{projDelta}, as well as the operator $\downarrow^{L}_{K}$ from \eqref{arrowops}, and then we define, like in \eqref{liftdef}, the theta lift
\begin{equation} \label{liftK}
\Phi_{m,p_{K}\mu}^{K}(u_{1})=\int_{\operatorname{Mp}_{2}(\mathbb{Z})\backslash\mathcal{H}}^{\mathrm{reg}}\big\langle\Theta_{K}(\tau,u_{1}),\downarrow^{L}_{K}(F_{m,\mu})(\tau)\big\rangle_{K}d\mu(\tau).
\end{equation}
Since when $n\geq3$ the latter function is defined as a Poincar\'{e} series, we deduce that $\downarrow^{L}_{K}(F_{m,\mu})=0$ when $\mu \not\in I^{\perp}_{L^{*}}/I^{\perp}_{L}$, and with it $\Phi_{m,p_{K}\mu}^{K}$ vanishes as well. A similar argument will show that if $n\leq2$ and $\mu \not\in I^{\perp}_{L^{*}}/I^{\perp}_{L}$ then the analogue of $\Phi_{m,p_{K}\mu}^{K}$ can only be a lift of an element of $M_{1-\frac{n}{2}}(\overline{\rho}_{K})$.

The theta lift $\Phi_{m,\mu}^{L}$ from \eqref{liftdef} is described explicitly in the formula preceding Definition~3.11 of \cite{[Br1]}, which we now put in our notation and conventions. We recall that \[\mathcal{V}_{\kappa}(A,B)=\int_{0}^{\infty}\Gamma(1-\kappa,A^{2}t)e^{-B^{2}t-1/t}\frac{dt}{t^{3/2}}\] is the special function from Equation (3.25) of \cite{[Br1]}, which decays as $C_{\varepsilon}e^{-2(1-\varepsilon)\sqrt{A^{2}+B^{2}}}$ as $A^{2}+B^{2}\to\infty$. We define $\Phi_{m,p_{K}\mu}^{K}$ as in \eqref{liftK}, and we will describe it in Proposition \ref{PhimmuK} below.
\begin{prop} \label{PhimmuL}
Let $Z_{L} \in P_{+}$ be such that $(Z_{L},z)=1$ and the image $Y$ of $\Im Z_{L}$ in $C \subseteq K_{\mathbb{R}}$ lies in a Weyl chamber $W$ that is associated with $F_{m,\mu}$ and satisfies $|Y^{2}|>2m$. Then there is a constant $C_{m,\mu}$ such that the equality \[\Phi_{m,\mu}^{L}(Z_{L})=\frac{|Y|}{\sqrt{2}}\Phi_{m,p_{K}\mu}^{K}\bigg(\frac{Y}{|Y|}\bigg)+C_{m,\mu}-b_{m,\mu}^{+}(0,0)\log|Y^{2}|+\] \[-4\sum_{\substack{\nu\in[I^{\perp}_{L^{*}}\cap(\pm\mu+L)]/I \\ \nu^{2}=2m \\ (\nu+I_{\mathbb{R}},W)>0}}\log\big|1-\mathbf{e}\big((\nu,Z_{L})\big)\big|-4\sum_{\substack{\nu \in I^{\perp}_{L^{*}}/I \\ \nu^{2}\leq0 \\ (\nu+I_{\mathbb{R}},W)>0}}b_{m,\mu}^{+}\big(\nu+L,-\tfrac{\nu^{2}}{2}\big)\log\big|1-\mathbf{e}\big((\nu,Z_{L})\big)\big|+\] \[+\frac{2}{\sqrt{\pi}}\sum_{\substack{\nu \in I^{\perp}_{L^{*}}/I \\ \nu^{2}>0}}b_{m,\mu}^{-}\big(\nu+L,-\tfrac{\nu^{2}}{2}\big)\sum_{h=1}^{\infty}\frac{\mathbf{e}\big(h(\nu,\Re Z_{L})\big)}{h}\mathcal{V}_{1+\frac{n}{2}}\big(\pi h|\nu|\cdot|Y|,\pi h(\nu,Y)\big)\] holds, where the first term and the first sum do not appear if $\mu \not\in I^{\perp}_{L^{*}}/I^{\perp}_{L}$.
\end{prop}

\begin{proof}
Let $Z=X+iY \in K_{\mathbb{R}}+iC$ be the coordinate $Z_{L}-\zeta+I_{\mathbb{C}}$ corresponding to $Z_{L}$. Then the formula from \cite{[Br1]} in question involves expressions like $\mathbf{e}\big((\delta,\zeta)+(\lambda,Z)\big)$ for $\lambda \in K^{*}$ and $\delta \in L^{*}_{I}/L$ with $p_{K}\delta=\lambda+K$. The expression for $Z$ and the fact that $\lambda$ is of the form $\nu+I_{L^{*}}$ for some $\nu \in I^{\perp}_{L^{*}}$ imply that $\frac{\nu^{2}}{2}=\frac{\lambda^{2}}{2}$ and $(\lambda,Z)$ is the same as $(\nu,Z_{L}-\zeta)$. In addition, the choice of $\delta$ determines $\nu$ up to $I$ rather than up to $I_{L^{*}}$, and the pairing $(\delta,\zeta)$ is then the image in $\mathbb{Q}/\mathbb{Z}$ of $(\nu,\zeta)$, so that the expression in question is just $\mathbf{e}\big((\nu,Z_{L})\big)$, which is defined for $\nu \in I^{\perp}_{L^{*}}/I$. Substituting all this, and normalizing according to our convention for the signature, gives the result. This proves the proposition.
\end{proof}

\begin{rmk} \label{sumpm}
In the first sum in Proposition \ref{PhimmuL}, the sum is over the relevant elements of $[I^{\perp}_{L^{*}}\cap(\mu+L)]/I$ plus those of $[I^{\perp}_{L^{*}}\cap(-\mu+L)]/I$ (see Equation (3.2) of \cite{[Br1]}), meaning that when $2\mu=0$ in $\Delta_{L}$, and the two sets coincide, we consider every summand twice. The constant $C_{m,\mu}$ there is given, in the regularization from \cite{[Br1]}, by \[b_{m,\mu}^{+}(0,0)(\log2\pi-\gamma)+b_{m,\mu}^{+'}(0,0)-2\sum_{0\neq\ell \in I_{L^{*}}/I}b_{m,\mu}^{+}(\ell,0)\log\big|1-\mathbf{e}\big((\ell,\zeta)\big)\big|,\] where $\gamma$ is the Euler-Mascheroni constant and $b_{m,\mu}^{+'}(0,0)$ is the derivative with respect to the spectral parameter (see Theorem 3.9 of \cite{[Br1]}). Note that by Proposition 2.11 of \cite{[Br1]}, this derivative does not appear in $C_{m,\mu}$ when one works with the regularization from \cite{[Bo1]}. The proof of Proposition \ref{PhimmuL} also shows that the exponents in the formula for $C_{m,\mu}$ do not depend on the choice of $\zeta$ and $\tilde{I}$ as well.
\end{rmk}

As we are interested in the behavior near toroidal boundary divisors, we first need to divide by $\mathcal{U}_{\mathbb{Z}}(\Xi)=K$ and consider $\Phi_{m,\mu}^{L}$ as a function on the subset $(K_{\mathbb{R}}^{1}/K)+iC$ of the principal homogenous space $T_{K}^{1}(\mathbb{C})$ from \eqref{GVRmodK}. This is done as follows.
\begin{lem} \label{PhiLonTK1}
The expression from Proposition \ref{PhimmuL} is invariant under the action of $K$ by translations on the real part, and thus represents the value $\Phi_{m,\mu}^{L}(Z_{L}+K)$ with the variable $Z_{L}+K\in(K_{\mathbb{R}}^{1}/K)+iC \subseteq T_{K}^{1}(\mathbb{C})$ as well.
\end{lem}

\smallskip

For evaluating the first term in Proposition \ref{PhimmuL}, we assume that $L$ has Witt rank 2 and hence $K$ is also isotropic, and we choose an isotropic sublattice of $K$, which is of the form $J/I$ for an isotropic sublattice $J$ of rank 2 in $L$. We extend the basis $z$ of $I$ to an oriented basis $(z,w)$ for $J$, so that the image of $w$ in $K$, which we shall also denote by $w$, is a generator for $J/I$ whose additive inverse lies on the boundary of $C$. A unimodular complement $\tilde{J}$ is spanned by the dual basis, which we can write, if $\tilde{I}$ is chosen appropriately, as consisting of $\zeta\in\tilde{I}$ from above and another element $\upsilon\in\tilde{J} \cap I^{\perp}_{L^{*}}$ pairing to 1 with $w$. We would like to express $\Phi_{m,p_{K}\mu}^{K}$ via Equation (3.11) of \cite{[Br1]}, up to a small correction as follows. Recall the lattice $D$ from \eqref{JDdef} and the map $p_{D}:J^{\perp}_{L^{*}}/J^{\perp}_{L}=L^{*}_{J}/L\to\Delta_{D}$ for $D$ analogous to \eqref{projDelta}. Since $D$ is positive definite, its associated Grassmannian from \eqref{Grassdef} is trivial, and the theta function from \eqref{Thetadef} has just the variable $\tau$, in which it is a holomorphic modular form of the weight $\frac{n}{2}-1$ and representation $\rho_{D}$. Explicitly we get
\begin{equation} \label{ThetaD}
\Theta_{D}(\tau)=\sum_{\beta \in D^{*}}q^{\beta^{2}/2}\mathfrak{e}_{\beta+D} \in M_{\frac{n}{2}-1}(\rho_{D})\quad\mathrm{and}\quad\uparrow^{L}_{D}(\Theta_{D})(\tau)=\sum_{\beta \in J^{\perp}_{L^{*}}/J}q^{\beta^{2}/2}\mathfrak{e}_{\beta+L}\in M_{\frac{n}{2}-1}(\rho_{L}).
\end{equation}
Then the analogue of \eqref{liftdef} and \eqref{liftK} defines the constant $\Phi_{m,p_{D}\mu}^{D}$ as the regularized integral of a non-holomorphic modular function on $\operatorname{Mp}_{2}(\mathbb{Z})\backslash\mathcal{H}$ arising from $\Theta_{D}$ and $\downarrow^{L}_{D}(F_{m,\mu})$, which vanishes by the same argument if $n\geq3$ and $\mu \not\in J^{\perp}_{L^{*}}/J_{L^{*}}$.

Let $\mathbb{B}_{2}$ be the second Bernoulli function, which is 1-periodic and defined by \[\mathbb{B}_{2}(x)=x^{2}-x+\tfrac{1}{6}\qquad\mathrm{for}\quad0 \leq x \leq 1,\] and let $K_{n/2}$ denote the modified Bessel function of the third kind with the parameter $\frac{n}{2}$. Recall that the pairing $(w,Y)$, which is well-defined for $Y \in K_{\mathbb{R}}$ because $w \in I^{\perp}$, is positive by the orientation on $J/I$ and $C$, and we obtain the following expression.
\begin{prop} \label{PhimmuK}
If $\mu \in I^{\perp}_{L^{*}}/I^{\perp}_{L}$ and the inequality $(w,Y)<\frac{|Y|}{2\sqrt{m}}$ holds, then the first term $\frac{|Y|}{\sqrt{2}}\Phi_{m,p_{K}\mu}^{K}\big(\frac{Y}{|Y|}\big)$ in Proposition \ref{PhimmuL} equals \[\frac{|Y^{2}|}{2(w,Y)}\Phi_{m,p_{D}\mu}^{D}+4\pi(w,Y)\sum_{\ell \in J_{L^{*}}/J}b_{m,\mu}^{+}(\ell,0)\mathbb{B}_{2}\big((\ell,\upsilon)\big)+4\pi(w,Y)\!\!\sum_{\substack{\lambda\in[J^{\perp}_{L^{*}}\cap(\mu+L)]/J \\ \lambda^{2}=2m}}\!\mathbb{B}_{2}\bigg(\frac{(\lambda,Y)}{(w,Y)}\bigg)+\] \[+\frac{4(w,Y)}{\pi}\sum_{0\neq\lambda \in J^{\perp}_{L^{*}}/J}b_{m,\mu}^{-}\big(\lambda+L,-\tfrac{\lambda^{2}}{2}\big)\sum_{h=1}^{\infty}\bigg(\frac{\pi h|\lambda|\cdot|Y|}{(w,Y)}\bigg)^{n/2}\mathbf{e}\bigg(h\frac{(\lambda,Y)}{(w,Y)}\bigg)K_{n/2}\bigg(\frac{2\pi h|\lambda|\cdot|Y|}{(w,Y)}\bigg),\] where the third term also appears only when $\mu \in J^{\perp}_{L^{*}}/J_{L^{*}}$.
\end{prop}

\begin{proof}
We know that $\mu \in I^{\perp}_{L^{*}}/I^{\perp}_{L}$, and we take some element $u_{1} \in C$ with $u_{1}^{2}=-1$ and $(w,u_{1})<\frac{1}{2\sqrt{m}}$, which allows us to express $\Phi_{m,p_{K}\mu}^{K}(u_{1})$ via Proposition 3.1 and Equation (3.3) of \cite{[Br1]}. This formula involves arguments like $(\beta,\tilde{\mu})+(\alpha,\upsilon)$, with $\beta \in D^{*}$ and $\alpha\in\Delta_{K}$ that is perpendicular to the image of $J_{L^{*}}$ in $\Delta_{K}$ and whose image in $\Delta_{D}$ is $\beta+D$, and where $\tilde{\mu}$ is the vector $\frac{u_{1}}{(w,u_{1})}-\upsilon+\frac{w}{2(w,u_{1})^{2}}$ in $(J/I)^{\perp}_{\mathbb{R}} \subseteq K_{\mathbb{R}}$ in our signature convention. Each such $\alpha$ is the image of several elements $\delta \in J^{\perp}_{L^{*}}/J^{\perp}_{L}$ with $p_{D}\delta=\beta+D$, for each such $\delta$ we can write $(\alpha,\upsilon)$ as $(\delta,\upsilon)$, and for the Fourier coefficients we recall that $\downarrow^{L}_{K}(F_{m,\mu})(\tau)$ expands as the principal part $q^{-m}(\mathfrak{e}_{p_{K}\mu}+\mathfrak{e}_{-p_{K}\mu})$ plus \[\sum_{\alpha\in\Delta_{K}}\sum_{0 \leq l\in-\frac{\alpha^{2}}{2}+\mathbb{Z}}\sum_{\substack{\delta \in I^{\perp}_{L^{*}}/I^{\perp}_{L} \\  p_{K}\delta=\alpha}}b_{m,\mu}^{+}(\delta,l)q^{l}\mathfrak{e}_{\alpha}+
\sum_{\alpha\in\Delta_{K}}\sum_{0>l\in-\frac{\alpha^{2}}{2}+\mathbb{Z}}\sum_{\substack{\delta \in I^{\perp}_{L^{*}}/I^{\perp}_{L} \\  p_{K}\delta=\alpha}}b_{m,\mu}^{-}(\delta,l)\Gamma(1-k,4\pi|l|y)q^{l}\mathfrak{e}_{\alpha}.\] We can thus work with pairs $(\beta,\delta)$ for such $\beta$ and $\delta$, which are in correspondence with elements $\lambda \in J^{\perp}_{L^{*}}/J$ like in the proof of Proposition \ref{PhimmuL}. The pairing of the $K^{*}$-image such an element, which is its class in $J^{\perp}_{L^{*}}/(J+I_{L^{*}})$, with $w \in J/I$ is 0, and its pairing with $u_{1}$ is well-defined up to integral multiples of $(w,u_{1})$. Since we have $(\beta,\tilde{\mu})=(\lambda,\tilde{\mu})$ and $\frac{\beta^{2}}{2}=\frac{\lambda^{2}}{2}$, and $(\delta,\upsilon)$ is the image of $(\lambda,\upsilon)$ in $\mathbb{Q}/\mathbb{Z}$ once again, it follows that the argument in question is just $(\lambda,\tilde{\mu}+\upsilon)$. By the formula for $\tilde{\mu}$ and the vanishing of $(\lambda,w)$, this is the same as $\frac{(\lambda,u_{1})}{(w,u_{1})}$, which is indeed defined up to integers and is substituted into 1-periodic functions. Recalling that in the first term in Proposition \ref{PhimmuL} we have $k=1-\frac{n}{2}$ and the coefficient $\frac{|Y|}{\sqrt{2}}$, and that the variable $u_{1}=\frac{Y}{|Y|}$ satisfies the inequality $(u_{1},w)<\frac{1}{2\sqrt{m}}$ from Proposition 3.1 of \cite{[Br1]} by our assumption on $Y$, this proves the proposition.
\end{proof}

\subsection{Determining the Growth Towards the Boundary}

As mentioned above, we need to see how automorphic Green functions behave near a boundary component. Lemma \ref{PhiLonTK1} already allows us to work in $T_{K}^{1}$, and we consider the rank 1 isotropic lattice $I$, the associated 0-dimensional boundary component $\Xi$ of $\mathcal{D}^{\mathrm{BB}}$, the lattice $K$ from \eqref{IKdef}, and the cone $C$. We take a cone $\sigma\in\Sigma(\Xi)$ that does not come from $\Sigma(\Upsilon)$ for any 1-dimensional cusp $\Xi\leq\Upsilon$, i.e., such that the orbit $O(\sigma)$, as a subset of $\mathcal{D}_{\Xi}$ from \eqref{fibDXi}, maps to $\Xi\in\mathcal{D}^{\mathrm{BB}}$, and we denote the primitive elements of $K$ that span $\sigma$ as a cone by $\omega_{j}$ with $1 \leq j \leq d$.

Some of the ray generators $\omega_{j}$ of $\sigma$ lie in $C$, while some others can be of the form $\omega_{J,I}$ for some primitive isotropic rank 2 lattice $I \subseteq J \subseteq I^{\perp}_{L}$. We order the indices such that there is some $0 \leq e \leq d$ such that $\omega_{j}$ is $\omega_{J_{j},I}$ for such a lattice $J_{j}$ when $1 \leq j \leq d$, while if $e+1 \leq j \leq d$ then $\omega_{j}$ is a primitive element of $K \cap C$. We take some small $0<\varepsilon<1$, and consider the \emph{open trimmed cone} \[\sigma^{\varepsilon}=\bigg\{\sum_{j=1}^{d}b_{j}\omega_{j}\in\sigma\bigg|\;b_{j}>0\mathrm{\ for\ every\ }1 \leq j \leq d\mathrm{\ and\ }\sum_{j=1}^{d}b_{j}>\varepsilon\bigg\}.\] We can write the open trimmed cone $\sigma^{\varepsilon}$ as
\begin{equation} \label{sigmapar}
\sigma^{\varepsilon}=\bigcup_{k=1}^{e}\sigma_{k}^{\varepsilon}\cup\sigma_{0}^{\varepsilon}=\bigcup_{k=1}^{e}\bigg\{\sum_{j=1}^{d}b_{j}\omega_{j}\in\sigma^{\varepsilon}\bigg|\;\sum_{j \neq k}b_{j}<\varepsilon b_{k}\bigg\}\cup\Big[\textstyle{\sigma^{\varepsilon}\setminus\bigcup_{k=1}^{e}\overline{\sigma_{k}^{\varepsilon/2}}}\Big],
\end{equation}
the union of sets that are open in $\sigma$, where our condition that $\varepsilon<1$ implies that the only non-trivial intersections can be between $\sigma_{k}^{\varepsilon}$ with $k\geq1$ and $\sigma_{0}^{\varepsilon}$. Note that our assumption that $\sigma$ is not the image of any $\Sigma(\Upsilon)$ precisely excludes the case where $d=e=1$. It follows that $\sigma_{0}^{\varepsilon}\neq\emptyset$, and it is clear that when $\sigma \subseteq C\cup\{0\}$, i.e., when $e=0$, the set $\sigma_{0}^{\varepsilon}$ is just the full trimmed cone $\sigma^{\varepsilon}$.

\smallskip

Let $\kappa$ be a subspace of $K_{\mathbb{R}}$ of dimension $\dim\sigma$ that is defined over $\mathbb{Q}$ such that the pairing between $\kappa$ and $\mathbb{R}\sigma$ is non-degenerate. It will be notationally convenient to assume that $\kappa$ contains vectors of negative norm, and therefore its orthogonal complement $\kappa^{\perp}$ is positive definite (or trivial if $\dim\sigma=n$). This yields decompositions \[K_{\mathbb{R}}=\kappa^{\perp}\oplus\mathbb{R}\sigma,\ K_{\mathbb{C}}=\kappa^{\perp}_{\mathbb{C}}\oplus\mathbb{C}\sigma,\mathrm{\ and\ }T_{K}(\mathbb{C})=K_{\mathbb{C}}/K=[\kappa^{\perp}_{\mathbb{C}}/(\kappa^{\perp} \cap K)]\oplus[\mathbb{C}\sigma/(\sigma \cap K)],\] where allowing the real parts to be affine yields decompositions also for $K_{\mathbb{C}}^{1}$ and the affine model $T_{K}^{1}(\mathbb{C})$ from \eqref{GVRmodK}. For an element $Z+K \in T_{K}(\mathbb{C})$ we denote the projections onto the parts from these decompositions as $Z_{\kappa^{\perp}}$ and $Z_{\sigma}$, so that the imaginary part $Y \in K_{\mathbb{R}}$ is $Y_{\kappa^{\perp}}+Y_{\sigma}$ accordingly, and similarly for the affine variable $Z_{L}+K \in T_{K}^{1}(\mathbb{C})$. Note that while at this point is seems natural to take $\kappa=\mathbb{R}\sigma$, which makes these decompositions orthogonal, see the proof of Corollary \ref{prelog} below for the reason why we allow a more general space $\kappa$ here and decompositions that are not orthogonal. The additional boundary components $O(\tau)$ of the open subset $X_{K,\sigma}$ of $X_{K,\Sigma(\Xi)}$, where $\{0\}\neq\tau\in\Sigma(\Xi)$ is a face of $\sigma$, are obtained by considering points $Z+K$ in which the complex coordinates of $Z_{\sigma}$ that are associated with $\tau$ equal the limit $i\infty$, independently of $\Re Z_{\sigma}$.

Now, adding to $\sigma$ boundary points in which some of the coefficients of the generators $\omega_{j}$ may attain the value $+\infty$ yields a compact set, which we denote by $\bar{\sigma}$. Thus extending the range of $Y_{\sigma}$ to $\bar{\sigma}$ makes it defined on all of $X_{K,\sigma}$. For small $\varepsilon>0$ we write $\bar{\sigma}^{\varepsilon}$ for $\sigma^{\varepsilon}\cup(\bar{\sigma}\setminus\sigma)$, and for every $0 \leq k \leq e$ we denote by $\bar{\sigma}_{k}^{\varepsilon}$ the relative interior of the closure of the set $\sigma_{k}^{\varepsilon}$ from \eqref{sigmapar} in $\bar{\sigma}^{\varepsilon}$. Recall that the open subset $\mathcal{D}_{\Xi}$ of $X_{K,\Sigma(\Xi)}$ can be written as in \eqref{XKSigma}, and set \[\mathcal{D}_{\Xi,\sigma}=\mathcal{D}_{\Xi} \cap X_{K,\sigma}\quad\mathrm{and}\quad\mathcal{D}_{\Xi,\sigma}^{\varepsilon}=\big\{Z_{L}+K\in\mathcal{D}_{\Xi,\sigma}\big|\;(\Im Z_{L})_{\sigma}\in\bar{\sigma}^{\varepsilon}\big\},\] where the full orbit $O(\sigma)$ is contained in $\mathcal{D}_{\Xi,\sigma}^{\varepsilon}$. Recall that if $\{0\}\neq\tau\leq\sigma$ is not of the form $\Omega(\Upsilon)\cup\{0\}$ for some $\Upsilon$ then $O(\tau)\subseteq\mathcal{D}_{\Xi,\sigma}$, while if $\tau\neq\sigma$ then some points of $O(\tau)$ will not lie in $\mathcal{D}_{\Xi,\sigma}^{\varepsilon}$. On the other hand, all the points in $O(\tau)$ for such $\tau$ will be covered when we work in $\mathcal{D}_{\Xi,\tau}^{\varepsilon}$. The decomposition from \eqref{sigmapar} allows us to write
\begin{equation} \label{DXisigmak}
\mathcal{D}_{\Xi,\sigma}^{\varepsilon}=\bigcup_{k=0}^{e}\mathcal{D}_{\Xi,\sigma}^{\varepsilon,k}=\bigcup_{k=0}^{e}\big\{Z_{L}+K\in\mathcal{D}_{\Xi,\sigma}^{\varepsilon}\big|\;(\Im Z_{L})_{\sigma}\in\bar{\sigma}_{k}^{\varepsilon}\big\}.
\end{equation}
Recall that if $Y$ is the imaginary part of an element of $\mathcal{D}_{\Xi,\sigma}^{\varepsilon}$ then $Y_{\sigma} \in C$ and hence $Y_{\sigma}^{2}<0$, with $|Y_{\sigma}^{2}|$ being very small only when $Y_{\sigma}$ is close either to 0 or to $\mathbb{R}_{+}\omega_{k}$ for some $1 \leq k \leq e$.

Consider now a compact subset $M$ of $\kappa^{\perp}$, and then for $R>\varepsilon$ we denote by $A_{I}^{R}=A_{\Xi,\sigma,\kappa}^{M,R}$ the set of those $Y=\sum_{j=1}^{d}b_{j}\omega_{j}+Y_{\kappa^{\perp}} \in C$ with $Y_{\kappa^{\perp}} \in M$, $b_{j}\geq0$ for $1 \leq j \leq d$, and $b_{j}>R$ for some $1 \leq j \leq d$. We now prove the following lemma.
\begin{lem} \label{logbound}
If $R$ is large enough then the inequality $|Y^{2}|>2m$ holds for every $Y \in A_{I}^{R}$, and $\log|Y^{2}|$ equals $\log|Y_{\sigma}^{2}|$ plus a smooth function on $A_{I}^{R}$ that decays as $R\to\infty$.
\end{lem}

\begin{proof}
For any $Y \in C$ we can write
\begin{equation} \label{Y2decom}
|Y^{2}|=-Y^{2}=-(Y_{\sigma}+Y_{\kappa^{\perp}})^{2}=|Y_{\sigma}^{2}|-2(Y_{\sigma},Y_{\kappa^{\perp}})-Y_{\kappa^{\perp}}^{2}.
\end{equation}
Since the quadratic form $Y_{\sigma}^{2}$ is bounded away from 0 on $\sigma_{0}^{\varepsilon}$, we deduce that the first term in \eqref{Y2decom} grows quadratically with each $b_{j}$. On the other hand, as $Y_{\kappa^{\perp}}$ is confined to the compact set $M$, the second term is linear in these parameters with bounded coefficients, and the third one is bounded. The first assertion follows directly, and as applying the logarithm to \eqref{Y2decom} produces
\begin{equation} \label{Y2log}
\log|Y^{2}|=\log|Y_{\sigma}^{2}|+\log\bigg[1-\frac{2(Y_{\sigma},Y_{\kappa^{\perp}})}{|Y_{\sigma}^{2}|}-\frac{Y_{\kappa^{\perp}}^{2}}{|Y_{\sigma}^{2}|}\bigg],
\end{equation}
this yields the second one as well. This proves the lemma.
\end{proof}

For the next step, we prove the following result, in which we again identify $\mathcal{D}$ with the set of $Z_{L} \in P_{+}$ with $(Z_{L},z)=1$ as in \eqref{GVRKRiC}.
\begin{lem}
\label{GreenPhiK}
Define a function on the set of $Z_{L}\in\mathcal{D}$ with $Y=\Im Z_{L}$ satisfying $|Y^{2}|>2m$ by putting \[\widetilde{\psi}_{m,\mu}^{L}(Z_{L})=\frac{|Y|}{\sqrt{2}}\Phi_{m,p_{K}\mu}^{K}\bigg(\frac{Y}{|Y|}\bigg)-4\sum_{\substack{\nu\in[I^{\perp}_{L^{*}}\cap(\pm\mu+L)]/I \\ \nu^{2}=2m \\ (\nu+I_{\mathbb{R}},W)>0}}\log\big|1-\mathbf{e}\big((\nu,Z_{L})\big)\big|,\] where $W=W(Y)$ denotes any Weyl chamber associated with $(m,\mu)$ whose closure contains $Y$. Then $\widetilde{\psi}_{m,\mu}^{L}$ is smooth on the complement of $Z(m,\mu)$ with a logarithmic singularity along $2Z(m,\mu)$.
\end{lem}

\begin{proof}
In the terminology of \cite{[Br1]} (adapted to our notation), this is the function $\psi_{m,\mu}^{L}(Z_{L})$ from Equation (3.38) of that reference, plus the expression $\frac{|Y|}{\sqrt{2}}\xi_{m,p_{K}\mu}^{K}\big(\frac{Y}{|Y|}\big)-C_{m,\mu}$, minus the sum from the penultimate term in Proposition \ref{PhimmuL}. The function $\psi_{m,\mu}^{L}$ has the desired properties, by the relation with the holomorphic product $\Psi_{m,\mu}^{L}$ from Definition 3.14 and Equation (3.40) of that reference. Moreover, the condition $(\nu+I_{\mathbb{R}},W)>0$ for $\nu \in I^{\perp}_{L^{*}}$ with $\nu^{2}\leq0$ amounts to $\nu+I_{\mathbb{R}}$ lying in the closure of the cone $-C$, independently of $W$. Hence the sum over these elements appearing in Proposition \ref{PhimmuL} is smooth, and as the same holds for the remaining expression (this follows directly from its definition in Equation (3.12) and Definition 3.3 of \cite{[Br1]}), $\widetilde{\psi}_{m,\mu}^{L}$ also has the required properties. This proves the lemma.
\end{proof}
%
The behavior of our automorphic Green function on $\mathcal{D}_{\Xi,\sigma}^{\varepsilon,0}$ can now be described as follows.
\begin{prop} \label{Phinearomega}
Take $M$ and $R$ as in Lemma \ref{logbound}. Then on the set of $Z_{L}+K \in K\backslash\mathcal{D}$ with $Y=\Im Z_{L}$ lying in $A_{I}^{R}$, the function $\Phi_{m,\mu}^{L}(Z_{L}+K)$ equals
\begin{equation} \label{expPhiinn}
\frac{|Y_{\sigma}|}{\sqrt{2}}\Phi_{m,p_{K}\mu}^{K}\bigg(\frac{Y_{\sigma}}{|Y_{\sigma}|}\bigg)-b_{m,\mu}^{+}(0,0)\log|Y_{\sigma}^{2}|-2\log|\phi|^{2},
\end{equation}
plus a remainder that is smooth on $\mathcal{D}_{\Xi,\sigma}^{\varepsilon,0}$, where $\phi$ is a holomorphic function defining the divisor $Z(m,\mu)$ there.
\end{prop}

\begin{proof}
Lemma \ref{logbound} implies that when $R$ is large enough, we can evaluate $\Phi_{m,\mu}^{L}(Z_{L}+K)$, via Lemma \ref{PhiLonTK1}, as the expression from Proposition \ref{PhimmuL}, where $W$ is an appropriate Weyl chamber for $Y$. This lemma also produces the logarithmic term, up to a smooth function which also contains the constant.

Next, the expressions involving $\mathcal{V}_{1+\frac{n}{2}}$ in Proposition \ref{PhimmuL} are smooth on $K\backslash\operatorname{Gr}(L)$ and decay with each $b_{j}$, and are thus smooth on all of $\mathcal{D}_{\Xi,\sigma}^{\varepsilon,0}$, and even on the whole set $\mathcal{D}_{\Xi,\sigma}$. The sum from the penultimate term in Proposition \ref{PhimmuL} was seen to be smooth on $K\backslash\operatorname{Gr}(L)$ as well, and we note that for $\nu \in I^{\perp}_{L^{*}}$ such that $\nu^{2}\leq0$ and $\nu+I_{\mathbb{R}}$ in the closure of $-C$, the imaginary part of $(\nu,Z_{L})$ grows linearly with every $b_{j}$, $1 \leq j \leq d$, except when $j \leq e$ and $\nu$ is a multiple of $\omega_{j}$. But since in this case other parameters must also grow (since we are $\mathcal{D}_{\Xi,\sigma}^{\varepsilon,0}$), it follows that all the summands in that term are smooth on $\mathcal{D}_{\Xi,\sigma}^{\varepsilon,0}$ as well.

The remaining terms combine to the function $\widetilde{\psi}_{m,\mu}^{L}$ from Lemma \ref{GreenPhiK}. Recalling the formula from Definitions 3.3 and 3.5 of \cite{[Br1]} for $\Phi_{m,p_{K}\mu}^{K}$, we write
\[\widetilde{\psi}_{m,\mu}^{L}(Z_{L})=\frac{|Y|}{\sqrt{2}}\xi_{m,p_{K}\mu}^{K}\bigg(\frac{Y}{|Y|}\bigg)+8\pi\big(\rho_{m,p_{K}\mu}(W),Y\big)-4\sum_{\substack{\nu\in[I^{\perp}_{L^{*}}\cap(\pm\mu+L)]/I \\ \nu^{2}=2m \\ (\nu+I_{\mathbb{R}},W)>0}}\log\big|1-\mathbf{e}\big((\nu,Z_{L})\big)\big|,\] with $W=W(Y)$ as in that lemma and $\rho_{m,p_{K}\mu}(W)$ the corresponding Weyl vector. The function $\xi_{m,p_{K}\mu}^{K}$ is smooth, and the combination of the remaining terms is independent of $W$. We choose some Weyl chamber $W$, and define a function $\Phi_{m,p_{K}\mu}^{K,W}$ on $C$ by setting \[\frac{1}{\sqrt{2}}\Phi_{m,p_{K}\mu}^{K,W}(Y)=\frac{|Y|}{\sqrt{2}}\xi_{m,p_{K}\mu}^{K}\bigg(\frac{Y}{|Y|}\bigg)+8\pi\big(\rho_{m,p_{K}\mu}(W),Y\big).\] This means that up to the homogenization to all of $C$, the function $\Phi_{m,p_{K}\mu}^{K,W}$ is the smooth extension to all of $\operatorname{Gr}(K_{\mathbb{R}})$ of the restriction of $\Phi_{m,p_{K}\mu}^{K}$ to the Weyl chamber $W$. We then express the function from Lemma \ref{GreenPhiK} as
\begin{equation} \label{psidecom}
\widetilde{\psi}_{m,\mu}^{L}(Z_{L})=\frac{1}{\sqrt{2}}\Phi_{m,p_{K}\mu}^{K,W}(Y)-4\sum_{\substack{\nu\in[I^{\perp}_{L^{*}}\cap(\pm\mu+L)]/I \\ \nu^{2}=2m \\ (\nu+I_{\mathbb{R}},W)>0}}\log\big|1-\mathbf{e}\big((\nu,Z_{L})\big)\big|
\end{equation}
for every $Z_{L}\in\mathcal{D}$, regardless of whether $Y$ lies in the closure of $W$ or not. Note that the sum from \eqref{psidecom} easily extends to all of $\mathcal{D}_{\Xi,\sigma}$, and since every element of $\mathcal{D}_{\Xi,\sigma}$ can be perpendicular to only finitely many $\nu\in[I^{\perp}_{L^{*}}\cap(\pm\mu+L)]/I$ with $\nu^{2}=2m$, this function is indeed $-2$ times the logarithm of the square of the absolute value of a holomorphic defining function $\phi$ for the divisor $Z(m,\mu)$ there.

Now, the function $\Phi_{m,p_{K}\mu}^{K,W}$ is smooth and homogenous of degree 1 on $C$. We thus write \[\Phi_{m,p_{K}\mu}^{K,W}(Y)=\Phi_{m,p_{K}\mu}^{K,W}(Y_{\sigma})+|Y_{\sigma}|\bigg[\Phi_{m,p_{K}\mu}^{K,W}\bigg(\frac{Y}{|Y_{\sigma}|}\bigg)-\Phi_{m,p_{K}\mu}^{K,W}\bigg(\frac{Y_{\sigma}}{|Y_{\sigma}|}\bigg)\bigg],\] and express the first term in the brackets by the Taylor expansion of $\Phi_{m,p_{K}\mu}^{K,W}$ at $\frac{Y_{\sigma}}{|Y_{\sigma}|}$. Note that the difference between the arguments is $\frac{Y_{\kappa^{\perp}}}{|Y_{\sigma}|}$, with the numerator lying in the compact set $M$ since we assume that $Y \in A_{I}^{R}$. As the set of norm $-1$ vectors inside the closure of $\sigma_{0}^{\varepsilon}$ in $C$ is compact, and the point $\frac{Y_{\sigma}}{|Y_{\sigma}|}$ around which we expand lies in that set, the coefficients of this expansion are smooth and bounded on all of $\mathcal{D}_{\Xi,\sigma}^{\varepsilon,0}$. Combining the external coefficient with the denominator of the Taylor variable $\frac{Y_{\kappa^{\perp}}}{|Y_{\sigma}|}$, we deduce that this difference is smooth on the desired set. Finally, if we take $W$ to be some Weyl chamber that contains $Y_{\sigma}$ in its closure, the remaining term $\Phi_{m,p_{K}\mu}^{K,W}(Y_{\sigma})$ gives, by the properties of $\Phi_{m,p_{K}\mu}^{K,W}$, the first desired term in \eqref{expPhiinn}. This completes the proof of the proposition.
\end{proof}

\begin{rmk} \label{indepW}
Note that the holomorphic function $\phi$ in Proposition \ref{Phinearomega} depends on the choice of a Weyl chamber containing $Y_{\sigma}$ in its closure. Moreover, it is, in general, not possible to take a Weyl chamber that will work for the entire set there simultaneously. However, changing the Weyl chamber alters $-2\log|\phi|^{2}$ by a linear function of $Y$ (hence a linear function of $Y_{\sigma}$ plus a smooth function on $A_{I}^{R}$), and combining this expression with the first term in Proposition \ref{Phinearomega} yields, as in Lemma \ref{GreenPhiK}, a function with the required logarithmic singularities that no longer depends on choosing a Weyl chamber.
\end{rmk}

Recalling that if all the rays of $\sigma$ are generated by elements of $C$ then $\sigma_{0}^{\varepsilon}=\sigma^{\varepsilon}$ in \eqref{sigmapar}. Hence $\mathcal{D}_{\Xi,\sigma}^{\varepsilon,0}=\mathcal{D}_{\Xi,\sigma}^{\varepsilon}$ in \eqref{DXisigmak}, Proposition \ref{Phinearomega} and Remark \ref{indepW} describe the behavior of $\Phi_{m,\mu}^{L}$ in all of the vectors of $\mathcal{D}_{\Xi,\sigma}^{\varepsilon}$ with large enough imaginary part, and in particular on all the boundary points there. It therefore gives the complete description near the boundary in case $L$ has Witt rank 1. We also remark that when $\dim\sigma=n$ we have $\kappa^{\perp}=0$ and $Y_{\sigma}=Y$, and the expression with $Y$ instead of $Y_{\sigma}$ throughout is also valid in Proposition~\ref{Phinearomega}, and even with a shorter proof. However, $Y_{\sigma}$ grasps the important part of $Y$ near the boundary of $\mathcal{D}_{\Xi,\sigma}^{\varepsilon,0}$, and will be more convenient for our consequences below.

\smallskip

Assuming now that the Witt rank of $L$ is 2, take a primitive isotropic lattice $I \subseteq J \subseteq I^{\perp}_{L}$ of rank 2 for our rank 1 primitive isotropic lattice $I$, with the associated 0-dimensional cusp $\Xi$, 1-dimensional cusp $\Upsilon$ with $\Xi\leq\Upsilon$, and the lattice $K$ from \eqref{IKdef}. Let $\rho$ be the ray $\Omega(\Upsilon)\cup\{0\}=\mathbb{R}_{+}\omega_{J,I}$, considered as an external ray that lies in $\Sigma(\Xi)$, and note that the natural extension of \eqref{sigmapar} to the case where $\sigma=\rho$ is with $\rho_{1}^{\varepsilon}=\rho^{\varepsilon}$ and $\rho_{0}^{\varepsilon}=\emptyset$. Consider $\rho\neq\sigma\in\Sigma(\Xi)$ with $\rho\leq\sigma$, and let $k\geq1$ be the index such that $\omega_{k}=\omega_{J,I}$. Recall that if $\mathcal{U}_{Z}(\Upsilon)$ is the cyclic group associated with the 1-dimensional cusp $\Upsilon$ then in the coordinates from \eqref{GVRKRiC} it corresponds to translations by $J/I=\mathbb{Z}\omega_{J,I}$, and that the quotient of the set $\mathcal{D}_{\Upsilon}$ by $K/\mathbb{Z}\omega_{J,I}$ is an open subset of $\mathcal{D}_{\Xi,\sigma}$. We now turn to investigating $\Phi_{m,\mu}^{L}$ on the set $\mathcal{D}_{\Upsilon}^{\varepsilon}$ defined analogously to the one on the right hand side of \eqref{DXisigmak}, as well as on $\mathcal{D}_{\Xi,\sigma}^{\varepsilon,k}$.

To do so we choose an oriented basis $(z,w)$ for $J$ such that $I=\mathbb{Z}z$ as always, and recall that $\omega_{k}=\omega_{J,I}=-w+I$ in this case. Moreover, note that if $\varepsilon$ is small enough then any Weyl chamber $W$ that has a non-trivial intersection with $\sigma_{k}^{\varepsilon}$ for our $\sigma$ and $k$ must contain $\omega_{k}=\omega_{J,I}$ in its closure, and this is clearly the case also when $\sigma=\rho$ and $\sigma_{1}^{\varepsilon}=\rho^{\varepsilon}$. We take some 1-dimensional rational negative definite vector space $\kappa_{J}$, that is therefore not perpendicular to $\omega_{J,I}$, where for our given $\sigma$, for which we already chose $\kappa$ above in a way that it contains negative norm vectors, we assume that $\kappa_{J}\subseteq\kappa$. It follows that each vector $Y \in C$ can be written as $t\omega_{J,I}+Y_{\kappa_{J}^{\perp}}$ where $Y_{\kappa_{J}^{\perp}}$ lies in a positive definite subspace of $K_{\mathbb{R}}$, and since $\kappa_{J}\subseteq\kappa$ and hence $\kappa^{\perp}\subseteq\kappa_{J}^{\perp}$ we can also write $Y_{\sigma}$ as $t\omega_{J,I}+Y_{\sigma,\kappa_{J}^{\perp}}$.

Take a compact subset $M_{J} \subseteq \kappa_{J}^{\perp}$, and for $T>\varepsilon$ set $A_{J}^{T}=A_{\Upsilon,\kappa_{J}}^{M_{J},T}$ to be those elements $Y=t\omega_{J,I}+Y_{\kappa_{J}^{\perp}} \in C$ for which $Y_{\kappa_{J}^{\perp}} \in M_{J}$ and $t>T$. Similarly, if $\rho=\mathbb{R}_{+}\omega_{J,I}\leq\sigma\in\Sigma(\Xi)$ with $\rho\neq\sigma$, $\kappa$ is as above (with $\kappa_{J}\subseteq\kappa$), $k$ is such that the generator $\omega_{k}$ of $\sigma$ is $\omega_{J,I}$, $M$ is a compact subset of $\kappa^{\perp}$, and $T>\varepsilon$ then we define $A_{I,J}^{T}=A_{\Xi,\Upsilon,\sigma,\kappa,\kappa_{J}}^{M,T}$ to be the set of $Y=Y_{\sigma}+Y_{\kappa^{\perp}} \in C$ with $Y_{\kappa^{\perp}} \in M$ and $Y_{\sigma}=t\omega_{J,I}+Y_{\sigma,\kappa_{J}^{\perp}}$ lying in $\sigma_{k}^{\varepsilon}$ with $t>T$. The analogue of Lemma \ref{logbound} here is as follows.
\begin{lem} \label{loglin}
If $\varepsilon$ is small enough and $T$ is large enough then for $Y$ in $A_{J}^{T}$ or in $A_{I,J}^{T}$ we have $|Y^{2}|>2m$ and $(w,Y)<\frac{|Y|}{2\sqrt{m}}$, and the value of $\log|Y^{2}|$ is $\log t$ plus a smooth function on that set. This function again decays with $T$ on $A_{J}^{T}$, and is bounded by $\log(\varepsilon t)$ on $A_{I,J}^{T}$.
\end{lem}

\begin{proof}
The isotropy of $\omega_{J,I}$ and its relation with $w$ imply that our analog of \eqref{Y2decom} is
\begin{equation} \label{Y2JI}
|Y^{2}|=-Y^{2}=-2t(\omega_{J,I},Y_{\kappa_{J}^{\perp}})-Y_{\kappa_{J}^{\perp}}^{2}=2t(w,Y_{\kappa_{J}^{\perp}})-Y_{\kappa_{J}^{\perp}}^{2},\quad\mathrm{and}\quad(w,Y)=(w,Y_{\kappa_{J}^{\perp}}).
\end{equation}
Considering first the set $A_{J}^{T}$, the two expressions $(w,Y_{\kappa_{J}^{\perp}})$ and $Y_{\kappa_{J}^{\perp}}^{2}$ from \eqref{Y2JI} are bounded, where the former is bounded away from 0 from below because $(w,Y)>0$ for every $Y \in C$. This proves the first two assertions, and as the logarithm of the expression from \eqref{Y2JI} is
\begin{equation} \label{Y2t}
\log|Y^{2}|=\log t+\log(w,Y)+\log\bigg(2-\frac{Y_{\kappa_{J}^{\perp}}^{2}}{t(w,Y)}\bigg),
\end{equation}
the third one follows as well.

For $A_{I,J}^{T}$, we first observe that $Y_{\sigma}=\sum_{j=1}^{k}b_{j}\omega_{j}$, where $b_{k}$ equals $t$ plus a linear combination of the other $b_{j}$'s, and $0 \leq b_{j}<\varepsilon t$ when $j \neq k$. It follows that $(w,Y_{\sigma})=\sum_{j \neq k}b_{j}(w,\omega_{j})$ is a positive linear combination of the $b_{j}$'s with $j \neq k$. Recalling that $(w,Y_{\kappa^{\perp}})$ is bounded for $Y_{\kappa^{\perp}} \in M$, we can define $M_{\sigma}$ to be the set of those combinations $\sum_{j \neq k}b_{j}\omega_{j}$ for which $\sum_{j \neq k}b_{j}(w,\omega_{j})\leq2(w,Y_{\kappa^{\perp}})$ for every $Y_{\kappa^{\perp}} \in M$, and set $M_{J}=M_{\sigma} \times M$. This is a compact subset of $\kappa_{J}^{\perp}$, for elements of which we have already proved our result. We shall therefore consider only $Y \in A_{I,J}^{T} \setminus A_{J}^{T}$ for the resulting $A_{J}^{T}$.

Now, applying \eqref{Y2JI} to $Y_{\sigma}$ and taking the logarithm yields an expression like \eqref{Y2log} for $\log|Y_{\sigma}^{2}|$, where we write $(w,Y_{\sigma})$ as a linear combination of the $b_{j}$'s with $j \neq k$ as above, and note that the numerator, which is now $Y_{\sigma,\kappa_{J}^{\perp}}^{2}$, is quadratic in these variables. Thus for small enough $\varepsilon$ the last logarithm in \eqref{Y2log} for $Y_{\sigma}$ is smooth as well, and the ratio between $|Y_{\sigma}^{2}|$ and $t(w,Y_{\sigma})$ is close to 1. Substituting these expressions into \eqref{Y2log}, where $|Y_{\kappa^{\perp}}^{2}|$ is bounded for $Y_{\kappa^{\perp}} \in M$, implies that $\log|Y^{2}|$ equals $\log|Y_{\sigma}^{2}|$ plus a smooth function on $A_{I,J}^{T}$ (as in Lemma \ref{logbound}), so that up to bounded functions, we are interested in $\log t+\log(w,Y_{\sigma})$. But as the positive linear combination $(w,Y_{\sigma})$ is bounded away from 0 in $A_{J}^{T}$, and the coefficients are all $O(\varepsilon t)$, this proves the lemma.
\end{proof}

Our analogue of Lemma \ref{GreenPhiK} is the following one, in which we recall the 1-periodic Bernoulli function $\mathbb{B}_{2}$.
\begin{lem} \label{GreenPhiD}
Let $C_{J}$ be the closure in $C$ of the union of the Weyl chambers whose closures contain $\omega_{J,I}$. Then, as a function of $Z_{L}\in\mathcal{D}$ with $Y=\Im Z_{L}$ lying in $C_{J}$ and satisfying the inequalities $|Y^{2}|>2m$ and $(w,Y)<\frac{|Y|}{2\sqrt{m}}$, the expression \[\sum_{\substack{\beta\in[J^{\perp}_{L^{*}}\cap(\mu+L)]/J \\ \beta^{2}=2m}}4\pi(w,Y)\mathbb{B}_{2}\bigg(\frac{(\beta,Y)}{(w,Y)}\bigg)-4\sum_{\substack{\nu\in[J^{\perp}_{L^{*}}\cap(\pm\mu+L)]/I \\ \nu^{2}=2m \\ (\nu+I_{\mathbb{R}},W)>0}}\log\big|1-\mathbf{e}\big((\nu,Z_{L})\big)\big|\] is smooth up to a logarithmic singularity along $2Z(m,\mu)$.
\end{lem}

\begin{proof}
The first inequality allows us to invoke Lemma \ref{GreenPhiK}, implying that the function $\widetilde{\psi}_{m,\mu}^{L}(Z_{L})$ has the required properties. Moreover, the second inequality allows us to express the first term in the expression for this function via Proposition \ref{PhimmuK}. This gives the first asserted sum, plus expressions that are smooth on $C_{J}$. The sum from Lemma \ref{GreenPhiK} gives the required sum, plus summands that are associated with vectors $\nu$ such that $(\nu,w)\neq0$. Since for such $\nu$ the sign of $(\nu+I_{\mathbb{R}},W)$ is the same for all $W \subseteq C_{J}$, all the corresponding summands are also smooth on $Z_{L}$ with $Y \in C_{J}$. This proves the lemma.
\end{proof}

We can now determine the behavior of of the automorphic Green function on the remaining domains.
\begin{prop} \label{PhinearJ}
For $M_{J}$ and $T$ as in Lemma \ref{loglin}, and for $Z_{L}+\mathbb{Z}\omega_{J,I}\in\mathbb{Z}\omega_{J,I}\backslash\mathcal{D}$ with $Y=\Im Z_{L}$ in $A_{J}^{T}$, the value of $\Phi_{m,\mu}^{L}(Z_{L}+\mathbb{Z}\omega_{J,I})$ is
\begin{equation} \label{expPhiJ}
t\Phi_{m,p_{D}\mu}^{D}-b_{m,\mu}^{+}(0,0)\log t-2\log|\phi|^{2}
\end{equation}
plus a function that is smooth on all of $\mathcal{D}_{\Upsilon}$, where $\phi$ is once again a holomorphic function defining $Z(m,\mu)$ in this set. In addition, consider $\sigma$, $\kappa$, $\varepsilon$, $k$, $M$, and $T$ as in that lemma. Then for $Z_{L}+K \in K\backslash\mathcal{D}$ with $Y \in A_{I,J}^{T}$, the formula from \eqref{expPhiinn} for $\Phi_{m,\mu}^{L}(Z_{L}+K)$ is still valid, with the remainder being smooth on all of $\mathcal{D}_{\Xi,\sigma}^{\varepsilon,k}$. Moreover, the difference between the first term there and that of \eqref{expPhiJ} is bounded by a constant multiple of $\varepsilon t$, and the difference between the second respective terms is bounded by $b_{m,\mu}^{+}(0,0)\log(\varepsilon t)$.
\end{prop}

\begin{proof}
In both cases Lemma \ref{loglin} allows us to express the $\Phi_{m,\mu}^{L}$-value, when $T$ is large enough, via Proposition \ref{PhimmuL}. Moreover, the terms involving $\mathcal{V}_{1+\frac{n}{2}}$ and the constant were seen in the proof of Proposition \ref{Phinearomega} to be smooth on all of $\mathcal{D}_{\Xi,\sigma}$, with the former decaying at the boundary. Now, since we work with a Weyl chamber $W$ containing $\omega_{J,I}$ in its closure, all the vectors $\nu$ from Proposition \ref{PhimmuL} satisfy $(\nu+I_{\mathbb{R}},\omega_{J,I})\geq0$. Moreover, if $\nu$ satisfies $(\nu+I_{\mathbb{R}},\omega_{J,I})>0$, then for large $T$ and small $\varepsilon$ we get $(\nu,Y)>0$ for every relevant $Y$, regardless of $W$, and as $t\to\infty$ all these terms are smooth on our domain ($\mathcal{D}_{\Upsilon}$ or $\mathcal{D}_{\Xi,\sigma}^{\varepsilon,k}$) and vanish on its boundary.

This covers all the summands in the penultimate term from Proposition \ref{PhimmuL}, except for those that are associated with $J_{L^{*}}/I$ with $J_{L^{*}}$ from \eqref{JDdef}. But these project modulo $I_{L^{*}}$ onto positive multiples of $w$, so that the corresponding summands produce constant coefficients of $F_{m,\mu}$ times expressions of the form $\log\big|1-\mathbf{e}\big(h(w,Z_{L})\big)\big|$ for positive $h\in\mathbb{Q}$ with bounded denominators. Since the latter are also smooth functions of $Z_{L}+\mathbb{Z}\omega_{J,I}$ or $Z_{L}+K$ and are independent of the coordinate $t$, they are smooth on $\mathcal{D}_{\Upsilon}$ and $\mathcal{D}_{\Xi,\sigma}^{\varepsilon,k}$. From the first sum in that proposition, only the summands from $J^{\perp}_{L^{*}}/I$ remain.

Now, Lemma \ref{loglin} transforms the logarithmic term in Proposition \ref{PhimmuL} into the second term in \eqref{expPhiJ} on $\mathcal{D}_{\Upsilon}$, and establishes the statement about the difference from the second term in \eqref{expPhiinn}. Moreover, the second inequality in that lemma implies that we can express the first term from Proposition \ref{PhimmuL} via Proposition \ref{PhimmuK}. As the arguments of the Bessel $K$-functions were seen in that lemma to grow towards the boundary, these terms are smooth on $\mathcal{D}_{\Upsilon}$ and on $\mathcal{D}_{\Xi,\sigma}^{\varepsilon,k}$, and decay at the boundary. Moreover, the terms involving the Bernoulli function combine with the sum over the elements from $J^{\perp}_{L^{*}}/I$ to yield the expression from Lemma \ref{GreenPhiD}, which extends to $\mathcal{D}_{\Upsilon}$ or $\mathcal{D}_{\Xi,\sigma}^{\varepsilon,k}$ and gives the last term of \eqref{expPhiJ} as in the proof of Proposition \ref{Phinearomega}.

Next, the second term from Proposition \ref{PhimmuK} is a constant times $(w,Y)$, which is a smooth and bounded function on $A_{J}^{T}$ and therefore also on $\mathcal{D}_{\Upsilon}$. Moreover, Lemma \ref{loglin} expresses the coefficient in front of $\Phi_{m,p_{D}\mu}^{D}$ as $t-\big|Y_{\kappa_{J}^{\perp}}^{2}\big|\big/2(w,Y)$, producing the first term of \eqref{expPhiJ}, with the quotient being again smooth and bounded on $A_{J}^{T}$ thus on $\mathcal{D}_{\Upsilon}$. On the other hand, on $A_{J,I}^{T}$ both $(w,Y)$ and this quotient were seen in the proof of that lemma to be $O(\varepsilon t)$, yielding the statement about the difference between the first terms. This completes the proof of the proposition.
\end{proof}
As in Remark \ref{indepW}, the function $\phi$ in Proposition \ref{PhinearJ} also depends on the choice of a Weyl chamber $W$, and the effect of changing $W$ just goes into the smooth remainder.

\smallskip

We would like to view the function $\Phi_{m,\mu}^{L}$ in the setting of \cite{[BBK]}. As this paper concerns smooth varieties, we shall consider smooth toroidal compactifications. The fact that $X_{L}$ itself may be not smooth can be dealt with by taking covers as in Remark \ref{congGamma}. For smoothness at the boundary, we take a refinement as in Remark \ref{smexist}, and the resulting variety $X_{L}^{\mathrm{tor}}$ is then smooth up to finite quotient singularities from $\Gamma_{L}$, with the toroidal boundary being a divisor with normal crossings. In this case our main result about $\Phi_{m,\mu}^{L}$, namely Theorem \ref{GreenPhi} below, can also be described using terminology from \cite{[BBK]}---see Remark \ref{logBBK} below.

Indeed, there are some simplifications when a fan $\sigma\in\Sigma(\Xi)$ is smooth in the sense of Definition \ref{smoothdef}. As the vectors $\omega_{j}$, $1 \leq j \leq d$ generate over $\mathbb{Z}$ the primitive sublattice $K\cap\mathbb{R}\sigma$ of $K$, we can take elements $\kappa_{i} \in K^{*}$ with $1 \leq i \leq d$ such that $(\kappa_{i},\omega_{j})$ is the Kronecker $\delta$-symbol $\delta_{ij}$, and such that the space $\kappa \subseteq K_{\mathbb{R}}$ that they span contains negative norm vectors. Note that we cannot always take such a dual $\mathbb{Z}$-basis inside $\mathbb{R}\sigma$, which is why $\kappa$ can be a different rational subspace of $K_{\mathbb{R}}$. We consider $\kappa^{\perp} \subseteq K_{\mathbb{R}}$ and the projections as before, and then in our decomposition of $Z \in T_{K}(\mathbb{C})$ as $Z_{\sigma}+Z_{\kappa^{\perp}}$ we have
\begin{equation} \label{cjdef}
Z_{\kappa^{\perp}}\in[\kappa^{\perp}_{\mathbb{C}}/(\kappa^{\perp} \cap K)]\quad\mathrm{and}\quad Z_{\sigma}=\sum_{j=1}^{d}c_{j}\omega_{j}\in[\mathbb{C}\sigma/(\sigma \cap K)]
\end{equation}
with uniquely defined $c_{j}\in\mathbb{C}/\mathbb{Z}$ for $1 \leq j \leq d$, and we write each such $c_{j}$ as $a_{j}+ib_{j}$ with $a_{j}\in\mathbb{R}/\mathbb{Z}$ and $b_{j}\in\mathbb{R}$. The smoothness condition is required for the vectors $\{\omega_{j}\}_{j=1}^{d}$ to span $[\mathbb{C}\sigma/(\sigma \cap K)]$ with free coefficients from $\mathbb{C}/\mathbb{Z}$, without additional quotients. Similarly, assume that $I \subseteq J \subseteq I^{\perp}_{L}$ with $D$ from \eqref{JDdef} and $\Upsilon$ as before, and take some basis $(z,w)$ for $J$ with $z \in I$ as well as a basis $(\zeta,\upsilon)$ for some unimodular complement $\tilde{J}$ for $J$ in $K^{*}$ that is dual to $(z,w)$. Assume that $\sigma$ contains the ray $\rho=\mathbb{R}_{+}\omega_{J,I}=\Omega(\Upsilon)\cup\{0\}$, or that $\sigma=\rho$ and we work in $\mathcal{D}_{\Upsilon}$. Then the dual basis $(\zeta,\upsilon)$ produces, as in \cite{[K3]}, the coordinates
\begin{equation} \label{coorJ}
Z_{L}=\zeta+\tau\upsilon+Z_{0}-\eta w+\big(\eta\tau-\tfrac{Z_{0}^{2}+(\zeta+\tau\upsilon)^{2}}{2}\big)z
\end{equation}
for $Z_{L} \in P_{+}$ that lies in the image of the section associated with $z$. Here $\tau\in\mathcal{H}$ is the coordinate from \eqref{secwithz}, $Z_{0}=X_{0}+iY_{0} \in D_{\mathbb{C}}$, and $\eta=s+it\in\mathbb{C}$ is another coordinate, which satisfies $t>\frac{Y_{0}^{2}}{2y}$ for $\Im Z_{L}$ to be of negative norm (this also gives a more explicit, though less canonical, description of Proposition \ref{struc1dim}). The cyclic group $\mathcal{U}_{Z}(\Upsilon)$ operates as translations of the real part of $\sigma$ by integers, and we denote, following \cite{[K3]}, the resulting coordinate $\mathbf{e}(\eta)$ on $\mathcal{U}_{\mathbb{Z}}(\Upsilon) \backslash \operatorname{Gr}(L)$ by $q_{2}$. Considering the imaginary part $Y$ of $Z$ or of $Z_{L}$, if $Z$ is presented as in \eqref{cjdef} then $Y_{\sigma}$ can be written uniquely as $\sum_{j=1}^{d}b_{j}\omega_{j}$. For $Z_{L}$ as in \eqref{coorJ} we can write $Y$ as $y\upsilon+Y_{0}-tw$, or equivalently $y\upsilon+Y_{0}+t\omega_{J,I}$, with $Y_{0}=\Im Z_{0} \in D_{\mathbb{R}}$. We recall from \eqref{znudef} that the function $z_{\kappa_{j}}$ on $\mathcal{D}_{\Xi,\sigma} \subseteq X_{K,\sigma}$ is just $\mathbf{e}(c_{j})$ for $Z_{L}$ as in \eqref{cjdef}, and from Lemma \ref{DXidim1} that $q_{2}=\mathbf{e}(\eta)$ is the local defining equation $z_{\nu_{\Upsilon}}$ for the complement of $\mathcal{U}_{\mathbb{Z}}(\Upsilon) \backslash \operatorname{Gr}(L)$ in $\mathcal{D}_{\Upsilon}$. We deduce that
\begin{equation} \label{bjtlogcoor}
b_{j}=\Im c_{j}=\frac{1}{2\pi}\log\frac{1}{|z_{\kappa_{j}}|}\mathrm{\ for\ }1 \leq j \leq d,\quad\mathrm{and}\quad t=\Im\eta=\frac{1}{2\pi}\log\frac{1}{|q_{2}|}.
\end{equation}
We remark again that if $\rho\leq\sigma$ for $\sigma\neq\rho$, with $\omega_{J,I}=\omega_{k}$ for $1 \leq k \leq d$, then the two descriptions from \eqref{bjtlogcoor} are related via the fact that $t-b_{k}$ is a linear combination of $\{b_{j}\}_{j \neq k}$, so that $q_{2}=z_{\nu_{\Upsilon}}$ is $z_{\kappa_{k}}$ times a toric function that is defined and does not vanish on the boundary of $\mathcal{D}_{\Upsilon}$.

We can now state and prove the relation to \cite{[BBK]} in the smooth case.
\begin{cor} \label{prelog}
If $X_{L}^{\mathrm{tor}}$ is smooth then the function $\Phi_{m,\mu}^{L}$ is pre-$\log$ along $Z(m,\mu)$ and along the boundary divisor $X_{L}^{\mathrm{tor}} \setminus X_{L}$ in the sense of Definition~1.5 of \cite{[BBK]}, and has a component that is pre-$\log$-$\log$ along the boundary.
\end{cor}

\begin{proof}
Recalling that $X_{L}^{\mathrm{tor}}$ is covered by the images of the sets of the form $\mathcal{D}_{\Xi,\sigma}$ or $\mathcal{D}_{\Upsilon}$, it suffices to consider the behavior of $\Phi_{m,\mu}^{L}$ on these sets. We use the description from Propositions \ref{Phinearomega} and \ref{PhinearJ}, in which we substitute the expression for $Y_{\sigma}$ as $\sum_{j=1}^{d}b_{j}\omega_{j}$ into \eqref{expPhiinn}, and apply \eqref{bjtlogcoor} for the coordinates $b_{j}$ there as well as for the coordinate $t$ in \eqref{expPhiJ}. Now, in the closure of $\mathcal{D}_{\Xi,\sigma}^{\varepsilon,0}$ the argument of $\Phi_{m,p_{K}\mu}^{K}$ in \eqref{expPhiinn} lies in a compact subset of $C$ and is thus bounded, and $|Y_{\sigma}|$ is the square root of
\begin{equation} \label{Ysigma2}
|Y_{\sigma}^{2}|=-Y_{\sigma}^{2}=
-\sum_{i=1}^{d}\sum_{j=1}^{d}(\omega_{i},\omega_{j})b_{i}b_{j}=\sum_{i=1}^{d}\sum_{j=1}^{d}\frac{\big|(\omega_{i},\omega_{j})\big|}{4\pi^{2}}\log\frac{1}{|z_{\kappa_{i}}|}\log\frac{1}{|z_{\kappa_{j}}|},
\end{equation}
with all the coefficients $-(\omega_{i},\omega_{j})$ being strictly positive. Noting that \eqref{Ysigma2} reduces to $\frac{|\omega_{1}|^{2}}{4\pi^{2}}\log^{2}\frac{1}{|z_{\kappa_{1}}|}$ when $d=1$ and is bounded by a multiple of $\sum_{j=1}^{d}\log^{2}\frac{1}{|z_{\kappa_{j}}|}$ in general, and that the smoothness of $\sigma$ implies that the $z_{\kappa_{j}}$s are adapted to the boundary divisor of $X_{K,\sigma}$, this gives the $\log$-growth from Definition 1.4 of \cite{[BBK]} (with $M=1$) of the first term in \eqref{expPhiinn} as well as the $\log$-$\log$-growth from Definition 1.2 of that reference (with $M=2$) of the second term there, and the third term clearly has the $\log$-growth along $Z(m,\mu)$ (again with $M=1$). The fact that the second term satisfies the condition to be pre-$\log$-$\log$ from Definition 1.3 there is now clear, and recalling that by appropriately normalizing the defining function $\phi$ of $Z(m,\mu)$ we can replace $\Phi_{m,p_{K}\mu}^{K}$ by its smooth analog $\Phi_{m,p_{K}\mu}^{K,W}$ associated with a Weyl chamber $W$ in $C$, the conditions for pre-$\log$ in Definition 1.5 of that reference are seen to be satisfied by the first and third terms of that equation. In addition, substituting the value of $t$ from \eqref{bjtlogcoor} into \eqref{expPhiJ} immediately establishes the required properties of $\Phi_{m,\mu}^{L}$ along $\mathcal{D}_{\Upsilon}$, and the bounds from Proposition \ref{PhinearJ} on the difference between \eqref{expPhiinn} and \eqref{expPhiJ} on $\mathcal{D}_{\Xi,\sigma}^{\varepsilon,k}$ combine with this value of $t$ to establish the desired properties also on this set. This completes the proof of the corollary.
\end{proof}
In fact, the proof of Corollary \ref{prelog} establishes the linear logarithmic growth of $\Phi_{m,\mu}^{L}$ and the $\log$-$\log$-growth of the corresponding term also when $X_{L}^{\mathrm{tor}}$ is not smooth, but expressing this in coordinates becomes more complicated.

\smallskip

For the next important property of $\Phi_{m,\mu}^{L}$ we shall need the following technical lemma.
\begin{lem} \label{intoflog}
Take $1 \leq d\in\mathbb{N}$, $d<\alpha\in\mathbb{R}$, and $0<\varepsilon_{j}<1$ for every $1 \leq j \leq d$. We then have \[\int_{0}^{\varepsilon_{d}}\cdots\int_{0}^{\varepsilon_{1}}\frac{dr}{\prod_{j=1}^{d}r_{j}\cdot\log^{\alpha}\prod_{j=1}^{d}\frac{1}{r_{j}}}<\infty\quad\mathrm{and}\quad
\int_{0}^{\varepsilon_{d}}\cdots\int_{0}^{\varepsilon_{1}}\frac{\log\log\prod_{j=1}^{d}\frac{1}{r_{j}} \cdot dr}{\prod_{j=1}^{d}r_{j}\cdot\log^{\alpha}\prod_{j=1}^{d}\frac{1}{r_{j}}}<\infty,\] where $dr$ is the standard Lebesgue measure on $\mathbb{R}^{d}$.
\end{lem}

\begin{proof}
For $d=1$ a direct evaluation gives the finite values $\frac{1}{(\alpha-1)\log^{\alpha-1}\frac{1}{\varepsilon}}$ and $\frac{\log\log\frac{1}{\varepsilon}+\frac{1}{\alpha-1}}{(\alpha-1)\log^{\alpha-1}\frac{1}{\varepsilon}}$ respectively, with $\varepsilon=\varepsilon_{1}$. If $d>1$ then we set $R=\prod_{j=1}^{r_{j}}$ inside the internal integral, which transforms it into $\prod_{j=1}^{d-1}\frac{1}{r_{j}}$ times the 1-dimensional integral associated with $\varepsilon=\varepsilon_{d}\prod_{j=1}^{d-1}r_{j}$. Substituting the explicit values from the case $d=1$, and then replacing $r_{d-1}$ by $\varepsilon_{d}r_{d-1}$, transforms the integral of the first form into $\frac{1}{\alpha-1}$ times the one with $d-1$, $\alpha-1$, and last parameter $\varepsilon_{d-1}\varepsilon_{d}$. For the integral of the second form this gives $\frac{1}{\alpha-1}$ times the integral of the second form plus $\frac{1}{(\alpha-1)^{2}}$ times the one of the first form, both with the same parameters again. A simple induction on $d$ thus yields the finiteness of both integrals. This proves the lemma.
\end{proof}

We can now prove the following result, which holds also in the non-smooth setting.
\begin{prop} \label{locint}
If $n\geq2$ then $\Phi_{m,\mu}^{L}$ is locally integrable on any set of the form $\mathcal{D}_{\Xi,\sigma}^{\varepsilon,0}$. In case $n\geq3$ it is also locally integrable on each $\mathcal{D}_{\Upsilon}$ and on every $\mathcal{D}_{\Xi,\sigma}^{\varepsilon,k}$ with $k\geq1$.
\end{prop}

\begin{proof}
We recall from Section 4.1 of \cite{[Br1]} that the invariant measure on $\operatorname{Gr}(L)$ presented as in \eqref{GVRKRiC}, as well as on its quotient modulo $K$, is a multiple of $\frac{dXdY}{|Y^{2}|^{n}}$. The local integrability of $\Phi_{m,\mu}^{L}$ on the open subset $K\backslash\operatorname{Gr}(L)$ follows from the fact that there it is smooth up to logarithmic singularities. In addition, since all the boundaries in question have measure 0, the local integrability does not depend on the choice of the toroidal compactification, and it therefore suffices to consider the case where each fan $\Sigma(\Xi)$ is smooth as in Definition \ref{smoothdef}, where we can use the expressions from the proof of Corollary \ref{prelog}.

Consider next a point in $\xi\in\mathcal{D}_{\Xi,\sigma}^{\varepsilon,0}$ for some non-trivial cone $\sigma\in\Sigma(\Xi)$, which is therefore not an external ray in $C$, and the argument preceding \eqref{sigmapar} shows that it suffices to consider $\xi \in O(\sigma)$. Our choice of the space $\kappa$ as above allows us to take a neighborhood of $\xi$ which is the product of a neighborhood of the projection $\xi_{\kappa^{\perp}}$ in $\kappa^{\perp}_{\mathbb{C}}/(\kappa^{\perp} \cap K)$ and a neighborhood of $\xi_{\sigma}$ in the completion of $\mathbb{C}\sigma/(\sigma \cap K)$, with both having compact closures in the respective spaces. The smoothness of $\sigma$ allows one to use the coordinates from \eqref{cjdef}, and the relation between $c_{j}=a_{j}+ib_{j}$ and $z_{\kappa_{j}}$, which is visible in \eqref{bjtlogcoor}, implies that each factor $da_{j}db_{j}$ in $dXdY$ becomes $\frac{1}{4\pi^{2}}$ times $\frac{dz_{\kappa_{j}}d\overline{z_{\kappa_{j}}}}{-2i|z_{\kappa_{j}}|^{2}}$. The coordinates arising from $\xi_{\kappa^{\perp}}$ can be expressed similarly, thus giving a smooth bounded $(2n-2d)$-form on our neighborhood, because the coordinates that do not come from $\sigma$ do not vanish on $\mathcal{D}_{\Xi,\sigma}^{\varepsilon,0}$. We can also assume that the neighborhood of $\xi_{\sigma}$, on which all the $z_{\kappa_{j}}$s vanish, is a product of balls, one for each coordinate, and when we consider the $j$th multiplier in polar coordinates, we see that the $2d$-form $da_{j}db_{j}$ is $\frac{1}{(4\pi^{2})^{d}}$ times $\prod_{j=1}^{d}\frac{dr_{j}d\theta_{j}}{r_{j}}$, with $0\leq\theta_{j}<2\pi$ and $0<r_{j}<\varepsilon_{j}$ for every $j$. As for the denominator $|Y^{2}|^{n}$, the proof of Proposition \ref{Phinearomega} allows us to write it as $|Y_{\sigma}^{2}|^{n}$ times a function that is smooth and bounded on our neighborhood of $\xi$, and since \eqref{Ysigma2} shows that the latter expression is quadratic in $\log\frac{1}{|z_{\kappa_{j}}|}=\log\frac{1}{r_{j}}$ with $1 \leq j \leq d$ and positive wherever each one of them is positive, it is comparable with any other quadratic form with these properties, such as $\big(\sum_{j=1}^{d}\log\frac{1}{r_{j}}\big)^{2}$, in these expressions. Altogether, we have written our invariant measure on neighborhoods of points in $O(\sigma)$ as $\prod_{j=1}^{d}\frac{dr_{j}}{r_{j}}$ times $\log^{-2n}\prod_{j=1}^{d}\frac{1}{r_{j}}$ on $\prod_{j=1}^{d}[0,\varepsilon_{j})$, times a smooth and bounded $(2n-d)$-form with no singularities.

We therefore need to establish the integrability of $\Phi_{m,\mu}^{L}$ in this neighborhood with respect to this measure, and we express it via \eqref{expPhiinn}. For the smooth part of this function, as well as the term with $\log|\phi|$, it suffices to consider the integrability of the measure itself, which now follows from the case of the first integral with $\alpha=2n \geq 2d>d$ in Lemma \ref{intoflog}. The integral of the second term of \eqref{expPhiinn} is just the second integral in Lemma \ref{intoflog} with $\alpha=2n$, while in the first term, the fact that the argument of $\widetilde{\Phi}_{m,p_{K}\mu}^{K}$ lies in a compact subset of $C$ when $Z_{L} \in \mathcal{D}_{\Xi,\sigma}^{\varepsilon,0}$ yields the boundedness of that multiplier. Since the multiplier $|Y_{\sigma}|$ simply replaces the power $2n$ in the denominator by $2n-1$, and we have $2n-1 \geq d+n-1>d$ because $n\geq2$, Lemma \ref{intoflog} establishes the integrability of this term as well. This proves the first assertion.

For the second one we take a boundary point $\xi$ of either $\mathcal{D}_{\Upsilon}$ or $\mathcal{D}_{\Xi,\sigma}^{\varepsilon,k}$ for a cone $\sigma\in\Sigma(\Xi)$ with $\rho=\mathbb{R}_{+}\omega_{J,I}=\Omega(\Upsilon)\cup\{0\}\leq\sigma$ and $\sigma\neq\rho$, where in the latter case we can again assume that $\xi \in O(\sigma)$. Recalling the coordinate $\eta=s+it$ from \eqref{coorJ}, the same argument shows that on a similar product neighborhood, the part $dXdY$ of the invariant measure is again $\prod_{j=1}^{d}\frac{dr_{j}}{r_{j}}$ times a smooth $(2n-d)$-form with no singularities, where $r=r_{k}=|q_{2}|$ and the other $r_{j}$s come from the other rays of $\sigma$ in case $\sigma\neq\rho$. We now express $|Y^{2}|$ as in \eqref{Y2JI}, and the proof of Proposition \ref{PhinearJ} shows that if $\sigma=\rho$ then it is a smooth and bounded function times $\log\frac{1}{r}$ while when $\sigma\neq\rho$ it can be written as a smooth and bounded function times $t(w,Y_{\sigma,\kappa_{J}^{\perp}})$. As the second multiplier is a linear form with positive coefficients in $\log\frac{1}{r_{j}}$ with $j \neq k$, it is a smooth function times the linear form $\sum_{j \neq k}\log\frac{1}{r_{j}}$, and our invariant measure here is the product of $\log^{-n}\frac{1}{r}\cdot\frac{dr}{r}$ on $[0,\varepsilon)$ for $\varepsilon=\varepsilon_{k}$, of $\prod_{j=1}^{d}\frac{dr_{j}}{r_{j}}$ times $\log^{-n}\prod_{j \neq k}\frac{1}{r_{j}}$ on $\prod_{j \neq k}[0,\varepsilon_{j})$ when $\sigma\neq\rho$, and of a smooth, bounded, non-singular $(2n-d)$-form.

For the function $\Phi_{m,\mu}^{L}$ we now use the expression from \eqref{expPhiJ}, and recall from Proposition \ref{PhinearJ} that while when $\sigma\neq\rho$ we still need to use \eqref{expPhiinn}, the difference is bounded by a constant times $t$ and can thus be considered like the first term in \eqref{expPhiJ}. Since the only parameter appearing in \eqref{expPhiJ} is $t$, or equivalently $\log\frac{1}{r}$ via \eqref{bjtlogcoor}, the integrability of the $(d-1)$-form arising from $r_{j}$ with $j \neq k$ follows from Lemma \ref{intoflog} because $n \geq d$ and the dimension here is $d-1$. The integrability of the smooth part of $\Phi_{m,\mu}^{L}$ and of $\log|\phi|$, follow from the case $d=1$ and $\alpha=n\geq2$ in the first integral in Lemma \ref{intoflog}, and that of the second term from the second integral with $d=1$ and $\alpha=n$ there. As for the first term we obtain the first integral with $d=1$ and $\alpha=n-1$, and that lemma is again applicable since we assume that $n\geq3$, the second assertion is also established. This proves the proposition.
\end{proof}

\smallskip

Note that Proposition \ref{Phinearomega} becomes simpler when $d=1$. Then $\sigma\in\Sigma(\Xi)$ is the ray $\rho=\mathbb{R}_{+}\omega$ with $\omega$ primitive in $K \cap C$, the boundary points lie in $O(\rho) \subseteq
X_{K,\rho}$, and if $\kappa \in K^{*}$ satisfies $(\kappa,\omega)=1$ then Proposition \ref{ordTdiv} implies that the function $z_{\kappa}$ is the defining equation for $O(\rho)$ in $X_{K,\rho}$. Then the vector $Y_{\sigma}$ appearing in \eqref{expPhiinn} is just $b\omega$ with $b=b_{1}$, so that $|Y_{\sigma}|=b|\omega|$, and it is clear that $\mathcal{D}_{\Xi,\rho}^{\varepsilon,0}=\mathcal{D}_{\Xi,\rho}$ because $\rho$ has no faces that are rays lying on the boundary of $C$. This reduces \eqref{expPhiinn} to
\begin{equation} \label{expPhiomega}
\frac{b|\omega|}{\sqrt{2}}\Phi_{m,p_{K}\mu}^{K}\bigg(\frac{\omega}{|\omega|}\bigg)-2b_{m,\mu}^{+}(0,0)\log b-2\log|\phi|^{2}
\end{equation}
up to an additive constant, and we also recall that $O(\sigma)$ and the subset $O(\Upsilon)\subseteq\mathcal{D}_{\Upsilon}$ from Lemma \ref{DXidim1} map onto the generic part of the divisors $B_{I,\omega}$ and $B_{J}$ from \eqref{BJBIomega} respectively. Combining these observations with \eqref{bjtlogcoor}, we make the following definition.
\begin{defn} \label{mults}
For a rank 1 primitive isotropic sublattice $I$ of $L$, with the associated lattice $K$ from \eqref{IKdef} and cone $C$ as in \eqref{GVRKRiC}, and a primitive element $\omega \in K \cap C$, we define the multiplicity \[\operatorname{mult}_{I,\omega}(m,\mu)=\frac{|\omega|}{8\sqrt{2}\pi}\Phi_{m,p_{K}\mu}^{K}\bigg(\frac{\omega}{|\omega|}\bigg).\] Given a rank 2 primitive isotropic sublattice $J$ in $L$, set $D$ to be as in \eqref{JDdef}, and we set \[\operatorname{mult}_{J}(m,\mu)=\frac{1}{8\pi}\Phi_{m,p_{D}\mu}^{D}.\]
\end{defn}
In fact, Proposition \ref{PhiDeval} below will show that the multiplicities $\operatorname{mult}_{J}(m,\mu)$ from Definition \ref{mults} are always rational numbers. On the other hand, the multiplicities $\operatorname{mult}_{I,\omega}(m,\mu)$ are real numbers, about the rationality of which we have no information in general.

\smallskip

Consider again the variety $X_{L}^{\mathrm{tor}}=X^{\mathrm{tor}}_{\{\Sigma(\Xi)\}_{\Xi}}$ from Theorem \ref{formtc}, which is defined by a collection of fans $\{\Sigma(\Xi)\}_{\Xi}$ that satisfies the admissibility conditions from Corollary \ref{admcdecom} and Lemma \ref{GammaL0dim} with respect to $\Gamma_{L}$. Using Definition \ref{mults}, the required extension of Theorem \ref{PhiGreen} can now be established.
\begin{thm} \label{GreenPhi}
For every pair $(m,\mu)$ with $\mu\in\Delta_{L}$ and $0<m\in\frac{\mu^{2}}{2}+\mathbb{Z}$ we define, as in \eqref{ZtorInt}, the divisor $Z^{\mathrm{tor}}(m,\mu)\in\operatorname{Div}(X_{L}^{\mathrm{tor}})\otimes\mathbb{R}$ by \[Z^{\mathrm{tor}}(m,\mu)=Z(m,\mu)+\sum_{\substack{J\subseteq L \\ \mu \perp J_{L^{*}}/J}}\operatorname{mult}_{J}(m,\mu) \cdot B_{J}+\sum_{\substack{I\subseteq L \\ \mu \perp I_{L^{*}}/I}}\sum_{\mathbb{R}_{+}\omega}\operatorname{mult}_{I,\omega}(m,\mu) \cdot B_{I,\omega}.\] Here $J$ (resp.~$I$) runs over a set of representatives of rank 2 (resp.~rank 1) primitive isotropic sublattices of $L$ modulo $\Gamma_{L}$ and given $I$, the index $\mathbb{R}_{+}\omega$ runs over representatives for the inner rays in $\Sigma(\Xi)$ with $\Xi=\Xi_{I}$, and the coefficients $\operatorname{mult}_{J}(m,\mu)$ and $\operatorname{mult}_{I,\omega}(m,\mu)$ are the ones from Definition \ref{mults}. Assuming that $n$ is larger than the Witt rank of $L$, the function $\frac{1}{2}\Phi_{m,\mu}^{L}$ is integrable on $X_{L}^{\mathrm{tor}}$, grows logarithmically along the divisor $Z^{\mathrm{tor}}(m,\mu)$ with a term with $\log$-$\log$-growth along the boundary, and satisfies the $dd^{c}$-equation \eqref{ddceq}.
\end{thm}

\begin{proof}
First, as the invariant measure is smooth on $X_{L}$ and our function has only the logarithmic singularities along $Z(m,\mu)$, it is locally integrable on $X_{L}$, which also covers the case where $L$ has Witt rank 0 since then there is no boundary. For the local integrability at boundary points, we recall that for each such point one of the maps $\iota_{\Xi}$ and $\iota_{\Upsilon}$ appearing in \eqref{BJBIomega} is local homeomorphism, so that it suffices to consider local integrability on the sets $\mathcal{D}_{\Xi,\sigma}$ or $\mathcal{D}_{\Upsilon}$. Now, if $L$ has Witt rank 1 then there are no rays of the form $\mathbb{R}_{+}\omega_{J,I}$, so that $\mathcal{D}_{\Xi,\sigma}^{\varepsilon,0}=\mathcal{D}_{\Xi,\sigma}$ for every $\sigma\in\sigma(\Xi)$, and as the boundary of $X_{L}^{\mathrm{tor}}$ is then covered by such sets, the assumption that $n\geq2$ yields the desired local integrability by Proposition \ref{locint}. On the other hand, when $L$ has Witt rank 2 then $n\geq3$, and we also have local integrability on the remaining covering sets $\mathcal{D}_{\Upsilon}$ by the same proposition. This proves the local integrability of $\frac{1}{2}\Phi_{m,\mu}^{L}$ on $X_{L}^{\mathrm{tor}}$, and the global one follows from compactness.

Now, the assertion about the growth of $\frac{1}{2}\Phi_{m,\mu}^{L}$ was established in Corollary \ref{prelog}, and the integrability just proved allows us to substitute $\frac{1}{2}\Phi_{m,\mu}^{L}$ inside current equations. Moreover, the fact that the logarithmic growth is linear, namely we had $M=1$ in the proof of Corollary \ref{prelog} in the notation of Definition 1.4 of \cite{[BBK]}, implies that it satisfies the desired $dd^{c}$-equation from \eqref{ddceq} for some divisor on $X_{L}^{\mathrm{tor}}$. Since the required properties of $\eta(m,\mu)$ are also clear, all that remains it to determine this divisor to be $Z^{\mathrm{tor}}(m,\mu)$.

Consider first a divisor of the form $B_{I,\omega}$ for some primitive vector $\omega \in K \cap C$, which exists only if $n\geq2$ in our assumption, and for which we denote $\rho=\mathbb{R}_{+}\omega\in\Sigma(\Xi)$. Then the case $d=1$ of Proposition \ref{Phinearomega} discussed above evaluates $\Phi_{m,\mu}^{L}$ on $\mathcal{D}_{\Xi,\rho}$ as \eqref{expPhiomega} plus a smooth function. Substituting \eqref{bjtlogcoor}, dividing by 2, and recalling the behavior of $-\log$ of the square of the absolute value of a coordinate under $dd^{c}$, establishes the equality \[dd^{c}[\tfrac{1}{2}\Phi_{m,\mu}^{L}]+\delta_{Z_{I,\omega}(m,\mu)}=[\eta_{I,\omega}(m,\mu)]\] as currents on $\mathcal{A}_{c}^{2n-2}(\mathcal{D}_{\Xi,\rho})$, where $\eta_{I,\omega}(m,\mu)$ is the sum of a smooth 2-form on $\mathcal{D}_{\Xi,\rho}$ and a multiple of $dd^{c}\log\log\frac{1}{|z_{\kappa}|}$, and
\begin{equation} \label{ZIomegammu}
Z_{I,\omega}(m,\mu)=Z(m,\mu)+\frac{|\omega|}{8\sqrt{2}\pi}\Phi_{m,p_{K}\mu}^{K}\bigg(\frac{\omega}{|\omega|}\bigg) \cdot O(\rho),
\end{equation}
where $O(\rho)$ is the pullback of $B_{I,\omega}$ to $\mathcal{D}_{\Xi,\rho}\subseteq\mathcal{D}_{\Xi}$ via $\iota_{\Xi}\circ\pi_{\Xi}$.

We now turn to the divisors $B_{J}$ under the assumption that $n\geq3$, where Proposition \ref{PhinearJ} describes $\Phi_{m,\mu}^{L}$ on $\mathcal{D}_{\Upsilon}$ via \eqref{expPhiJ} up to a smooth function. We again substitute \eqref{bjtlogcoor} and divide by 2, and we recall from Proposition \ref{torbcdim1} that the divisor $O(\Upsilon)$ inside $\mathcal{D}_{\Upsilon}$ is characterized by the vanishing of $q_{2}$. The same argument now produces the equality \[dd^{c}[\tfrac{1}{2}\Phi_{m,\mu}^{L}]+\delta_{Z_{J}(m,\mu)}=[\eta_{J}(m,\mu)]\] on $\mathcal{A}_{c}^{2n-2}(\mathcal{D}_{\Upsilon})$, in which $\eta_{J}(m,\mu)$ is the sum of a smooth 2-form on $\mathcal{D}_{\Upsilon}$ and a multiple of $dd^{c}\log\log\frac{1}{|q_{2}|}$, and
\begin{equation} \label{ZJmmu}
Z_{J}(m,\mu)=Z(m,\mu)+\frac{\Phi_{m,p_{D}\mu}^{D}}{8\pi} \cdot O(\Upsilon),
\end{equation}
with $O(\Upsilon)$ being the pullback of $B_{J}$ to $\mathcal{D}_{\Upsilon}$ under $\iota_{\Upsilon}\circ\pi_{\Upsilon}$

Now, the union of $X_{L}$ and the generic parts of the toroidal divisors has a complement in $X_{L}^{\mathrm{tor}}$ whose complex codimension there is 2, and therefore the restriction map from divisors on $X_{L}^{\mathrm{tor}}$ to divisors on this open subvariety is surjective. It follows that the divisor on $X_{L}^{\mathrm{tor}}$ which appears in the global $dd^{c}$-equation for $\frac{1}{2}\Phi_{m,\mu}^{L}$ has to have the pullback from \eqref{ZIomegammu} to $\mathcal{D}_{\Xi,\mathbb{R}_{+}\omega}$ for every $I$ and $\omega$, as well as the pullback given in \eqref{ZJmmu} to $\mathcal{D}_{\Upsilon}$ for every $J$. Since Definition \ref{mults} implies that the only divisor having this property is $Z^{\mathrm{tor}}(m,\mu)$, we have also established \eqref{ddceq} with the correct divisor. This completes the proof of the theorem.
\end{proof}
Note that in the case $n=1$ in Theorem \ref{GreenPhi} the assumption about the Witt rank of $L$ implies that there are no cusps, but this case is therefore covered in Theorem \ref{PhiGreen}.
\begin{rmk} \label{logBBK}
When $X_{L}^{\mathrm{tor}}$ is smooth, the result of Theorem \ref{GreenPhi} means that $\frac{1}{2}\Phi_{m,\mu}^{L}$ is an integrable logarithmic Green function for $Z^{\mathrm{tor}}(m,\mu)$ on $X_{L}^{\mathrm{tor}}$ with additional $log$-$log$-growth along the boundary divisor in the sense of Definition 1.12 of \cite{[BBK]}.
\end{rmk}

\subsection{Explicit Formulas for the Multiplicities}

The multiplicities $\operatorname{mult}_{J}(m,\mu)$ and $\operatorname{mult}_{I,\omega}(m,\mu)$ are given in Definition \ref{mults} as regularized theta lifts. However, it can be useful to have more explicit expressions for the values of these multiplicities.

\smallskip

We begin by assuming that $n\geq2$ and considering the divisor $B_{J}$ that is associated via \eqref{BJBIomega} with the rank 2 isotropic sublattice $J$ of $L$, where $D$ is the positive definite lattice from \eqref{JDdef} as usual. Note that $\Gamma_{L}$ preserves the condition $\mu \in J^{\perp}_{L^{*}}/J^{\perp}_{L}\subseteq\Delta_{L}$ (by definition) and the isomorphism class of $D$, so that the expressions below are independent of the choice of the lattice $J$ that represents $B_{J}$. We recall the theta functions $\Theta_{D}$ and $\uparrow^{L}_{D}(\Theta_{D})$ from \eqref{ThetaD}, set $\sigma_{1}(n)=\sum_{d|n}d$ for $1 \leq n\in\mathbb{N}$, denote by \[E_{2}^{*}(\tau)=E_{2}(\tau)-\frac{3}{\pi y}=1-24\sum_{n=1}^{\infty}\sigma_{1}(n)q^{n}-\frac{3}{\pi y}\] the nearly holomorphic Eisenstein series of weight 2 (with $E_{2}$ being the weight 2 holomorphic quasi-modular Eisenstein series), and introduce the notation $\operatorname{CT}(G)$ for the constant term of a Laurent series $G$ in the variable $q$. Recalling that the arrow operators commute with differential operators and are dual with respect to the pairings via
\begin{equation} \label{arrowpar}
\langle a,\downarrow^{L}_{D}b \rangle_{D}=\langle\uparrow^{L}_{D}a,b\rangle_{L}\qquad\mathrm{for}\qquad a=\sum_{\alpha\in\Delta_{D}}a_{\alpha}\mathfrak{e}_{\alpha}\quad\mathrm{and}\quad b=\sum_{\delta\in\Delta_{L}}b_{\delta}\mathfrak{e}_{\delta},
\end{equation}
and observing that when $n=2$ and $L$ has Witt rank 2 the lattice $D$ is trivial and $\Theta_{D}$ from \eqref{ThetaD} is just the constant function 1, Theorem 4.4.1 of \cite{[BHY]} evaluates the theta lift $\Phi_{m,p_{D}\mu}^{D}$ as follows.
\begin{prop} \label{PhiDeval}
The number $\Phi_{m,p_{D}\mu}^{D}$ vanishes when $\mu \not\in J^{\perp}_{L^{*}}/J^{\perp}_{L}$, and if $\mu \in J^{\perp}_{L^{*}}/J^{\perp}_{L}$ then we have \[\Phi_{m,p_{D}\mu}^{D}=\begin{cases}\frac{8\pi}{n-2}\operatorname{CT}\Big(\big\langle q\frac{d}{dq}\uparrow^{L}_{D}(\Theta_{D}),F_{m,\mu}^{+}\big\rangle_{L}\Big) & n\geq3 \\[1ex] \frac{\pi}{3}\operatorname{CT}\big(E_{2}\langle\uparrow^{L}_{D}(1),F_{m,\mu}^{+}\rangle_{L}\big) & n=2.\end{cases}\]
\end{prop}
We remark that the value of $\Phi_{m,p_{D}\mu}^{D}$ for $n\geq3$ in Proposition \ref{PhiDeval} also equals, by modifying the proof of Theorem 9.2 of \cite{[Bo1]} as in Proposition \ref{PhiKeval} below, to the expression
\begin{equation} \label{PetPhiD} \frac{\pi}{3}\operatorname{CT}\big(E_{2}\cdot\langle\uparrow^{L}_{D}(\Theta_{D}),F_{m,\mu}^{+}\rangle_{L}\big)-\frac{\pi}{3}\big(E_{2}^{*}\uparrow^{L}_{D}(\Theta_{D}),\xi_{1-\frac{n}{2}}F_{m,\mu}\big)_{\mathrm{Pet}}, \end{equation}
where the latter pairing is the Petersson inner product of the nearly holomorphic modular form $E_{2}^{*}\uparrow^{L}_{D}(\Theta_{D})$ and the cusp form $\xi_{1-\frac{n}{2}}F_{m,\mu}$, both of weight $1+\frac{n}{2}$ and representation $\rho_{L}$. This gives back the expression for $n=2$ in Proposition \ref{PhiDeval}, since as in the proof of Theorem 4.4.1 of \cite{[BHY]}, the Petersson inner product there vanishes when $\Theta_{D}$ reduces to a constant. In addition, since $\big\langle q\frac{d}{dq}\uparrow^{L}_{D}(\Theta_{D})$ is holomorphic in $q$ with integral coefficients and without a constant term, the value of the multiplicity $\operatorname{mult}_{J}(m,\mu)=\frac{1}{8\pi}\Phi_{m,p_{D}\mu}^{D}$ from Definition \ref{mults}, with $n\geq3$, is rational. More explicitly, the symmetry of $\uparrow^{L}_{D}(\Theta_{D})$ and the simple principal part of $F_{m,\mu}^{+}$ reduce its value to
\begin{equation} \label{eltsmmu}
\frac{2m}{n-2}\big|\{\beta \in p_{D}\mu|\;\beta^{2}=2m\}\big|.
\end{equation}
This is not necessarily true for $n=2$, where the latter expression vanishes by the triviality of $D$, since when $L$ has Witt rank 2 multiples of $\uparrow^{L}_{D}(1)$ are $\overline{\rho}_{L}$-invariant vectors by which modifying $F_{m,\mu}^{+}$ can affect the value from Proposition \ref{PhiDeval}.

\smallskip

We now turn to the evaluation of $\Phi_{m,p_{K}\mu}^{K}\big(\frac{\omega}{|\omega|}\big)$ for a primitive vector $\omega \in K \cap C$, for which $\omega^{2}=-2N_{\omega}$ as usual, and we write $N$ for $N_{\omega}$ for short. Recall that this expression is defined in \eqref{liftK}, where the theta function $\Theta_{K}$ is defined as in \eqref{Thetadef} with the negative definite subspace of $K_{\mathbb{R}}$ that is spanned by $\omega$. Denote the orthogonal complement of $\omega$ in $K$ by $\omega^{\perp}_{K}$ and in $K_{\mathbb{R}}$ by just $\omega^{\perp}$, and define the $\rho_{K}\otimes\rho_{N}$-valued function
\begin{equation} \label{Thetaomegadef}
\Theta_{K,\omega}(\tau)=\sum_{r\in\mathbb{Z}/2N\mathbb{Z}}\sum_{\substack{\lambda \in K^{*}/\mathbb{Z}\omega \\ (\lambda,\omega)\equiv-r\pmod{2N}}}\mathbf{e}\Big(\tau\tfrac{\lambda_{\omega^{\perp}}^{2}}{2}\Big)\mathfrak{e}_{\lambda+K}\otimes\mathfrak{e}_{r}.
\end{equation}
Let $\Lambda_{N}$ be the 1-dimensional lattice spanned by a vector of norm $2N$, so that its discriminant $\Delta_{N}=\Delta_{\Lambda_{N}}$ is cyclic of order $2N$, and we denote by $\rho_{N}$ the Weil representation $\rho_{\Lambda_{N}}$ and by $\Theta_{N}$ the theta function $\Theta_{\Lambda_{N}} \in M_{\frac{1}{2}}(\rho_{N})$, which is defined by
\begin{equation} \label{ThetaNdef}
\Theta_{\Lambda_{N}}(\tau)=\sum_{l=-\infty}^{\infty}q^{l^{2}/4N}\mathfrak{e}_{l+2N\mathbb{Z}}=\sum_{r\in\mathbb{Z}/2N\mathbb{Z}}\bigg(\sum_{l \equiv r\pod{4N}}q^{l^{2}/2N}\bigg)\mathfrak{e}_{r}.
\end{equation}
If $\overline{\Theta}_{N}$ is the complex conjugate of the function from \eqref{ThetaNdef}, and $\langle\cdot,\cdot\rangle_{N}$ is the pairing associated with the lattice $\Lambda_{N}$, then we obtain the following decomposition.
\begin{lem} \label{thetadecom}
The function $\Theta_{K,\omega}$ from \eqref{Thetaomegadef} is modular of weight $\frac{n-1}{2}$ and representation $\rho_{K}\otimes\rho_{N}$ with respect to the full group $\operatorname{Mp}_{2}(\mathbb{Z})$, and we have the equality \[\Theta_{K}\big(\tau,\tfrac{\omega}{|\omega|}\big)=\sqrt{y}\big\langle\Theta_{K,\omega}(\tau),\overline{\Theta}_{N}(\tau)\big\rangle_{N}.\]
\end{lem}
Lemma \ref{thetadecom} is well-known. For a proof of both assertions in a more general setting see Theorem 1.11 and Proposition 1.8 of \cite{[Ze4]}.

\smallskip

Expressions like $\Phi_{m,p_{K}\mu}^{K}\big(\frac{\omega}{|\omega|}\big)$ are evaluated, in a slightly different context, in Theorem~10.6 of \cite{[Bo1]}, which uses the incorrect Lemma~9.5 of that reference. Explicit pre-images of $-\frac{\sqrt{N}}{8\pi}\Theta_{N}$ under the operator $\xi_{3/2}$ are constructed in Theorem~4.3 of \cite{[BS]} as follows.
\begin{thm} \label{xiGtheta}
For every natural $N$ there exists a harmonic Maass form $G_{N}$ of weight $\frac{3}{2}$ and representation $\overline{\rho}_{N}$, with $\xi_{3/2}G_{N}=-\frac{\sqrt{N}}{8\pi}\Theta_{N}$, and such that its holomorphic part $G_{N}^{+}$ has rational Fourier coefficients with bounded denominators.
\end{thm}
It is clear that $G_{N}$ is unique only up to adding elements of $M_{3/2}^{!}(\overline{\rho}_{N})$ with rational Fourier coefficients, but our results below will be independent of the choice of $G_{N}$. The evaluation of $\Phi_{m,p_{K}\mu}^{K}\big(\frac{\omega}{|\omega|}\big)$ can now be carried out as follows.
\begin{prop} \label{PhiKeval}
The value of $\Phi_{m,p_{K}\mu}^{K}\big(\frac{\omega}{|\omega|}\big)$ is \[-\frac{8\pi}{\sqrt{N}}\operatorname{CT}\big(\big\langle\langle\uparrow^{L\oplus\Lambda_{N}}_{K\oplus\Lambda_{N}}(\Theta_{K,\omega}),G_{N}^{+}\rangle_{N},F_{m,\mu}^{+}\big\rangle_{L}\big)+
\frac{8\pi}{\sqrt{N}}\big(\big\langle\uparrow^{L\oplus\Lambda_{N}}_{K\oplus\Lambda_{N}}(\Theta_{K,\omega}),G_{N}\big\rangle_{N},\xi_{1-\frac{n}{2}}F_{m,\mu}\big)_{\mathrm{Pet}}^{\mathrm{reg}},\] where the second term is the regularized Petersson inner product of these modular forms.
\end{prop}

\begin{proof}
We follow the proof of Theorem 10.6 of \cite{[Bo1]}. As in \eqref{liftK}, we need the regularized integral of $\big\langle\Theta_{K}(\tau,\mathbb{R}\omega),\downarrow^{L}_{K}\big(F_{m,\mu}(\tau)\big)\big\rangle_{K}d\mu(\tau)$, and Lemma \ref{thetadecom} presents this theta function as $\sqrt{y}$ times the pairing $\big\langle\Theta_{K,\omega}(\tau),\overline{\Theta}_{N}(\tau)\big\rangle_{N}$ of the holomorphic theta function $\Theta_{K,\omega}$ and the anti-holomorphic theta function $\overline{\Theta}_{N}$. Via \eqref{arrowpar} we can express this measure as
\[\big\langle\uparrow^{L}_{K}\big(\langle\Theta_{K,\omega}(\tau),\overline{\Theta}_{N}(\tau)\rangle_{N}\big),F_{m,\mu}(\tau)\big\rangle_{L}=
\Big\langle\big\langle\uparrow^{L\oplus\Lambda_{N}}_{K\oplus\Lambda_{N}}\big(\Theta_{K,\omega}(\tau)\big),\overline{\Theta}_{N}(\tau)\big\rangle_{N},F_{m,\mu}(\tau)\Big\rangle_{L}\] times $y^{-3/2}dxdy$. Theorem \ref{xiGtheta} allows us to write the non-holomorphic multiplier $\overline{\Theta}_{N}y^{-3/2}$ as $\frac{16\pi i}{\sqrt{N}}\partial_{\overline{\tau}}G_{N}$, and if we assume for the moment that $F_{m,\mu}$ is holomorphic on $\mathcal{H}$ then the differential 2-form that we have in the integrand is the exact $\operatorname{Mp}_{2}(\mathbb{Z})$-invariant form
\begin{equation} \label{exact}
d\bigg(\frac{8\pi}{\sqrt{N}}\Big\langle\big\langle\uparrow^{L\oplus\Lambda_{N}}_{K\oplus\Lambda_{N}}\big(\Theta_{K,\omega}(\tau)\big),G_{N}(\tau)\big\rangle_{N},F_{m,\mu}(\tau)\Big\rangle_{L}d\tau\bigg).
\end{equation}
Now, the integral from \eqref{liftK} can be carried out over the fundamental domain \[\mathcal{F}=\big\{\tau=x+iy\in\mathcal{H}\big|\;|x|\leq\tfrac{1}{2},\ |\tau|\geq1\big\}\] of $\operatorname{Mp}_{2}(\mathbb{Z})$, and in the regularization of this divergent integral, we evaluate the integral on $\mathcal{F}_{R}=\{\tau\in\mathcal{F}|\;y \leq R\}$ for large $R$, and take the limit $R\to\infty$. We can thus apply Stokes' Theorem for the integral over $\mathcal{F}_{R}$, and the $\operatorname{Mp}_{2}(\mathbb{Z})$-invariance of the 1-form appearing as the argument of $d$ on the right hand side of \eqref{exact} implies that on the integrals on the boundary parts of $\mathcal{F}_{R}$ that are also part of the boundary of $\mathcal{F}$ cancel. All that remains is the integral along $\big[-\frac{1}{2},\frac{1}{2}\big]+iR$ oriented negatively, which produces the constant term of the Fourier expansion of the differential form in question, at $y=R$. The part arising from $G_{N}^{+}$ and $F_{m,\mu}^{+}$ yields the first asserted term, and since the Fourier coefficients of $G_{N}^{-}$ and $F_{m,\mu}^{-}$ multiply $q^{l}$ for negative $l$ by decreasing functions of $y$, the contribution from this part vanishes at the limit $R\to\infty$.

However, as $F_{m,\mu}$ is not holomorphic, the formula from \eqref{exact} has to be corrected by subtracting the 2-form \[\frac{16\pi i}{\sqrt{N}}\Big\langle\big\langle\uparrow^{L\oplus\Lambda_{N}}_{K\oplus\Lambda_{N}}\big(\Theta_{K,\omega}(\tau)\big),G_{N}(\tau)\big\rangle_{N},\partial_{\overline{\tau}}F_{m,\mu}(\tau)\Big\rangle_{L}dxdy,\] and the relation between $\xi_{1-\frac{n}{2}}$ and $\partial_{\overline{\tau}}$ allows us to write \[\partial_{\overline{\tau}}F_{m,\mu}=\frac{i}{2}y^{\frac{n}{2}-1}\overline{\xi_{1-\frac{n}{2}}F_{m,\mu}}.\] Hence the additional term yields the asserted regularized Petersson inner product. This completes the proof of the proposition.
\end{proof}

We remark that if $n=3$ then in Example \ref{Siegel} the theta function $\Theta_{D}$ is equal to $\Theta_{\Lambda_{1}}$, the maps $p_{K}$ and $p_{D}$ are trivial, and then the multiplicity $\operatorname{mult}_{J}(m,\mu)$ with which the 1-dimensional toroidal boundary component appears in $Z(m,\mu)$ is $4m$ when $m$ is a square and 0 otherwise (this is most easily visible from \eqref{eltsmmu}). In the second compactification appearing there, the orbit of inner rays yields another boundary divisor, for which $N_{\omega}=3$ and we can take $G_{3}$ to be the function formed from Zagier's Eisenstein series, with $G_{3}^{+}$ based on Hurwitz class numbers (Lemma 9.5 of \cite{[Bo1]} does hold for prime $N$). Moreover, since the space $S_{3/2}(\rho_{L})$ is trivial, the evaluation of the corresponding multiplicity $\operatorname{mult}_{I,\omega}(m,\mu)$ from Definition \ref{mults} via Proposition \ref{PhiKeval} gives $-\frac{1}{12}$ times the constant term of $F_{m,\mu}^{+}=F_{m,\mu}$ plus $\sum_{\lambda \in K^{*}/\mathbb{Z}\omega}H(12m-6\lambda_{\omega^{\perp}}^{2})$, which is therefore rational.

\section{Modularity of Heegner divisors \label{Modularity}}

Here we consider some basic parts from the theory of Borcherds products, and show how to modify them in order to prove our main result, Theorem \ref{Ztormod}, which is the modularity of the appropriately modified divisors from Theorem \ref{GreenPhi} in the Chow group of $X_{L}^{\mathrm{tor}}$.

\smallskip

Consider again the Grassmannian $\operatorname{Gr}(L)$ that is associated via \eqref{Grassdef} to the even lattice $L$. Recall that the map from \eqref{fiboverGVR} describes the chosen component $P^{+}$ of the set $P$ from \eqref{bunGV} as a $\mathbb{C}^{\times}$-bundle over $\operatorname{Gr}(L)$, and we define an \emph{automorphic form of weight $m$ on $\operatorname{Gr}(L)$ with respect to $\Gamma_{L}$} to be a function from $P^{+}$ to $\mathbb{C}$ that is $\Gamma_{L}$-invariant and homogenous of degree $-m$. Given an oriented sublattice $I$ of rank 1 in $L$, one can identify $\operatorname{Gr}(L)$ with the affine tube domain $K_{\mathbb{R}}^{1}+iC$ from \eqref{GVRKRiC}, and choosing a complement for $I^{\perp}_{L^{*}}$ in $L^{*}$ maps this affine tube domain onto the linear one $K_{\mathbb{R}}+iC \subseteq K_{\mathbb{C}}$. Combining these maps with the section from $\operatorname{Gr}(L)$ to $P^{+}$, one identifies an automorphic form of weight $m$ on $\operatorname{Gr}(L)$ with a function on $K_{\mathbb{R}}+iC \subseteq K_{\mathbb{C}}$ which transforms under $\Gamma_{L}$ with the $m$th power of an appropriate factor of automorphy. In any case, we are interested here only in such automorphic forms that are meromorphic.

An important class of meromorphic automorphic forms on $X_{L}$ is given by \emph{Borcherds products}, see Theorem~13.3 of \cite{[Bo1]}. The result is as follows.
\begin{thm} \label{Bor13.3}
Let $F \in M_{1-n/2}^{!}(\overline{\rho}_{L})$ be given, with a Fourier expansion as in the holomorphic part of \eqref{Fourier}, and assume that $c^{+}(\delta,l)\in\mathbb{Z}$ if $l<0$. Then there exists a meromorphic automorphic form $\Psi$ on $\operatorname{Gr}(L)$ of weight $c^{+}(0,0)/2$ and some multiplier system $\chi$ of finite order with respect to $\Gamma_{L}$, whose divisor of $X_{L}$ is given by
\[\operatorname{div}_{X_{L}}(\Psi)=\frac{1}{2}\sum_{\mu\in\Delta_{L}}\sum_{0<m\in\mu^{2}/2+\mathbb{Z}}c^{+}(\mu,-m)Z(m,\mu).\]
\end{thm}

The proof of Theorem \ref{Bor13.3} is based on the fact that the theta lift of $F$, given in \eqref{liftdef} with $F_{m,\mu}$ replaced by $F$, is, up to an additive constant, 2 times the logarithm of the Petersson metric of a meromorphic function $\Psi$ on $\operatorname{Gr}(L)$ with the asserted divisor, which is then shown to be automorphic. By the theory of Baily--Borel, the power of $\Psi$ for which the multiplier system $\chi$ becomes trivial is defined on the Baily--Borel compactification $X_{L}^{\mathrm{BB}}$, and hence its pullback is defined also on $X_{L}^{\mathrm{tor}}$. We can thus determine its divisor on $X_{L}^{\mathrm{tor}}$.

\begin{thm} \label{divbctor}
Assume that $n\geq3$ or that $n=2$ and $\mathbb{C}[\Delta_{L}]$ contains no $\overline{\rho}_{L}$-invariant vectors, take $F \in M_{1-\frac{n}{2}}^{!}(\overline{\rho}_{L})$ that satisfies the condition of Theorem \ref{Bor13.3}, and let $\Psi$ be the associated Borcherds product. If the character $\chi$ of $\Psi$ is trivial, then the divisor of $\Psi$ on $X_{L}^{\mathrm{tor}}$ is given by
\[\operatorname{div}_{X_{L}^{\mathrm{tor}}}(\Psi)=\frac{1}{2}\sum_{\mu\in\Delta_{L}}\sum_{0<m\in\mu^{2}/2+\mathbb{Z}}c^{+}(\mu,-m)Z^{\mathrm{tor}}(m,\mu).\] The difference between this divisor and the closure of $\operatorname{div}_{X_{L}}(\Psi)$ from Theorem \ref{Bor13.3} is \[\frac{1}{24}\sum_{J}\operatorname{CT}\big(E_{2}\cdot\langle\uparrow^{L}_{D}(\Theta_{D}),F\rangle_{L}\big) \cdot B_{J}-\sum_{I}\sum_{\mathbb{R}_{+}\omega}\operatorname{CT}\Big(\big\langle\langle\uparrow^{L\oplus\Lambda_{N}}_{K\oplus\Lambda_{N}}(\Theta_{K,\omega}),G_{N}^{+}\rangle_{N},F\big\rangle_{L}\Big) \cdot B_{I,\omega},\] where $J$, $D$, $I$, $K$, $\omega$, and $N$ are as in Theorem \ref{GreenPhi}, and the coefficients are rational.
\end{thm}

\begin{proof}
Recalling that our assumptions on $n$ and $\rho_{L}$ imply that the harmonic Maass form of weight $1-\frac{n}{2}$ and representation $\overline{\rho}_{L}$ with a given principal part is unique, we deduce that if $F \in H_{1-\frac{n}{2}}(\overline{\rho}_{L})$ has the Fourier expansion \eqref{Fourier} and $n\geq3$ then
\begin{equation} \label{Fprin}
F=\frac{1}{2}\sum_{\mu\in\Delta_{L}}\sum_{0<m\in\frac{\mu^{2}}{2}+\mathbb{Z}}c^{+}(\mu,-m)F_{m,\mu}.
\end{equation}
The linearity of the regularized theta lift from \eqref{liftdef} thus implies that the Petersson metric of $\Psi$ is given by \[-\log\|\Psi\|^{2}=\frac{1}{2}\sum_{\mu\in\Delta_{L}}\sum_{0<m\in\frac{\mu^{2}}{2}+\mathbb{Z}}c^{+}(\mu,-m)\Phi_{m,\mu}^{L}.\] Consequently, the divisor of $F$ on $X_{L}^{\mathrm{tor}}$ is given by the corresponding linear combination of the divisors $Z^{\mathrm{tor}}(m,\mu)$ of the Green functions. Hence the assertion about the divisor of $\Psi$ follows from Theorem~\ref{GreenPhi}.

The difference between the divisors on $X_{L}^{\mathrm{tor}}$ and on $X_{L}$ is the sum of the parts of the divisors $Z^{\mathrm{tor}}(m,\mu)$ that are supported on the boundary divisors. The multiplicity of each such boundary divisor is the linear combination of the corresponding multiplicities $\operatorname{mult}_{J}(m,\mu)$ or $\operatorname{mult}_{I,\omega}(m,\mu)$ from Definition \ref{mults}, according to the coefficients from \eqref{Fprin}. But the formulas for these multiplicities, given in \eqref{PetPhiD} and in Proposition \ref{PhiKeval} respectively, are linear in $F$, and we can thus just substitute $F$ instead of $F_{m,\mu}$ there. But since $F$ is weakly holomorphic, we have $F^{+}=F$ and the terms involving $\xi_{1-\frac{n}{2}}F$ vanish, and we remain with the asserted multiplicities. Finally, as the principal part of $F$ determines $F$ and has rational Fourier coefficients, it is therefore invariant under $\operatorname{Aut}(\mathbb{C}/\mathbb{Q})$, and therefore so is $F$ itself and all of its coefficients are rational (see \cite{[McG]}). Since Theorem \ref{xiGtheta} provides the same property for $G_{N}$, and it is clear that $E_{2}$ and the theta functions share this property, it follows that the asserted multiplicities are rational as desired. This proves the theorem.
\end{proof}
When $n\geq3$ we could have also written the multiplicity of $B_{J}$ in Theorem \ref{divbctor} as $\frac{1}{n-2}\operatorname{CT}\Big(\big\langle q\frac{d}{dq}\uparrow^{L}_{D}(\Theta_{D}),F\big\rangle_{L}\Big)$. We remark that the multiplicity of $B_{I,\omega}$ does not depend on the choice of the function $G_{N}$ from Theorem \ref{xiGtheta}, since this function pairs with the weakly holomorphic modular form $\langle\Theta_{K,\omega},\downarrow^{L}_{K}(F)\rangle_{K}$ of weight $\frac{1}{2}$, and Proposition 3.5 of \cite{[BFu]} shows that this constant term can be evaluated only using $\xi_{3/2}G_{N}$. We also remark that a bit of additional analysis shows that the results of Theorem \ref{divbctor} hold also without the assumption on $\overline{\rho}_{L}$ when $n=2$, and even when $n=1$, provided that $F$ is taken to have rational Fourier coefficients for the rationality of the multiplicities to remain valid.

\begin{rmk} \label{pfBo1K3}
Theorem \ref{divbctor} yields new proofs of some known results from the literature. First, for describing the behavior of $\Psi$ near a 0-dimensional cusp $\Xi=\Xi_{I}$, we express the corresponding lift $\Phi$ from \eqref{liftdef} as in Proposition \ref{PhimmuL}, where the parts involving the functions $\mathcal{V}_{1+\frac{n}{2}}$ vanish and all the terms of the form $\log\big|1-\mathbf{e}\big((\nu,Z_{L})\big)\big|$ give the logarithm of the expansion of $\Psi$ as a Borcherds product. The first term $\Phi^{K}$, the theta lift of $\downarrow^{L}_{K}(F)$ with respect to the lattice $K$, is a continuous piecewise linear function of $Y$, which is given, in a Weyl chamber $W$ associated with $F$, by $8\sqrt{2}\pi$ times the pairing of $Y$ with a vector $\rho \in K_{\mathbb{Q}}$ that is called the \emph{Weyl vector} associated with $L$, $I$, $F$, and $W$. Hence the product expansion of $\Psi$ is multiplied by $\mathbf{e}\big((Z,\rho)\big)$ in the variable $Z=X+iY$ in the linear tube domain $K_{\mathbb{R}}+iC$ as well as on the corresponding subset of the torus $T_{K}$, and if $\rho \in K^{*}$ then this function is just the toric function $z_{\rho}$. The order of $\Psi$ along $B_{I,\omega}$ for primitive $\omega \in K \cap C$ with $\mathbb{R}_{+}\omega\in\Sigma(\Xi)$ is thus the order of $z_{\rho}$, so that Proposition \ref{ordTdiv} reproduces the corrected result of Theorem 10.6 of \cite{[Bo1]}. Similarly, by considering the order of $\Psi$ along $B_{J}$, Equation (4.27) of \cite{[K3]} shows that we have reproduced also the parts of Theorem 2.1 and Corollary 2.3 of that reference involving the power of the coordinate $q_{2}$ defining this divisor. We also remark that combining $\sum_{l=0}^{\infty}\log|1-\alpha q^{l}|$ from the proof of Proposition \ref{PhinearJ} with the summands arising from $\mathbb{Z}w-\nu$ yields, by the Jacobi triple product identity, the logarithm of the absolute value of the quotients $\theta/\eta$ of the translated classical theta function over the Dedekind eta function appearing in the expansions from \cite{[K3]}.
\end{rmk}

\smallskip

We now recall Borcherds' modularity criterion, see Theorem 3.1 of \cite{[Bo2]}. The rationality statement appearing here is a consequence of the main result of \cite{[McG]}.
\begin{prop} \label{Serredual}
Let $a(\mu,m)$ for $\mu\in\Delta_{L}$ and $0 \leq m\in\frac{\mu^{2}}{2}+\mathbb{Z}$ be elements of some complex vector space $Y$, and let \[\sum_{\mu\in\Delta_{L}}\sum_{0 \leq m\in\frac{\mu^{2}}{2}+\mathbb{Z}}a(\mu,m) \cdot q^{m}\mathfrak{e}_{\mu}\] be the resulting formal $\rho_{L}(T)$-invariant power series with values in $Y\otimes_{\mathbb{C}}\mathbb{C}[\Delta_{L}]$. Assume that the coefficients $a(\mu,m)$ have the symmetry property that is required for this series to be invariant under the action of the center of $\operatorname{Mp}_{2}(\mathbb{Z})$ with weight $2-k$ and representation $\rho_{L}$. Then this power series belongs to $Y\otimes_{\mathbb{C}}M_{2-k}(\rho_{L})$ if and only if its pairing with every $F \in M_{k}^{!}(\overline{\rho}_{L})$ vanishes, that is, denoting the Fourier expansion by $F(\tau)=\sum_{\delta,l}c(\delta,l)q^{l}\mathfrak{e}_{\delta}$, we have the equality \[\sum_{\mu\in\Delta_{L}}\sum_{l\in\frac{\mu^{2}}{2}+\mathbb{Z}}c(\mu,-l)a(\mu,l)=0 \in Y.\] In addition, having these relations only for those elements $F \in M_{k}^{!}(\overline{\rho}_{L})$ whose Fourier coefficients are integral suffices for the series to be in $Y\otimes_{\mathbb{C}}M_{2-k}(\rho_{L})$.
\end{prop}
Our second main result, Theorem \ref{main} from the introduction, is the following one.
\begin{thm} \label{Ztormod}
Assume that $n$ is larger than the Witt rank of $L$, set $Z^{\mathrm{tor}}(0,0)$ to be any divisor on $X_{L}^{\mathrm{tor}}$ that represents the line bundle of modular forms of weight $-1/2$, and write $[Z^{\mathrm{tor}}(m,\mu)]$ for the image of $Z^{\mathrm{tor}}(m,\mu)$ in $\operatorname{CH}^{1}(X_{L}^{\mathrm{tor}})\otimes\mathbb{R}$ for every pair $(m,\mu)$. Then the formal power series \[\sum_{\mu\in\Delta_{L}}\sum_{m\in\frac{\mu^{2}}{2}+\mathbb{Z}}[Z^{\mathrm{tor}}(m,\mu)]\cdot q^{m}\mathfrak{e}_{\mu}\] is a modular form of weight $1+n/2$ and representation $\rho_{L}$ with coefficients in $\operatorname{CH}^{1}(X_{L}^{\mathrm{tor}})\otimes\mathbb{R}$.
\end{thm}

\begin{proof}
If $F \in M_{1-\frac{n}{2}}^{!}(\overline{\rho}_{L})$ has the Fourier expansion from \eqref{Fourier}, and if $c_{\mu}(n)\in\mathbb{Z}$ wherever $n<0$, then by Theorem~\ref{divbctor} we get the identity
\[\sum_{\mu\in\Delta_{L}}\sum_{m\in\frac{\mu^{2}}{2}+\mathbb{Z}}c_{\mu}(-m)Z^{\mathrm{tor}}(m,\mu)=0\] in $\operatorname{CH}^{1}(X_{L}^{\mathrm{tor}})\otimes\mathbb{R}$. The result thus follows from Proposition \ref{Serredual} as in \cite{[Bo2]}, with the rationality statement allowing one to omit the Galois action appearing in that reference. This proves the theorem.
\end{proof}

\begin{rmk} \label{modCH1}
According to Theorem \ref{GreenPhi}, the pairs $(Z^\mathrm{tor}(m,\mu), \frac{1}{2}\Phi_{m,\mu}^{L})$ define classes in the arithmetic Chow group $\widehat{\operatorname{CH}}^{1}(X_{L}^{\mathrm{tor}},\mathcal{D}_{\mathrm{pre}})$ of $X_{L}^{\mathrm{tor}}$ over $\mathbb{C}$ in the sense of Definition~1.15 of \cite{[BBK]}. The same argument as in the proof of Theorem \ref{Ztormod} shows that the generating series of these classes is modular with values in that arithmetic Chow group.
\end{rmk}

\end{document}